\newif\ifsiam

\siamfalse

\newif\ifarxiv
\ifsiam \arxivfalse \else \arxivtrue \fi

\ifsiam
\documentclass[review,onefignum,onetabnum]{siamonline250211}
\fi

\ifarxiv
\documentclass[english,final,twoside,pdftex]{myarticle}
\fi

\usepackage{lipsum}
\usepackage{amsfonts}
\usepackage{graphicx}
\usepackage{epstopdf}
\usepackage{algorithmic}
\usepackage{xr}
\ifpdf
  \DeclareGraphicsExtensions{.eps,.pdf,.png,.jpg}
\else
  \DeclareGraphicsExtensions{.eps}
\fi

\usepackage{enumitem}
\usepackage{tikz}
\usepackage{pgfplots}
\usepgfplotslibrary{groupplots}
\setlist[enumerate]{leftmargin=.5in}
\setlist[itemize]{leftmargin=.5in}

\definecolor{royalblue}{HTML}{4169E1}
\definecolor{crimson}{HTML}{DC143C}
\definecolor{forestgreen}{HTML}{228B22}
\definecolor{darkorange}{HTML}{FF8C00} % For the Hopf branches
\definecolor{gold}{HTML}{FFD700}       % For the other Hopf branch

\pgfplotsset{
    compat=1.18,
    scatter style/.style={only marks, mark size=1.2pt},
    bifurcation line style/.style={very thick, no marks},
    fit line style/.style={thick, no marks}
}

\usetikzlibrary{external} 
\tikzsetexternalprefix{Figures/} 
\tikzset{external/system call={lualatex -shell-escape -halt-on-error -interaction=batchmode -jobname "\image" "\texsource"}}
\ifsiam
    \newsiamremark{remark}{Remark}
    \newsiamremark{hypothesis}{Hypothesis}
    \crefname{hypothesis}{Hypothesis}{Hypotheses}
    \newsiamthm{claim}{Claim}
    \newsiamthm{prop}{Proposition}
\else
    \usepackage{amsthm}
    \usepackage{cleveref}
    
    \crefname{hypothesis}{Hypothesis}{Hypotheses}
   
   \newtheorem{prop}{Proposition}
\fi

\newcommand{\TheTitle}{Bifurcation Analysis of Generalized Hopf Bifurcation in Ordinary and Delay Differential Equations} 
\newcommand{\TheAuthors}{N.A.M. Delmeire, M.M. Bosschaert, Yu. A. Kuznetsov} 

\ifsiam
    \headers{Bifurcation Analysis of Generalized Hopf Bifurcation in ODEs and DDEs}{N.A.M. Delmeire, M.M. Bosschaert, Yu. A. Kuznetsov}
  \else
  \usepackage{fancyhdr}
\fi

\ifsiam \title{{\TheTitle}\thanks{Submitted to the editors DATE.}} \fi
\ifarxiv \title{\TheTitle} \fi
\author{Nathalie A. M. Delmeire\thanks{Department of Mathematics, Utrecht University, 3508 TA Utrecht, The Netherlands
  .}
\and Maikel M. Bosschaert\thanks{Data Science Institute, Hasselt University, Diepenbeek Campus, 3590 Diepenbeek, Belgium \email{(maikel.bosschaert@uhasselt.be)}.}
\and Yuri A. Kuznetsov\thanks{Department of Mathematics, Utrecht University, 3508 TA Utrecht, The Netherlands and Department of Applied
Mathematics, University of Twente, 7500 AE Enschede, The Netherlands \email{(i.a.kouznetsov@uu.nl)}.}
}

\usepackage{listings}

\usepackage{amsopn}

\usepackage{mathtools}
\newcommand{\N}{\mathbb{N}} % de natuurlijke getallen
\newcommand{\R}{\mathbb{R}} % de reele getallen
\newcommand{\C}{\mathbb{C}}
\allowdisplaybreaks
\graphicspath{ {Figures/} }
\ifpdf
\hypersetup{
  pdftitle={Bifurcation Analysis of Generalized Hopf Bifurcation in Ordinary and Delay Differential Equations},
  pdfauthor={N. A. M. Delmeire, M. M. Bosschaert, Yuri A. Kuznetsov}
}
\fi

\begin{document}

\maketitle

% REQUIRED
\begin{abstract}
The generalized Hopf (Bautin) bifurcation is a well-studied codimension two bifurcation characterized by an equilibrium with a pair of simple purely imaginary eigenvalues as the only critical eigenvalues and the vanishing first Lyapunov coefficient. This bifurcation arises in both ordinary differential equations (ODEs) and delay differential equations (DDEs). Generically, a codimension one bifurcation curve of nonhyperbolic (double) limit cycles ($LPC$ curve) emanates from a generalized Hopf point at the Hopf bifurcation curve $H$. By performing the parameter-dependent center manifold reduction near this point, the first-order predictors to initiate continuation of the $LPC$ curve have been derived in the literature. These predictors, however, do not distinguish the curves $H$ and $LPC$ in the parameter space. In this paper, we overcome this deficiency by deriving higher-order predictors for the $LPC$ curve in ODEs and DDEs for the first time. The new predictors, which require seventh-order derivatives of the system, have been implemented, and their effectiveness is demonstrated in several models.
\end{abstract}

% REQUIRED

% REQUIRED

\section{Introduction}
\begin{figure}[!ht]
    \centering
    \includegraphics[scale=0.65]{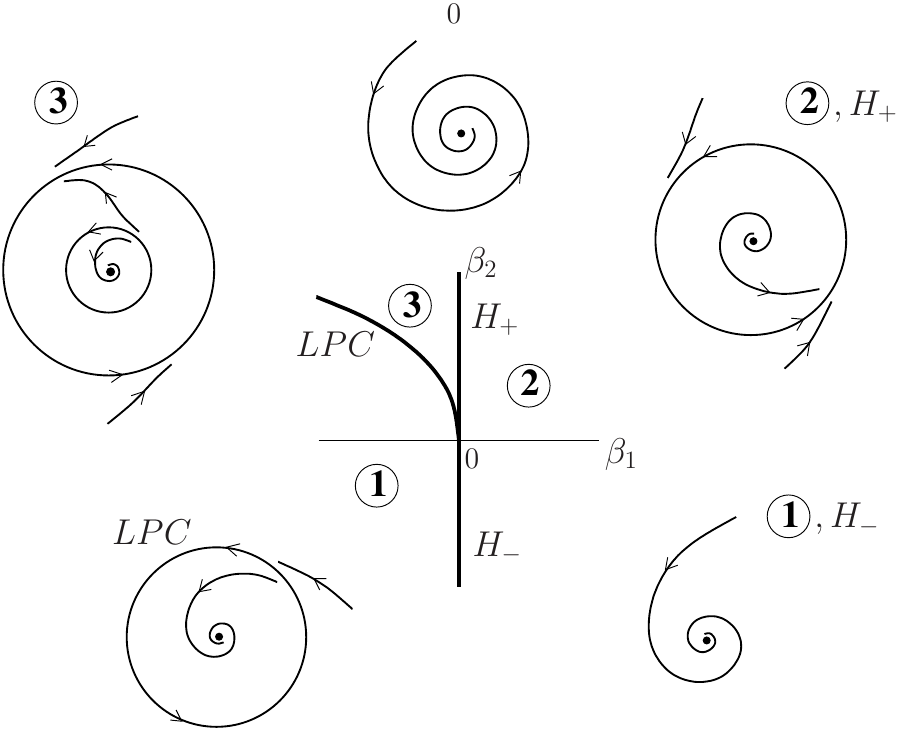}
    \caption{Bifurcation diagram near a generalized Hopf bifurcation for the case $l_2 < 0$. Figure adapted from \emph{\cite{kuznetsov2023elements}}.}
    \label{fig:GH_bifurcationdiagram}
\end{figure}

Many applications use models that consist of autonomous Ordinary Differential Equations (ODEs)
\begin{equation}\label{eq:1}
    \dot{x}(t) = F(x(t),\alpha), 
\end{equation}
where $x: \R \to \R^n$, $\alpha \in \R^p$ and $F:\R^n\times \R^p \to \R^n$ is a smooth mapping. %We use dot notation to indicate the derivative with respect to $t$, denoted as $\dot{x} = \frac{dx}{dt}$. 
One is often interested in the behaviour of such systems under parameter variations. If the system has an equilibrium point, the variation of one parameter might result in a local codim 1 bifurcation. A well known example is the \textit{$($Andronov-$)$Hopf bifurcation}. By allowing the variation of a second parameter, the first Lyapunov coefficient may vanish while the two purely imaginary eigenvalues persist. At such a point a \textit{generalized Hopf bifurcation} (or \textit{Bautin bifurcation}) occurs. This bifurcation, which requires two conditions to be satisfied to manifest, is a \textit{codimension two} bifurcation. There are four other well-known local codim 2 bifurcations. For example, a Hopf-Hopf bifurcation occurs when an equilibrium has two pairs of simple purely imaginary eigenvalues. Their normal forms and local behaviour are discussed in detail in \cite[Chapter 8]{kuznetsov2023elements}. In this paper, we will only deal with the generalized Hopf bifurcation. A sketch of the bifurcation diagram for the case where $l_2 < 0$ is shown in \cref{fig:GH_bifurcationdiagram}. In this illustration, we can discern three different regions near the generalized Hopf bifurcation. In region {\bf 1} there is only a stable focus at the origin. If we move from region {\bf 1} to region {\bf 2}, passing the supercritical Hopf bifurcation line $H_{-}$, the equilibrium at the origin becomes unstable and a unique stable limit cycle appears. Continuing through the subcritical Hopf bifurcation line $H_{+}$ into region {\bf 3}, the equilibrium recovers its stability, alongside the emergence of an unstable cycle nested within the pre-existing stable cycle. These two cycles collide and disappear at the $LPC$ curve separating regions {\bf 3} and {\bf 1}, where a limit point (or fold) bifurcation of cycles occurs.
%(Sketch emanating codimension 1 curves and include illustration of the bifurcation diagram.)

In $n$-dimensional systems where a local bifurcation arises, a smooth family of parameter-dependent invariant \textit{center manifolds} $W_{\alpha}^c$ exists for sufficiently small $\|\alpha\|$. All qualitative changes near the bifurcation point occur on this lower dimensional manifold. For the Hopf and the generalized Hopf bifurcations, this manifold is two-dimensional. If we restrict our system \cref{eq:1} to the center manifold, the system can locally be transformed to a normal form. The existence of an invariant center manifold can be used to derive equations for the normal form coefficients.

Another interesting class of differential equations is that of Delay Differential Equations (DDEs). Unlike ODEs, which RHS depends only on current state information, the RHS of DDEs involves the values of the phase variables at previous times.  This accounts for delays in system’s response, capturing phenomena where the evolution of a system depends not only on its current state but also on its history. A specific case of DDEs that are often encountered in applications is of the following form
\begin{equation}\label{eq:Discrete_DDE}
    \dot{x}(t) = F(x(t), x(t - \tau_1), \ldots, x(t - \tau_m),\alpha), 
\end{equation}
where $x:\R \to \R^n$, $\alpha \in \R^p$, $F:\R^{n \times (m+1)} \times \R^p \to \R^n$ is a smooth mapping, and $0 < \tau_1 < \ldots < \tau_m =:h$ are constant delays. This is also known as a \textit{discrete DDE}. 
Such systems allow for study of bifurcations, including the generalized Hopf bifurcation (i.e. \cite{ion2007example,zhen2010bautin,ion2013bautin}). In contrast to ODEs, DDEs belong to the class of infinite dimensional dynamical systems, since a natural space of initial data (state space) here is $X=C([-h,0], \R^n)$. However, this state space does lead to some technical complications. This can be resolved with the help of perturbation theory for dual semigroups, which has been developed in \cite{clement1987perturbation,clement1988perturbation,clement1989perturbation,clement1989perturbation4} and is also known as \textit{sun-star calculus}. Using this framework, the existence of a finite-dimensional smooth invariant center manifold has been rigorously established for DDEs \cite{diekmann2012delay,bosschaert2020switching}. On this center manifold, the infinite-dimensional system of DDEs can be reduced to a finite-dimensional system of ODEs. This in turn allows us to ``lift'' results from the bifurcation theory of ODEs to the theory of DDEs. 

In most cases, the codim 1 bifurcation curves near a generalized Hopf bifurcation cannot be computed analytically. Instead one uses numerical methods for the location of the codim 2 bifurcation and the continuation of the emanating codim 1 bifurcation curves. If the approximate location of a generalized Hopf bifurcation point is known, we would like to switch to the $LPC$ curve emanating from this point using only local information available at the generalized Hopf point. To start the continuation, we need an approximation of the $LPC$ curve in the original parameter space, an approximation of the corresponding periodic orbit and the period. Combined, this forms a \textit{predictor}. A general method to derive such approximations, using only local information at the codim 2 bifurcation point, has been introduced in \cite{beyn2002numerical}. Within the framework of the sun-star calculus, the normalization method developed for ODEs has been extended to DDEs \cite{janssens2010normalization}. These methods have already been applied to derive first-order predictors for both ODEs \cite{kuznetsov2008switching} and DDEs \cite{bosschaert2020switching}. Such predictors, however, do not distinguish the curves $H$ and $LPC$ in the parameter space. 

In this paper, we derive a higher-order predictor for both ODEs and DDEs that does not suffer from this drawback. The paper is organized as follows. In \cref{sec:theory} we review the parameter-dependent center manifold reduction technique in the context of ODEs and DDEs, including a brief overview of the sun-star calculus needed for DDEs. The existing solution cases for the linear operator equations used in \cite{bosschaert2020switching} are insufficient when computing higher-order coefficients. Therefore, we have derived new, more general solutions which are presented in \cref{sec:DDE_eqsolve}. In \cref{sec:predictor_normalform}, we derive higher-order approximations of the $LPC$ curve for the normal form, concluding with a discussion on the coefficients required in the parameter and center manifold approximations. The method previously used for deriving parameter-dependent coefficients in \cite{kuznetsov2008switching} and \cite{bosschaert2020switching} does not extend to higher orders. Therefore, we used a slightly different approach, resulting in new equations for the parameter-dependent coefficients. These new equations are derived in \cref{sec:coeffs_ODE} for ODEs. In \cref{sec:coeff_DDE} we explain how these equations will change in the DDE setting. Despite many similarities between the equations of both classes of differential equations, one must still be careful. For instance, one term vanishes in the ODE case, but persists for DDEs. Furthermore, we require an expression for the seventh-order critical normal form coefficient for the higher-order approximation. While this coefficient was previously derived in \cite{sotomayor2007bifurcation}, we rederived it here and identified a missing term in one of their published expressions, though it was included in their example calculations. Additionally, we extend this derivation to obtain the seventh-order coefficients for DDEs. All components of the new predictors are summarised in \cref{sec:predictors}. The new equations have been implemented and their performance is shown through some examples in \cref{sec:examples}.

All of the computations were performed in Julia to create the examples. The corresponding code can be obtained upon request. The higher-order predictor for ODEs has also been implemented in the Julia software package BifurcationKit. 

\section{Parameter-dependent center manifold reduction and normalization}
\label{sec:theory}
In this section, we review the general technique for computing the normal form coefficients on the parameter-dependent center manifold for ODEs and DDEs. For DDEs, this includes a summary of relevant results from sun-star calculus. The general methods used to solve the resulting linear systems are also discussed. This includes additional special cases of representations required to obtain the higher-order coefficients in DDEs. 

\subsection{The center manifold reduction and normalization method for ODEs}\label{subsec:method_ODE}
Consider a system of ODEs depending on two parameters
\begin{equation}\label{eq:ODE_sys}
    \dot{x} = F(x,\alpha),
\end{equation}
where $x: \R \to \R^n$, $\alpha \in \R^2$ and $F:\R^n\times \R^2 \to \R^n$. Throughout this paper, we assume that $F$ is $C^k$-smooth for some sufficiently large integer $k$.  Suppose that $x_0 = 0$ is an equilibrium at $\alpha_0 = 0$, i.e. $F(0,0) = 0$. We denote the Jacobian matrix at the equilibrium as $A = D_xF(0,0)$. Let $n_c$ be the number of eigenvalues with zero real part and $T^c$ the corresponding critical eigenspace. For the generalized Hopf bifurcation, we have $n_c = 2$. A procedure to switch to the codim 1 bifurcation curves emanating from a codim 2 bifurcation point was introduced in \cite{beyn2002numerical}. We will follow the same approach below. 

For each sufficiently small $\|\alpha\|$, system \cref{eq:ODE_sys} has a smooth local $n_c$-dimensional invariant center manifold $W^c_{\alpha}$. At $x = 0$, $W_0^c$ is tangent to the critical eigenspace $T^c$ of $A$. Restricted to the center manifold, we can transform the system into a certain normal form using only smooth coordinate and parameter transformations:
\begin{equation}\label{eq:G}
\dot{w} = G(w,\beta), \quad G:\R^{n_c} \times \R^{2} \to \R^{n_c}.
\end{equation}
These normal forms are known for all five codim 2 bifurcations, and details can be found, for example, in \cite[Chapter 8]{kuznetsov2023elements}. From the normal form, we can derive an approximation for the codim 1 curves emanating from the codim 2 point. To relate the local behaviour of the normal form on the center manifold to our original system, we need a relation between the original parameters $\alpha$ and the unfolding parameters $\beta$:
\begin{equation}\label{eq:K_gen}
\alpha = K(\beta), \quad K:\R^2 \to \R^2.
\end{equation}
Furthermore, we need a parameterisation of the center manifold depending on the new parameters $\beta$:
\begin{equation}\label{eq:H_gen}
    x = H(w,\beta), \quad H: \R^{n_c} \times \R^{2} \to \R^n.
\end{equation}
By substituting equations $\cref{eq:G},\cref{eq:K_gen}$ and $\cref{eq:H_gen}$ into equation $\cref{eq:ODE_sys}$, we find the \textit{homological equation}:
\begin{equation}\label{eq:HOM_ODE}
H_w(w,\beta)G(w,\beta) = f(H(w,\beta),K(\beta)).
\end{equation}
This is essentially a consequence of the invariance of the center manifold. In the following, we will use multi-indices $\nu$ and $\mu$ to simplify the notation in the expansions of $G, H$ and $K$. A multi-index $\mu$ is defined as $\mu = (\mu_1, \ldots, \mu_n)$ for $\mu_i \in \N_0$ and we have that $w^{\mu} = w_1^{\mu_1}\ldots w_n^{\mu_2}$, $\mu! = \mu_1!\mu_2!\ldots\mu_n!$ and $|\mu| = \mu_1 + \ldots +\mu_n$. The general form of the normal form expansion is known and is expanded as
\begin{equation}\label{eq:expansion_G}
    G(w,\beta) = \sum_{|\nu| + |\mu| \geq 1} \frac{1}{\nu!\mu!}g_{\nu\mu}w^{\nu}\beta^{\mu}.
\end{equation}
Meanwhile, $H$ and $K$ are unknown and admit the expansions
\begin{equation}\label{eq:expansions_HK}
    H(w,\beta) =  \sum_{|\nu| + |\mu| \geq 1} \frac{1}{\nu!\mu!}H_{\nu\mu}w^{\nu}\beta^{\mu}, \quad K(\beta) = \sum_{|\mu|\geq 1}\frac{1}{\mu!} K_{\mu}\beta^{\mu}.
\end{equation}
All of these expansions will be truncated at some sufficiently high order. For the Taylor expansion of $F$, we define the multilinear forms
\begin{align*}
    B(u,v) &= \sum_{i,j=1}^n \frac{\partial^2F(x_0,\alpha_0)}{\partial x_i \partial x_j}u_i v_j, \quad &J_1u = \sum_{i = 1}^2  \frac{\partial F(x_0,\alpha_0)}{\partial \alpha_i}u_i,\\
    C(u,v,w) &= \sum_{i,j,k=1}^n \frac{\partial^3F(x_0,\alpha_0)}{\partial x_i \partial x_j \partial x_k}u_i v_j w_k, \quad &A_1(u,v) = \sum_{i = 1}^n \sum_{j =1}^2 \frac{\partial^2F(x_0,\alpha_0)}{\partial x_i \partial \alpha_j}u_i v_j, \\
    B_1(u,v,w) &= \sum_{i,j=1}^n\sum_{k = 1}^2 \frac{\partial^3F(x_0,\alpha_0)}{\partial x_i \partial x_j \partial \alpha_k}u_i v_j w_k. \quad &\textrm{ etc. }
\end{align*}
We use the letters $A,B,C,D,E,K,$ and $L$ to denote (in increasing order) the derivatives of $F$ with respect to its first argument evaluated at the critical point. The subscript denotes the order of the derivative of $F$ with respect to the parameters.

Substituting the expansions $\cref{eq:expansion_G}, \cref{eq:expansions_HK}$ along with the Taylor expansion of $F$, into the homological equation $\cref{eq:HOM_ODE}$ and collecting coefficients of the $w^{\nu}\beta^{\mu}$-terms results in a set of linear equations. These can be solved recursively for the unknown coefficients $g_{\nu\mu}, H_{\nu\mu}$, and $K_{\mu}$ by applying the Fredholm solvability condition.

\subsection{The center manifold reduction and normalization method for DDEs}\label{subsec:method_DDE}
As presented in \cite{janssens2010normalization}, the normalization technique for local bifurcations in ODEs can be lifted to the infinite-dimensional setting of DDEs. This normalization method was extended to include parameters in \cite{bosschaert2020switching}. We will closely follow the procedure from \cite{bosschaert2020switching}, summarizing results from sun-star calculus needed to apply this technique for the generalized Hopf bifurcation. 

We take the nonreflexive Banach space $X := C([-h,0], \R^n)$ and define for each $t \geq 0$ the history function $x_t: [-h, 0] \to \R^n$ at time $t$ as
\begin{equation*}
    x_t(\theta) := x(t + \theta), \quad \textrm{ for all } \theta \in [-h, 0].
\end{equation*}
 Consider the classical parameter-dependent DDE
\begin{equation}\label{eq:Classical_DDE}
    \dot{x}(t) = F(x_t, \alpha), \quad t \geq 0,
\end{equation}
where $F:X \times \R^2 \to \R^n$ is a $C^k$-smooth operator for sufficiently large integer $k$ and $0 < h < \infty$.
Assume that system $\cref{eq:Classical_DDE}$ satisfies $F(0,0) = 0$ and that the trivial equilibrium undergoes a generalized Hopf bifurcation at $\alpha=0$. The existence of a center manifold for DDEs has been rigorously established within the framework of sun-star calculus \cite{diekmann2012delay}. As for ODEs, the proof relies on a variations-of-constants formula describing the solutions. For DDEs, it is convenient to formulate this in the larger space $X^{\odot\star}$, the sun-star dual space of $X$. Here, we provide only the definitions and results from sun-star calculus needed for the normalisation method. For a full introduction to sun-star calculus and proofs, we refer to \cite{diekmann2012delay}.

There exists a unique matrix-valued function of normalized bounded variation $\zeta:[0,h] \to \R^{n\times n}$ such that the linear part of $\cref{eq:Classical_DDE}$ at $\alpha=0$ can written as
\begin{equation}\label{eq:linrep}
    D_1F(0,0)\varphi = \langle \zeta, \varphi \rangle := \int_0^h d\zeta(\theta)\varphi(-\theta).
\end{equation}
The above integral is of Riemann-Stieltjes type. With this notation, we can express the right-hand side of $\cref{eq:Classical_DDE}$ in terms of its linear and nonlinear parts:
\begin{equation*}
    F(\varphi, \alpha) = \langle \zeta, \varphi \rangle + D_2F(0,0)\alpha + G(\varphi, \alpha),
\end{equation*}
where the nonlinear part $G$ is a smooth operator satisfying $G(0,0) = 0$, $D_1G(0,0) = 0$, and $D_2G(0,0) = 0$. Let $T$ be the unique $C_0$-semigroup on $X$ corresponding to the linear part of $\cref{eq:Classical_DDE}$ at $0 \in X$ for the critical parameter value $\alpha = 0$. The generator $A$ of the semigroup $T$ is given by
\begin{equation}
    A\varphi = \dot{\varphi}, \quad \textrm{ with } \quad D(A) = \{ \varphi \in C^1 | \dot{\varphi}(0) = \langle \zeta, \varphi \rangle \}.
\end{equation}
When we go to the dual space $X^{\star}$ of $X$, we lose strong continuity with the adjoint semigroup $T^{\star}$. Therefore, we consider the maximal subspace of strong continuity $X^{\odot}$. The space $X^{\odot}$ has the representation 
\begin{equation}\label{eq:xsun_rep}
    X^{\odot} = \R^{n} \times L^1([0,h], \R^n),
\end{equation}
with the duality pairing between $\varphi^{\odot} = (c, g) \in X^{\odot}$ and $\varphi \in X$ given by 
\begin{equation}\label{eq:xsun_dual}
    \langle \varphi^{\odot}, \varphi \rangle = c^T\varphi(0) + \int_0^h g(\theta)\varphi(-\theta) d\theta.
\end{equation}
On the Banach space $X^{\odot}$ we have a $C_0$-semigroup $T^{\odot}$ with generator $A^{\odot}$. The dual space $X^{\odot\star}$ of $X^{\odot}$ has the representation
\begin{equation}\label{eq:xsunstar_rep}
    X^{\odot\star} = \R^{n} \times L^{\infty}([0,h], \R^n).
\end{equation}
On this space, we have the generator 
\begin{equation}
     A^{\odot\star}(\alpha, \varphi)= (\langle \zeta,\varphi\rangle,  \dot{\varphi} ), \quad \textrm{ with } \quad D(A^{\odot\star}) = \{(\alpha, \varphi) | \varphi \in \textrm{Lip}(\alpha) \},
\end{equation}
where Lip$(\alpha)$ denotes the subset of $L^{\infty}([-h,0], \C)$ consisting of Lipschitz continuous functions which assume the value $\alpha$ at $\theta = 0$. The duality pairing between $\varphi^{\odot\star} = (a, \psi) \in X^{\odot\star}$ and $\varphi^{\odot} = (c, g) \in X^{\odot}$ is given by
\begin{equation}\label{eq:xsunstar_dual}
    \langle \varphi^{\odot\star}, \varphi^{\odot} \rangle = c^Ta + \int_0^h g(\theta)\psi(-\theta) d\theta.
\end{equation}
We look again at the maximal subspace of strong continuity $X^{\odot\odot}$. There exists an injection $j: X \to X^{\odot\odot}$ defined by 
\begin{equation}\label{eq:injection}
    j\varphi = (\varphi(0), \varphi) \in  X^{\odot\odot} \textrm{ for all } \varphi \in X.
\end{equation}
We have that $X^{\odot\odot} = j(X)$. This property is also known as \textit{sun-reflexivity}. We will often move back and forth between the space $X$ and its sun-dual space $X^{\odot\star}$. 

Assume that there are $n_0 \geq 1$ eigenvalues of the linearization of $\cref{eq:Classical_DDE}$ at $\alpha=0$ on the imaginary axis, with a corresponding real $n_0$-dimensional center eigenspace $X_0$. Then  \cite[Corollary 20]{bosschaert2020switching} will imply the existence of a parameter-dependent local center manifold $\mathcal{W}_{\textrm{loc}}^c(\alpha)$ for $\cref{eq:Classical_DDE}$. As in the ODE case, we want to include a relation $\alpha = K(\beta)$ between the original parameters $\alpha$ and some new unfolding parameters $\beta$. Let $u: I \to X$ with $u(t) := x_t \in \mathcal{W}_{\textrm{loc}}^c(\alpha)$ be as in \cite[Corollary 20]{bosschaert2020switching}. Then, $u$ is differentiable on $I$ and satisfies the equation
\begin{equation}\label{eq:udot}
    j\dot{u}(t) = A^{\odot\star}ju(t) + (D_2F(0,0)K(\beta))r^{\odot\star} + G(u(t), K(\beta))r^{\odot\star}, \quad \textrm{ for all } t \in I.
\end{equation}
Here $wr^{\odot\star} = (w, 0) \in X^{\odot\star}$, for $w \in \R^n$. Now, choose a basis $\Phi$ of $X_0$. With respect to $\Phi$ and in terms of the new parameter $\beta$, we can consider the locally defined $C^k$-smooth parameterization $H: \R^{n_c} \times \R^2 \to X$ of the center manifold $\mathcal{W}_{\textrm{loc}}^c(\alpha)$. Let $z(t)$ be the coordinate with respect to $\Phi$ of the projection of $u(t)$ onto the center subspace $X_0$. Then $z: I \to \R^{n_0}$ satisfies a parameter-dependent ODE where the right-hand side is a $C^k$-smooth vector field that can be expanded as
\begin{equation}\label{eq:z}
    \dot{z} = \sum_{|\nu| + |\mu| \geq 1} \frac{1}{\nu!\mu!} g_{\nu\mu}z^{\nu}\beta^{\mu}.
\end{equation}
We may assume that $\cref{eq:z}$ is a smooth normal form in terms of the unfolding parameters $\beta$. Furthermore, we have that
\begin{align}\label{eq:uH}
    u(t) = H(z(t), \beta), \quad t \in I.
\end{align}
If we substitute equation $\cref{eq:uH}$ into the equation $\cref{eq:udot}$, we find the following \textit{homological equation}:
\begin{equation}\label{eq:HOM_DDE}
    A^{\odot\star}jH(z,\beta) + (J_1K(\beta))r^{\odot\star} + G(H(z,\beta),K(\beta))r^{\odot\star} = jD_zH(z,\beta)\dot{z},
\end{equation}
where $\dot{z}$ is given by $\cref{eq:z}$ and we defined $J_1 := D_2F(0,0)$.  The mappings $H$ and $K$ allow for the expansions
\begin{equation}\label{eq:expansions_HK_DDE}
    H(z,\beta) =  \sum_{|\nu| + |\mu| \geq 1} \frac{1}{\nu!\mu!}H_{\nu\mu}z^{\nu}\beta^{\mu}, \quad K(\beta) = \sum_{|\mu|\geq 1}\frac{1}{\mu!} K_{\mu}\beta^{\mu}.
\end{equation}
Meanwhile, the nonlinear part $G(\varphi, \alpha)$ can be expanded as
\begin{equation}\label{eq:nonlnearity_expansion}
    G(\varphi,\alpha) = \sum_{r + s > 1} \frac{1}{r! s!} D_1^r D_2^s F(0,0)(\varphi^{(r)}, \alpha^{(s)}),
\end{equation}
where $\varphi^{(r)} = (\varphi, \ldots, \varphi) \in X^r$, $\alpha^{(s)} = (\alpha, \ldots, \alpha) \in [\R^2]^s$ and $D_1^r D_2^s F(0,0): X^r \times [\R^2]^s \to \R^n$ is the mixed Fréchet derivative of order $r + s$ evaluated at $(0,0) \in X \times \R^2$. Just as in the ODE-case, we will indicate the multilinear forms in the expansion of $F$ with $B, C, D, E, K, L$ and for the parameter-dependent derivatives $A_i, B_i, C_i, $ etc. where the subscript $i$ indicated the number of derivatives with respect to the parameter. Thus,
\begin{align*}
    &J_1 = D_2F(0,0), \quad &B(u,u) = D_1^2F(0,0)(u,u), \\
    &A_1(u,\alpha) = D_1^1D_2^1 F(0,0)(u,\alpha), \quad &C(u,u,u) = D_1^3F(0,0)(u,u,u), \\
    &B_1(u,u,\alpha) = D_1^2D_2^1F(0,0)(u,u,\alpha), \quad  &\textrm{ etc. }
\end{align*}
For the case of discrete DDEs, explicit formulas for the computation of the multilinear forms are presented in \cite[Section 5]{bosschaert2020switching}. 

We can now substitute the expansions $\cref{eq:z}, \cref{eq:expansions_HK_DDE}$ and $\cref{eq:nonlnearity_expansion}$ into the homological equation $\cref{eq:HOM_DDE}$ and collect terms of equal powers $z^{\nu}\beta^{\mu}$. Then it is possible to solve recursively for the coefficients $g_{\nu\mu}, H_{\nu\mu}$ and $K_{\mu}$. Before we discuss how we can solve the resulting equations, we need some facts about the spectrum of the generator $A$. 
\subsubsection{The spectrum} To determine if bifurcations are present, we need to analyse the spectrum of the generator $A$ of the semigroup corresponding to the linear part of $\cref{eq:Classical_DDE}$. For this, we need the \textit{characteristic matrix function} which is defined as
\begin{equation}
    \Delta(z) = z I - \int_0^h e^{-z\theta} d\zeta(\theta), \quad \Delta : \C \to \C^{n \times n},
\end{equation}
and contains all spectral information about $A$. The integral in the above expression is of Riemann-Stieltjes type and $\zeta$ is the same from $\cref{eq:linrep}$. The eigenvalues of $A$ are given by the roots of the \textit{characteristic equation} 
\begin{equation}\label{eq:char_eq}
    \det \Delta(\lambda) = 0.
\end{equation}
We only have to concern ourselves with simple eigenvalues, i.e. eigenvalues for which both the geometric and algebraic multiplicities equal one. If $\lambda \in \C$ is a simple eigenvalue of $A$, there exist an eigenfunction $\varphi$ and an adjoint eigenfunction $\varphi^{\odot}$ such that
\[
A\varphi = \lambda \varphi, \quad A^{\star}\varphi^{\odot} = \lambda \varphi^{\odot}.
\]
Let $q, p \in \C^n$ such that
\[
\Delta(\lambda)q = 0, \quad p^T \Delta(\lambda) = 0.
\]
Then the corresponding eigenfunctions are given by
\begin{equation}\label{eq:eigen_phi}
    \varphi(\theta) = e^{\lambda\theta} q, \quad \theta \in [-h, 0]
\end{equation}
and
\begin{equation}\label{eq:eigen_phi_sun}
    \varphi^{\odot} = \left(p , \theta \mapsto p \int_{\theta}^h e^{\lambda(\theta - \tau)} d\zeta(\tau) \right), \quad \theta \in [-h, 0].
\end{equation}
Furthermore, it is possible to normalize the eigenfunctions to
\begin{equation*}
    \langle \varphi^{\odot}, \varphi \rangle = p\Delta'(\lambda)q = 1.
\end{equation*}
Here $\Delta'(\lambda)$ is the derivative of $z \mapsto \Delta(z)$ evaluated at $z = \lambda$ 
\begin{equation}\label{eq:Delta_first_deriv}
    \Delta'(\lambda) = I + \int_0^h \theta e^{-\lambda\theta} d\zeta(\theta).
\end{equation}
Proofs of the above results on the spectrum of $A$ can be found in \cite[Chapter IV]{diekmann2012delay}. Note that for $k \geq 2$ we have the higher-order derivatives
\begin{equation}\label{eq:Delta_nth_deriv}
    \Delta^{(k)}(\lambda) = (-1)^{k+1}\int_0^h\theta^k e^{-\lambda\theta} d\zeta(\theta).
\end{equation}
In the following, we will denote the second, third, and fourth derivatives as $\Delta''(\lambda)$, $\Delta'''(\lambda)$, and  $\Delta''''(\lambda)$ respectively.  

In the special case of the discrete DDE $\cref{eq:Discrete_DDE}$, the characteristic matrix is given by 
\begin{equation*}
    \Delta(z) = z I  - \sum_{j = 0}^m M_j e^{-z\tau_j}, \quad z \in \C,
\end{equation*}
where $M_j : = D_{1,j}f(0,0) \in \R^{n \times n}$ is the partial derivative of $f$ with respect to its $j$th state argument evaluated at the origin \cite[Section 6]{bosschaert2020switching}.
\paragraph{Remark} When dealing with the spectrum of $A$, it is necessary to complexify all of the above spaces and the linear operators acting on them. As previously remarked by \cite[Remark 2.2]{janssens2010normalization} this is not a trivial task. Fortunately, this has already been carried out in detail in \cite[Section III.7]{diekmann2012delay} and therefore will not be discussed here. 

\subsubsection{Solving the linear operator equations}\label{sec:DDE_eqsolve}
From the homological equation $\cref{eq:HOM_DDE}$ we will find equations of the form
\begin{equation}\label{eq:lin_DDE}
    (\lambda I - A^{\odot\star})(v_0, v) = (w_0,w),
\end{equation}
for some known $(w_0,w) \in X^{\odot\star}$, $\lambda \in \C$ and an unknown $(v_0,v)\in D(A^{\odot\star})$. We will only have to consider two cases, either $\lambda$ is a simple eigenvalue or $\lambda$ is not an eigenvalue.
\\[3mm]If $\lambda$ is not an eigenvalue, then equation $\cref{eq:lin_DDE}$ has a unique solution
\begin{equation}
    (v_0, v) = (\lambda I - A^{\odot\star})^{-1}(w_0,w).
\end{equation}
To compute the solutions, we need a representation for the solution $(v_0,v)$. The following lemma gives the general result.
\begin{lemma}
    [\cite{janssens2010normalization}, Lemma 3.3]\label{lemma:inverse_sol} Suppose that $\lambda$ is not an eigenvalue. Then the unique solution of $\cref{eq:lin_DDE}$ is given by
    \begin{equation}\label{eq:lemma_inverse1}
    v(\theta)= e^{\lambda\theta}v_0 + \int_{\theta}^0 e^{\lambda(\theta - \sigma)}w(\sigma)d\sigma \quad (\theta \in [-h,0]), 
    \end{equation}
    with
    \begin{equation}\label{eq:lemma_inverse2}
    v_0 = \Delta^{-1}(\lambda)\left\{ w_0 + \int_0^h d\zeta(\tau)\int_0^{\tau}e^{-\lambda\sigma}w(\sigma - \tau) d\sigma\right\}.
    \end{equation}
\end{lemma}
We will only encounter equations where the right-hand side of $\cref{eq:lin_DDE}$ has a form that allows us to express the above representations in terms of derivatives of the characteristic matrix function $\Delta(z)$. The following Corollary presents some useful cases that follow from the application of \cref{lemma:inverse_sol}. 
\begin{corollary}
    \label{cor:inverse_cases} Suppose that $\lambda$ is not an eigenvalue. We have the following special cases:
    \begin{enumerate}
        \item \cite[Lemma 3]{bosschaert2020switching} Suppose that $(w_0,w) = \left(w_0, \theta \mapsto e^{\lambda\theta}\Delta^{-1}(\lambda)\eta\right)$ for some fixed $\eta \in \C^n$. Then the unique solution $(v_0,v)\in D(A^{\odot\star})$ of $\cref{eq:lin_DDE}$ has the representation
        \[
            v_0 = v(0), \quad v(\theta) = \Delta^{-1}(\lambda)\left(e^{\lambda\theta}w_0 + [\Delta'(\lambda) - I - \theta\Delta(\lambda)]w(\theta)\right) .
        \]
        \item Let $k \geq 1$. Suppose that $(w_0,w) = \left(0, \theta \mapsto \theta^k e^{\lambda\theta}\Delta^{-1}(\lambda)\eta\right)$ for some fixed $\eta \in \C^n$. Then the unique solution $(v_0,v)\in D(A^{\odot\star})$ of $\cref{eq:lin_DDE}$ has the representation
        \[
            v_0 = v(0), \quad v(\theta) = \frac{1}{k + 1}e^{\lambda\theta}\Delta^{-1}(\lambda)[\Delta^{(k+1)}(\lambda) - \theta^{k + 1}\Delta(\lambda)]\Delta^{-1}(\lambda)\eta.
        \]
    \end{enumerate}
\end{corollary}
The special cases from \cref{cor:inverse_cases} can be combined to arrive at the following result:
\begin{corollary}\label{cor:inverse_case1}
    Suppose that $\lambda$ is not an eigenvalue and that the right-hand side of $\cref{eq:lin_DDE}$ has the representation
    \[
    (w_0,w) = \left(w_0, \theta \mapsto e^{\lambda\theta}\Delta^{-1}(\lambda)[\eta + \theta \xi_1 + \theta^2\xi_2]\right),
    \]
    for some fixed $\eta, \xi_1, \xi_2 \in \C^n$. Then the unique solution $(v_0,v)\in D(A^{\odot\star})$ of $\cref{eq:lin_DDE}$ has the following representation
    \begin{align*}
        v(\theta) &= e^{\lambda\theta}\Delta^{-1}(\lambda)\Bigl(w_0 + [\Delta'(\lambda) - I - \theta\Delta(\lambda)]\Delta^{-1}(\lambda)\eta \\
        &+ \frac{1}{2}[\Delta''(\lambda) - \theta^2\Delta(\lambda)]\Delta^{-1}(\lambda)\xi_1 +\frac{1}{3}[\Delta'''(\lambda) - \theta^3\Delta(\lambda)]\Delta^{-1}(\lambda)\xi_2 \Bigr),
    \end{align*}
    and $v_0 = v(0)$.
\end{corollary}
\begin{proof}
    Write 
    \begin{align*}
        (w_0,w) &=  \left(w_0, \theta \mapsto e^{\lambda\theta}\Delta^{-1}(z)\eta \right) + \left( 0, \theta \mapsto  \theta e^{\lambda\theta}\Delta^{-1}(\lambda)\xi_1 \right) \\
        &+ \left( 0, \theta \mapsto  \theta^2 e^{\lambda\theta}\Delta^{-1}(\lambda)\xi_2 \right).
    \end{align*} Using the linearity of the inverse operator $(\lambda I - A^{\odot\star})^{-1}$ we can apply the cases from \cref{cor:inverse_cases}.
\end{proof}
Finally, from this result, we can derive the following special case, which will aid our calculations in \cref{sec:coeff_DDE}.
\begin{corollary}
    \label{cor:inverse_case2} Suppose that $\lambda$ is not an eigenvalue. We have the following special case:
        Suppose that $(w_0,w)$ has the form
        \begin{align*}
            w_0 = w(0), \quad w(\theta) = e^{\lambda\theta}\Delta^{-1}(\lambda)\left(M + [\Delta'(\lambda) - \theta\Delta(\lambda)]\hat{\eta} + [\Delta''(\lambda) - \theta^2\Delta(\lambda)]\hat{\xi}\right),
        \end{align*}
        for some fixed $M, \hat{\eta}, \hat{\xi} \in \C^n$. Then the unique solution $(v_0,v)\in D(A^{\odot\star})$ of $\cref{eq:lin_DDE}$ has the representation
        \begin{align*}
            v(\theta) = e^{\lambda\theta}\Delta^{-1}(\lambda)\left( [\Delta'(\lambda) - \theta\Delta(\lambda)]w(0) - \frac{1}{2}[\Delta''(\lambda) - \theta^2\Delta(\lambda)]\hat{\eta} - \frac{1}{3}[\Delta'''(\lambda) - \theta^3\Delta(\lambda)]\hat{\xi} \right).
        \end{align*}
        and $v_0 = v(0)$.
\end{corollary}
\begin{proof}
    We can rewrite $w(\theta)$ as
    \[
    w(\theta) = e^{\lambda\theta}\Delta^{-1}(\lambda)[M + \Delta'(\lambda)\hat{\eta} + \Delta''(\lambda)\hat{\xi} - \theta\Delta(\lambda)\hat{\eta} - \theta^2\Delta(\lambda)\hat{\xi}].
    \]
    Now we can apply \cref{cor:inverse_case1} by taking $\eta = M + \Delta'(\lambda)\hat{\eta} + \Delta''(\lambda)\hat{\xi}$, $\xi_1 = - \Delta(\lambda)\hat{\eta}$ and $\xi_2 = -\Delta(\lambda)\hat{\xi}$. Finally, use that $w_0 = \Delta^{-1}(\lambda)\eta$ to simplify the expression. 
\end{proof}
If $\lambda$ is an eigenvalue, system $\cref{eq:lin_DDE}$ does not have a unique solution. As in the ODE case, a variant of the Fredholm solvability condition exists for operator equations of the form $\cref{eq:lin_DDE}$:
\begin{lemma}
    [\cite{janssens2010normalization}, Lemma 3.2] For arbitrary $\lambda$, a solution $(v_0,v)\in D(A^{\odot\star})$ to system $\cref{eq:lin_DDE}$ exists if and only if
    \[
    \langle (w_0,w), \varphi^{\odot} \rangle = 0, \quad \textrm{ for all } \quad\varphi^{\odot} \in N(\lambda I - A^{\star}).
    \]
\end{lemma}
During our computations in \cref{sec:coeff_DDE} we will refer to it as the Fredholm solvability condition. Similarly to before, we can use a bordered \textit{operator} inverse 
\[(\lambda I - A^{\odot\star})^{INV}: R(\lambda I - A^{\odot\star}) \to D(A^{\odot\star})\] to find a unique solution to system $\cref{eq:lin_DDE}$ satisfying $\langle (v_0, v), \varphi^{\odot} \rangle = 0$ for all $(w_0,w)$ for which $\cref{eq:lin_DDE}$ is consistent, i.e.  satisfies the Fredholm solvability condition. A general representation for $(\lambda I - A^{\odot\star})^{INV}$ in the case that $\lambda$ is a simple eigenvalue is given by the following proposition
\begin{prop}
    [\cite{janssens2010normalization}, Proposition 3.6]\label{prop:bordered_rep} Let $\lambda$ be a simple eigenvalue of $A$ with eigenvector $\varphi$ and adjoint eigenvector $\varphi^{\odot}$, such that $\langle \varphi^{\odot},\varphi \rangle = 1$.  Suppose that $\cref{eq:lin_DDE}$ is consistent for some given $(w_0,w) \in X^{\odot\star}$. Then the unique solution $(v_0, v) = (\lambda I - A^{\odot\star})^{INV}(w_0,w)$ satisfying $\langle (v_0, v), \varphi^{\odot} \rangle = 0$ is given by
    \begin{equation}\label{eq:borderded_rep_v}
    v_0 = \xi + \gamma q, \quad v(\theta) = e^{\lambda\theta}v_0 + \int_{\theta}^0 e^{\lambda(\theta - \sigma)}w(\sigma) d\sigma \quad (\theta \in [-h,0]), 
    \end{equation}
    with 
    \begin{equation}\label{eq:borderded_rep_xi}
    \xi = \Delta(\lambda)^{INV}\left[ w_0 + \int_0^h d\zeta(\tau) \int_0^{\tau} e^{-\lambda\sigma} w(\sigma- \tau) d\sigma \right], 
    \end{equation}
    and the constant $\gamma$ is given by
    \begin{equation}\label{eq:borderded_rep_gamma}
    \gamma = - p\Delta'(\lambda) \xi - p\int_0^h \int_\tau^h e^{-\lambda s} d\zeta(s) \int_{-\tau}^0 e^{-\lambda \sigma} w(\sigma) d\sigma d\tau.
    \end{equation}
\end{prop}
The expression for $\xi$ contains the bordered matrix inverse $\Delta(\lambda)^{INV}$. A unique solution $x = \Delta(\lambda)^{INV}y$ satisfying $p^T x = 0$ can be found by solving the following bordered system 
\begin{equation*}
    \begin{pmatrix}
        \Delta(\lambda) & q \\
        p^T & 0 
    \end{pmatrix} \begin{pmatrix}
        x \\
        s
    \end{pmatrix} = \begin{pmatrix}
        y \\
        0
    \end{pmatrix},
\end{equation*}
for $(x,s) \in \C^{n+1}$. The derivation of the higher order coefficients requires additional cases for the representation of $(\lambda I - A^{\odot\star})^{INV}$ than the one presented in \cite{bosschaert2020switching}. A more general case follows from applying \cref{prop:bordered_rep}.
 \begin{corollary}\label{cor:bordered_case}
    Suppose in addition to the assumptions from \cref{prop:bordered_rep} that $(w_0,w)$ satisfies
    \[
    (w_0, w) = \left(w_0, e^{\lambda\theta}(\xi_0 + \theta \xi_1  + \ldots + \theta^m \xi_m) \right) , 
    \]
    for some $m \in \N$ and constant vectors $\xi_k \in \C^n$, $0 \leq k \leq m$.
    Then the unique solution $(v_0, v) = (\lambda I - A^{\odot\star})^{INV}(w_0,w)$ satisfying $\langle (v_0, v), \varphi^{\odot} \rangle = 0$ is given by
    \begin{equation}\label{eq:bordered_case_rep}
    v_0 = \xi+ \gamma q, \quad v(\theta) = e^{\lambda\theta}\left(v_0 - \theta \xi_0 - \frac{1}{2}\theta^2\xi_1  - \ldots - \frac{1}{m +1}\theta^{m+1} \xi_m\right),
    \end{equation}
    where
    \[
    \xi = \Delta(\lambda)^{INV}\left[w_0 + [\Delta'(\lambda) - I]\xi_0 + \frac{1}{2}\Delta''(\lambda)\xi_1 + \ldots + \frac{1}{m + 1}\Delta^{(m + 1)}(\lambda)\xi_m\right] ,
    \]
    and
    \[
    \gamma =- p\Delta'(\lambda) \xi + \frac{1}{2}p\Delta''(\lambda)\xi_0 + \frac{1}{6}p\Delta'''(\lambda)\xi_1 + \ldots + \frac{1}{(m + 1)(m + 2)}\Delta^{(m + 2)}(\lambda)\xi_m.
    \]
\end{corollary}
\paragraph{Remark} In \cref{sec:coeff_DDE}, we will encounter $(w_0,w)$ of the form
\[
w_0 = \eta + \xi_0, \quad w(\theta) = e^{\lambda\theta}(\xi_0 + \theta \xi_1  + \ldots + \theta^m \xi_m).
\]
For the bordered inverse of this case, we will use the notation $v = B_{\lambda}^{INV}(\eta, \xi_0, \xi_1, \ldots, \xi_m)$. For instance, if $(w_0,w) = (\eta, 0) + \kappa(q, \varphi)$, we write $v = B_{\lambda}^{INV}(\eta, \kappa q)$.

\section{Higher order \texorpdfstring{$LPC$}{LPC} curve approximation for the normal form}
\label{sec:predictor_normalform}
In this section, we first derive a parameter approximation of the $LPC$ curve for the normal form, followed by an approximation of the cycle period along this curve. We conclude with a detailed discussion of the necessary coefficients for the center manifold and parameter approximations when generalising to $n$-dimensional systems of ODEs.

\subsection{Approximmation of the \texorpdfstring{$LPC$}{LPC} curve derived from the normal form}
Assume that at $\alpha = 0$ there is an equilibrium at the origin with only one pair of purely imaginary simple eigenvalues
\[
\lambda_{1,2} = \pm i\omega_0, \quad \omega_0 > 0.
\]
Restricted to the center manifold near a generalized Hopf bifurcation, the system $\cref{eq:ODE_sys}$ can be transformed to the following parameter-dependent normal form
\[
\dot{w} = \lambda(\alpha)w + c_1(\alpha)w|w|^2 + c_2(\alpha)w|w|^4 + c_3(\alpha)w|w|^6 + O(|w|^8), \quad w \in \C,
\]
where $\lambda(0) = i\omega_0$, $d_1 = \Re{(c_1(0))} = 0$ and $d_2 = \Re{(c_2(0))} \neq 0$. Furthermore, we have the following \textit{transversaility} condition:
\begin{equation}\label{eq:trans_cond}
    \textrm{The map } \alpha \mapsto (\Re{(\lambda(\alpha))}, \Re{(c_1(\alpha))}) \textrm{ is regular at } \alpha = 0.
\end{equation}
If this condition is satisfied, we may introduce new parameters:
\[(\beta_1(\alpha), \beta_2(\alpha)) = (\Re{(\lambda(\alpha))}, \Re{(c_1(\alpha))}).\]
Then, for $\|\alpha\|$ small enough, the normal form can be expressed in terms of $\beta$ in the following way
\begin{equation}\label{eq:normalform_beta}
\dot{w} = (i\omega_0 + \beta_1 +ib_1(\beta))w + (\beta_2 + ib_2(\beta))w|w|^2 + c_2(\beta)w|w|^4 + c_3(\beta)w|w|^6 + O(|w|^8),
\end{equation}
where $b_1$ and $b_2$ are real valued functions with $b_1(0) = 0$ and $b_2(0) = \Im{(c_1(0))}$. Note that we write $c_i(\beta)$ instead of $c_i(\alpha(\beta))$ for convenience. To get a higher-order approximation of the $LPC$ curve, we first substitute $w = \rho e^{i\psi}$ into $\cref{eq:normalform_beta}$. Taking the real part yields the following amplitude equation
\begin{align}
    \dot{\rho} &= \rho(\beta_1 + \beta_2\rho^2 + \Re{}(c_2(\beta))\rho^4 + \Re{}(c_3(\beta))\rho^6 + O(\rho^7)). \label{eq:polar_normal_rho}
\end{align}
\subsubsection{LPC curve in the amplitude equation}\label{sec:LPC_appr_der}
To derive an approximation for the $LPC$ curve we expand
\begin{align*}
\Re{}(c_2(\beta)) &= d_2 + a_{3210}\beta_1 + a_{3201}\beta_2 + O(\|\beta\|^2),\\
\Re{}(c_3(\beta)) &= d_3 + O(\|\beta\|).
\end{align*}
The $LPC$ curve corresponds to a double equilibrium of the right-hand side of equation $\cref{eq:polar_normal_rho}$, which occurs when both the right-hand side and its first derivative with respect to $\rho$ vanish. This leads to the following system of equations:
\begin{equation}\label{eq:double_eq}
    \begin{cases}
        \beta_1 + \beta_2\rho^2 + \Re{}(c_2(\beta))\rho^4 + \Re{}(c_3(\beta))\rho^6 + O(\rho^7) = 0, \\
        \beta_2\rho + 2\Re{}(c_2(\beta))\rho^3 + 3 \Re{}(c_3(\beta))\rho^5 + O(\rho^6) = 0.
    \end{cases}
\end{equation}
Let us define $P:\R\times \R^2 \to \R^2$ by
\begin{equation}\label{eq:Prhobeta}
    P(\rho,\beta) = \begin{pmatrix}
        \beta_1 + \beta_2\rho^2 + \Re{}(c_2(\beta))\rho^4 + \Re{}(c_3(\beta))\rho^6 + O(\rho^7) \\
        \beta_2 + 2\Re{}(c_2(\beta))\rho^2 + 3 \Re{}(c_3(\beta))\rho^4  + O(\rho^5)
    \end{pmatrix}.
\end{equation}
Equation $\cref{eq:double_eq}$ is equivalent to $P(\rho,\beta) = 0$. Note that $P(0,0) = 0$ and $DP_{\beta}(0,0)$ is invertible. Thus, by the Implicit Function Theorem, there exists a unique, locally smooth function $\beta(\rho) = (\beta_1(\rho),\beta_2(\rho))$ near $\rho = 0$ such that $\beta(0) = 0$ and $P(\rho,\beta(\rho)) = 0$. Consequently, we can expand $\beta_1$ and $\beta_2$ as power series in $\rho$:
\begin{align*}
    \beta_i &= m_{i,1}\rho + m_{i,2}\rho^2 + m_{i,3}\rho^3 + m_{i,4}\rho^4 + m_{i,5}\rho^5 + m_{i,6}\rho^6 + O(\rho^7), \quad i = 1,2.
\end{align*}
 We can solve for the unknown coefficients by substituting these expansions into $P(\rho,\beta) = 0$ and collecting terms up to order $\rho^6$ in the first equation and $\rho^4$ in the second. This yields the following asymptotic expansions for the $\beta$ parameters:
\begin{align*}
    \beta_1 &= d_2\rho^4 + 2(d_3 - a_{3201}d_2)\rho^6 + O(\rho^7), \\
    \beta_2 &= -2d_2\rho^2 + (4a_{3201}d_2 - 3d_3)\rho^4 + O(\rho^5).
\end{align*}
%We see that the 7-th order coefficient $d_3$ has appeared as was already remarked in \cite[Remark 3.3.2]{meijer2006codimension} for the Chenciner bifurcation of iterated maps.
For small $\rho = \varepsilon > 0$ an approximation of the $LPC$ curve is given by:
\begin{equation}\label{eq:LPC}
    (\rho,\beta_1,\beta_2) = (\varepsilon, d_2\varepsilon^4 + 2(d_3 - a_{3201}d_2)\varepsilon^6+ O(\varepsilon^7), -2d_2\varepsilon^2 + (4a_{3201}d_2 - 3d_3)\varepsilon^4 + O(\varepsilon^5)).
\end{equation}

\subsubsection{Period approximmation}
To approximate the period $T$ of the cycles on the $LPC$ curve near the generalized Hopf bifurcation we use that
\begin{equation}\label{eq:period_int}
\int_0^T \dot{\psi} dt = 2\pi.
\end{equation}
To obtain the approximation of $\dot{\psi}$, we first take the imaginary part of equation $\cref{eq:normalform_beta}$ after substituting $w = \rho e^{i\psi}$ and partially truncate higher order terms. This results in the equation:
\begin{align}\label{eq:polar_normal_psi}
    \dot{\psi} &= \omega_0 + b_1(\beta) + b_2(\beta)\rho^2 + \Im{(c_2(0))}\rho^4 + \Im{(c_3(0))}\rho^6.
\end{align}
For $b_1(\beta)$ and $b_2(\beta)$ we make the following expansions
\begin{align*}
    b_1(\beta) &= b_{1,10}\beta_1 + b_{1,01}\beta_2 \textcolor{gray}{+ \frac{1}{2}b_{1,20}\beta_1^2 + b_{1,11}\beta_1\beta_2}+ \frac{1}{2}b_{1,02}\beta_2^2  \textcolor{gray}{+O(\|\beta\|^3)},\\
    b_2(\beta) &= \Im{(c_1(0))} + \textcolor{gray}{b_{2,10}\beta_1} + b_{2,01}\beta_2 + \textcolor{gray}{ O(\|\beta\|^2)}.
\end{align*}
Since we are only interested in an approximation of the period up to fourth order in $\varepsilon$, the terms in grey would not appear in the period approximation. Their inclusion here is meant to demonstrate this explicitly. Now, substitute the parameter expansions $\cref{eq:LPC}$ of the $LPC$ curve into equation $\cref{eq:polar_normal_psi}$ and solve the integral $\cref{eq:period_int}$. For small $\rho = \varepsilon > 0$ this yields the following period approximation:
\begin{align}\label{eq:period}
    T &= 2\pi/(\omega_0 + (\Im{(c_1(0))} -2d_2b_{1,01})\varepsilon^2 + [d_2b_{1,10} + (4a_{3201}d_2 - 3d_3)b_{1,01} + 2d_2^2b_{1,02} \nonumber \\
    &- 2d_2b_{2,01} + \Im{(c_2(0))}]\varepsilon^4 + O(\varepsilon^5)).
\end{align}

\subsection{Coefficients in the center manifold and para\-me\-ter transformation approximmations}\label{sec:necessary_coeffs}
To find all the coefficients that need to be included in the approximations $K$ and $H$, we must determine which terms in the normal form \cref{eq:normalform_beta} affect the $LPC$ curve approximation \cref{eq:LPC} up to the desired order. As noted in \cite{bosschaert2024bt}, it is equally important to account for the influence of terms \textit{not} included in the normal form. Any term affecting the $LPC$ curve approximation \cref{eq:LPC} must be eliminated, i.e. transformed into higher-order terms, to maintain the accuracy of our predictor at the intended order. 

For example, consider a term like $w\beta_2^2$ appearing in the normal form after reducing the system to the center manifold. If we derive the parameter approximation for the $LPC$ curve to the same order as \cref{eq:LPC}, following the method outlined in \cref{sec:LPC_appr_der}, this will introduce two additional terms in the expression for $\beta_1$. Thus, $w\beta_2^2$ impacts the $LPC$ curve approximation, indicating that it should be transformed away. To achieve this, we need to include the coefficient $K_{02}$ in the parameter approximation \cref{eq:K_gen}.

To derive our parameter approximation of the $LPC$ curve it is enough to consider the truncated normal form:
\begin{align*}
\dot{w} &= (i\omega_0 + \beta_1 +ib_1(\beta))w + (\beta_2 + ib_2(\beta))w|w|^2 + (c_2(0) + a_{3201}\beta_2)w|w|^4 + d_3w|w|^6.
\end{align*}
We now include terms of the form $g_{nmkl}w^n\bar{w}^m\beta_1^k \beta_2^l$ for $n,m,k,l \in \N$ with $n+m+k+l \geq 1$ to this normal form, and determine if they influence the parameter approximation. 

The amplitude equation $\dot{\rho}$ is derived by taking the real part of the normal form after substituting $w = \rho e^{i\psi}$. All terms for which $n -m - 1\neq 0$ will be resonant. These will certainly influence the predictor up to a certain order. Let us write $a_{nmkl}(\psi) = \Re{}\{g_{nmkl}e^{i\psi(n - m - 1)}\}$. Then the amplitude equation will have the following form
\begin{align*}
    \dot{\rho} = \beta_1\rho + \beta_2\rho^3 + (d_2  + a_{3201}\beta_2)\rho^5  + d_3\rho^7 + a_{nmkl}(\psi)\rho^{n + m}\beta_1^k \beta_2^l.
\end{align*}
To derive a parameter approximation of the $LPC$ curve, we look for a double zero in the equation
\begin{align}\label{eq:amp_zero_res}
      \beta_1\rho + \beta_2\rho^3 + (d_2  + a_{3201}\beta_2)\rho^5  + d_3\rho^7  + a_{nmkl}(\psi)\rho^{n + m}\beta_1^k \beta_2^l = 0.
\end{align}
 The goal is to determine the combinations of $n,m,k,l$ for which the coefficient $a_{nmkl}(\psi)$ will appear in our expansions of $\beta_1$ and $\beta_2$, consequently affecting our $LPC$ curve approximation. This will indicate which terms need to be removed in the normal form and thus which coefficients $H_{nmkl}$, $K_{kl}$ need to be included in the center manifold and parameter approximations \cite{bosschaert2024bt}. We will first discuss the case when  $n + m = 0$, and then provide the coefficients for $n +  m \geq 1$. 
 
 \textbf{Case $n + m = 0:$} The equations that need to be satisfied for a double zero are in this case given by
\begin{align}
    \beta_1\rho + \beta_2\rho^3 + (d_2 + a_{3201}\beta_2)\rho^5  + d_3\rho^7  + a_{00kl}(\psi)\beta_1^k \beta_2^l &= 0, \label{eq:case0-1}\\
    \beta_1 + 3\beta_2\rho^2 + 5(d_2  + a_{3201}\beta_2)\rho^4 + 7d_3\rho^6 = 0.\label{eq:case0-2}
\end{align}
As before, we expand
\begin{align*}
    \beta_i &= m_{i,1}\rho + m_{i,2}\rho^2 + m_{i,3}\rho^3 + m_{i,4}\rho^4  + m_{i,5}\rho^5 + m_{i,6}\rho^6 + O(\rho^7), \quad i = 1,2.
\end{align*}
From equation \cref{eq:case0-2} it follows that $m_{1,1} = m_{1,2} = 0$. We can simplify things a bit by setting $m_{i,1} = m_{i,3} = m_{i,5} =  0$ for $i = 1,2$. To see why, let us examine all the $\rho$-terms in equation \cref{eq:case0-1}. This yields 
\[
 a_{0001}(\psi)m_{2,1} + a_{0010}(\psi)m_{1,1} = 0.
\]
Since this must hold for all $\psi$, it follows that $m_{2,1}$ and $m_{1,1}$ have to equal zero. These terms are already zero in our $LPC$ approximation, so this does not present a problem. Now collect all $\rho^2$-terms in the first equation. Using that $m_{2,1}= m_{1,1} = 0$, we find the following equation
\[
a_{0001}(\psi)m_{2,2} + a_{0010}(\psi)m_{1,2}  = 0.
\]
This can only hold for all $\psi$ if $m_{1,2} = m_{2,2} = 0$. For $m_{1,2},$ this is not an issue since this term is already zero in the approximation. However, for $m_{2,2}$, it is a problem because this term is not zero in our $LPC$ curve approximation. Thus, if a term $g_{0001}\beta_2$ is present in the normal form, it will influence our predictor. Therefore, this term must be removed, implying that the coefficients $H_{0001}$ and $K_{01}$ must be included in the center manifold and parameter approximations. Whenever a coefficient like $m_{1,3}$, which is zero in our $LPC$ approximation, appears in front of a term $a_{00kl}(\psi)$, this term can be removed by simply setting $m_{1,3} = 0$ without altering the predictor. There is only a problem if the coefficient in front of $a_{00kl}(\psi)$ depends on coefficients that are nonzero in our predictor. Thus, it is enough to substitute
\begin{equation}\label{eq:coeffderiv_subs}
        \beta_1 = m_{1,4}\rho^4  + m_{1,6}\rho^6 , \textrm{ and } \beta_2 = m_{2,2}\rho^2 + m_{2,4}\rho^4,
\end{equation}
when collecting terms.
To arrive at the desired order of the predictor, we need to collect terms up to $\rho^7$ in the first equation and terms up to $\rho^6$ in the second equation. Therefore, all terms with $k,j \in \N$ for which $\beta_1^k\beta_2^l$ contains a term $\rho^j$, $j \leq 7$ need to be removed.  This implies that we need the coefficients $H_{0010}, H_{0001}, H_{0002}, H_{0011}$, and $H_{0003}$ in the center manifold approximation and the coefficients $K_{10},K_{01},K_{02},K_{11}$, and $K_{03}$ in the parameter approximation.

\textbf{Case $n + m \geq 1:$} Now consider the following two equations
\begin{align*}
    \beta_1\rho + \beta_2\rho^3 + (d_2 + a_{3201}\beta_2)\rho^5  + d_3\rho^7  + a_{nmkl}(\psi)\rho^{n+m}\beta_1^k \beta_2^l &= 0, \\
    \beta_1 + 3\beta_2\rho^2 + 5(d_2 + a_{3201}\beta_2)\rho^4 + 7d_3\rho^6 + a_{nmkl}(\psi)(n + m - 1)\rho^{n+m - 1}\beta_1^k \beta_2^l&= 0.
\end{align*}
By the same reasoning as in the previous case it is enough to substitute \cref{eq:coeffderiv_subs}. To derive the parameter approximation of the predictor up to sixth order in $\rho$, we need to collect all terms up to $\rho^7$ in the first equation and all terms up to $\rho^6$ in the second equation. If the term $a_{nmkl}(\psi)\rho^{n+m}\beta_1^k \beta_2^l$ contains terms of order $\rho^7$ or lower, this will influence the approximation. 

If $k = l = 0$, then a term $a_{nmkl}(\psi)$ will influence the predictor for all $1 \leq n + m \leq 7$ and thus all of these need to be removed. This implies that we need to include all coefficients $H_{nm00}$ with $1 \leq n + m \leq 7$ in the center manifold approximation. 
 
 Now we consider the cases for which $k$ or $l$ are nonzero. For $n + m = 1$, we need to remove all $a_{nmkl}(\psi)\rho\beta_1^k \beta_2^l$ terms for which $\beta_1^k \beta_2^l$ contains terms of $\rho^6$ or less. This is the case for $(kl) = (10),(01),(02), (11)$ and $(03)$. So we need $H_{nm10},H_{nm01},H_{nm02}, H_{nm11}$ and $H_{nm03}$ with $n + m = 1$. Continuing in this way, we find that all coefficients that need to be included in the center manifold approximation are
\begin{align*}
    &H_{nm10}, H_{nm01},H_{nm02},H_{nm11}, H_{nm03}, &n + m &= j, j \in \{0, 1\},\\
    &H_{nm10}, H_{nm01},H_{nm02}, &n + m &= j, j \in \{2,3\}, \\
    &H_{nm01},&n + m &= j, j \in \{4,5\},\\
    &H_{nm00}, &n + m &= j, j \in \{1,2,3,4,5,6,7\}.
\end{align*}
Finally, all coefficients that need to be included in the parameter transformation are $K_{10}$, $K_{01}$, $K_{02}$, $K_{11}$, and $K_{03}$.

\section{Coefficients of the parameter-dependent normal form and the predictor for ODEs}
\label{sec:coeffs_ODE}
In this section, we derive equations for the coefficients that we need in our predictor for ODEs. Assume that system $\cref{eq:ODE_sys}$ has an equilibrium at the origin at $\alpha = (0,0) \in \R^2$ with only one pair of purely imaginary simple eigenvalues
\[
\lambda_{1,2} = \pm i\omega_0, \quad \omega_0 > 0.
\]
Additionally, we assume that all the other eigenvalues have non-zero real parts. This allows us to introduce the complex eigenvectors $p,q \in \C^n$ satisfying
\[
Aq = i\omega_0 q, \quad A^Tp = -i\omega_0 p, \textrm{ and }\bar{q}^Tq = \bar{p}^Tq = 1.
\]
Furthermore, we assume that the first Lyapunov coefficient $l_1(0) = 0$ and the second Lyapunov coefficient $l_2(0) \neq 0$. The critical real eigenspace $T^c$ corresponding to $\lambda_{1,2}$ is now two dimensional and we can represent each $y \in T^c$ in terms of the complex coordinate $w = \langle p, y \rangle$ as 
\[
y = wq + \bar{w} \bar{q}.
\]Then, the homological equation $\cref{eq:HOM_ODE}$ has the form
\begin{equation}\label{eq:HOM_GH_ODE}
    D_wH(w,\bar{w},\beta) \dot{w} + D_{\bar{w}}H(w,\bar{w},\beta) \dot{\bar{w}} = F(H(w,\bar{w},\beta), K(\beta)).
\end{equation} Under the assumption that the transversality condition $\cref{eq:trans_cond}$ is satisfied, the truncated normal form restricted to the two-dimensional center manifold can be expressed in terms of the unfolding parameters $\beta = (\beta_1,\beta_2)$:
\begin{align}\label{eq:wdot}
    \dot{w} &= (i\omega_0 + \beta_1 +ib_1(\beta))w + (\beta_2 +  ib_2(\beta))w|w|^2 + (c_2(0) + g_{3201}\beta_2)w|w|^4 \nonumber \\
    &+ c_3(0)w|w|^6 , 
\end{align}
where it is enough to expand
\begin{align}
b_1(\beta) &= b_{1,10}\beta_1 + b_{1,01}\beta_2 + b_{1,11}\beta_1\beta_2+ \frac{1}{2}b_{1,02}\beta_2^2 + \frac{1}{6}b_{1,03}\beta_2^3,\label{eq:b1_expansion}\\
%\nonumber\\
%&+ O(|\beta_1|^2 + |\beta_1 ||\beta_2|^2 + |\beta_1|^2|\beta_2|  + \|\beta\|^4)
b_2(\beta) &= \Im{(c_1(0))} + b_{2,10}\beta_1 + b_{2,01}\beta_2  + b_{2,11}\beta_1\beta_2 + \frac{1}{2}b_{2,02}\beta_2^2  + \frac{1}{6}b_{2,03}\beta_2^3 .\label{eq:b2_expansion}%\nonumber\\
%&+  O(|\beta_1|^2 + |\beta_1 ||\beta_2|^2 + |\beta_1|^2|\beta_2|  + \|\beta\|^4).
\end{align}
Under the assumption that $F$ is sufficiently smooth, the truncated Taylor expansion of $F$ is taken as
\begin{align}\label{eq:F_Taylor}
    F(x,\alpha) &= Ax +J_1\alpha + A_1(x,\alpha) + \frac{1}{2}B(x,x) + \frac{1}{2}J_2(\alpha,\alpha) + \frac{1}{6}C(x,x,x) + \frac{1}{2}B_1(x,x,\alpha) \nonumber\\
    &+ \frac{1}{2}A_2(x,\alpha,\alpha) + \frac{1}{6}J_3(\alpha,\alpha,\alpha) +  \frac{1}{24}D(x,x,x,x) + \frac{1}{6}C_1(x,x,x,\alpha) + \frac{1}{4}B_2(x,x,\alpha,\alpha) \nonumber\\
    &+ \frac{1}{6}A_3(x,\alpha,\alpha,\alpha) + \frac{1}{120}E(x,x,x,x,x) + \frac{1}{24}D_1(x,x,x,x,\alpha) + \frac{1}{12}C_2(x,x,x,\alpha,\alpha) \nonumber\\
    &+ \frac{1}{12}B_3(x,x,\alpha,\alpha,\alpha)+\frac{1}{720}K(x,x,x,x,x,x) +  \frac{1}{120}E_1(x,x,x,x,x,\alpha) \nonumber\\
    &+ \frac{1}{36}C_3(x,x,x,\alpha,\alpha,\alpha) +\frac{1}{5040}L(x,x,x,x,x,x,x) .
\end{align}
The parameterization of the center manifold is expanded as
\begin{align}\label{eq:H_GH}
H(w,\bar{w},\beta)& = qw + \bar{q}\bar{w} + \sum_{n + m = 2}^7 \frac{1}{n!m!}H_{nm00}w^n\bar{w}^m + \sum_{n + m = 0}^5 H_{nm01}\frac{1}{n!m!}w^n\bar{w}^m\beta_2 \nonumber \\
&+ \sum_{n + m=0}^3 \frac{1}{n!m!}H_{nm10}w^n\bar{w}^m\beta_1 + \sum_{n + m=0}^3 \frac{1}{2n!m!}H_{nm02}w^n\bar{w}^m\beta_2^2 \nonumber\\
&+ \sum_{n + m=0}^1 \frac{1}{n!m!}H_{nm11}w^n\bar{w}^m\beta_1 \beta_2 + \sum_{n + m=0}^1 \frac{1}{6n!m!}H_{nm03}w^n\bar{w}^m\beta_2^3 \nonumber\\
&\textcolor{gray}{+ \frac{1}{6}H_{1103}w\bar{w}\beta_2^3 + \frac{1}{12}H_{2003}w^2\beta_2^3 + \frac{1}{12}H_{2103}w^2\bar{w}
\beta_2^3 }.
\end{align}
Note that the image of $H$ lies in $\R^n$ and thus we have that $H_{jk\mu} = \overline{H}_{kj\mu}$. The last three terms, marked in grey, are not needed to approximate the periodic orbit in phase space. However, as will be clear in \cref{sec:parameter_coef_ODE}, we do need the expressions for $H_{1103},H_{2003}$ and $H_{2103}$ to derive the coefficients $K_{03},H_{0003}$ and $H_{1003}$. The relation between the parameters is expanded as
\begin{align}\label{eq:K}
K(\beta) &= K_{10}\beta_1 + K_{01}\beta_2 + \frac{1}{2}K_{02} \beta_2^2 + K_{11}\beta_1\beta_2 + \frac{1}{6}K_{03}\beta_2^3. %\nonumber\\
%&+  O(|\beta_1|^2 + |\beta_1 ||\beta_2|^2 + |\beta_1|^2|\beta_2|  + \|\beta\|^4).
\end{align}
Although the coefficients $b_{1,11}$, $b_{1,03}$ in the expansion of $b_1(\beta)$ and $b_{2,02}, b_{2,11}, b_{2,03}$ in de expansion of $b_2(\beta)$ do not appear in our predictor, we will need their expressions from the homological equation to solve for $K_{02}, K_{11}$ and $K_{03}$. All the equations collected from the homological equation $\cref{eq:HOM_GH_ODE}$ are listed in the supplementary material (S2).

\begin{figure}[!ht]
    \centering
\includegraphics[width=\textwidth]{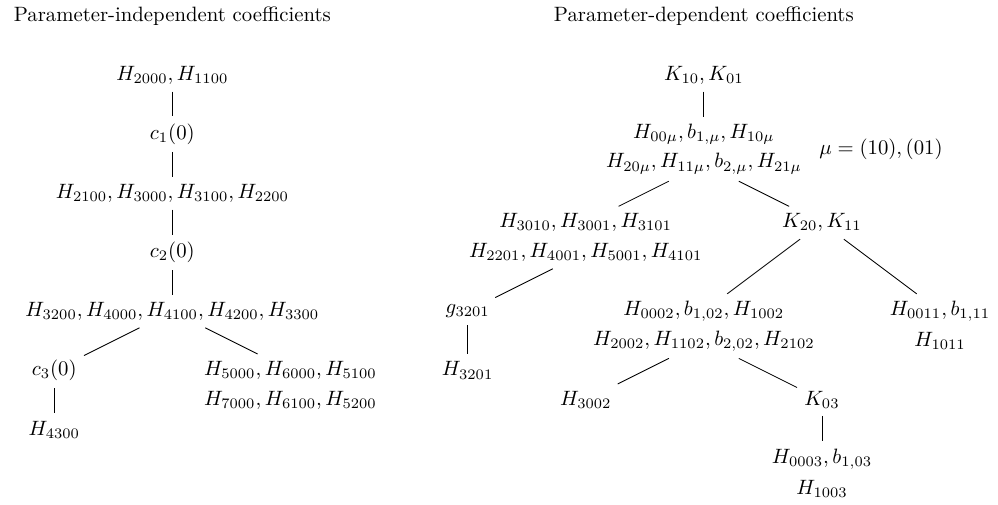}
    \caption{Schematic overview of all the coefficients that need to be determined.}
    \label{fig:coeffsdiagram}
\end{figure}

An overview of all the coefficients that must be determined for the higher-order predictor is given in \cref{fig:coeffsdiagram}. In the next two subsections, we derive all necessary critical normal form coefficients from $\cref{eq:wdot}$ and parameter-dependent coefficients from $\cref{eq:K}$. The remaining center manifold coefficients are presented in the supplementary material (S1.3).

\subsection{Critical normal form coefficients}
The critical normal form coefficients up to fifth order have been derived in previous works. For completeness, we include the derivation following \cite[Section 8.7.3]{kuznetsov2023elements}. The expression for the first Lyapunov coefficient $l_1 = \frac{1}{\omega_0}\Re{(c_1(0))}$ can be derived from the $w^2$, $w\bar{w}$ and $w^2\bar{w}$ terms in the homological equation. These terms yield the equations
\begin{align}
    H_{2000} &= (2i\omega_0I_n - A)^{-1}B(q,q), \\
    H_{1100} &= -A^{-1}B(q,\bar{q}), \\
    (i\omega_0I_n - A)H_{2100} &= 2B(q, H_{1100}) + B(\bar{q},H_{2000}) + C(q,q,\bar{q})- 2c_1(0)q \label{eq:H2100}.
\end{align}
Since the last equation is singular, we can use the Fredholm solvability condition to find
\begin{equation}
    c_1(0) = \frac{1}{2}\bar{p}^T[2B(q, H_{1100}) + B(\bar{q},H_{2000}) + C(q,q,\bar{q})].
\end{equation}
The vector $H_{2100}$ satisfying $ \bar{p}^T H_{2100} = 0$ can be found by solving the corresponding bordered matrix system, i.e. 
\begin{equation}
    H_{2100}=A^{INV}_{i\omega_0}\left(2B(q, H_{1100}) + B(\bar{q},H_{2000}) + C(q,q,\bar{q})- 2c_1(0)q\right),
\end{equation}
where we use the shorthand notation $A^{INV}_{i\omega_0} =(i\omega_0I_n - A)^{INV}$. The vector $w = A^{INV}_{i\omega_0}v$ can be found by solving the following non-singular bordered system
\begin{align}\label{eq:borderedsys}
    \begin{pmatrix}
    i\omega_0I_n - A & q \\
    \bar{p}^T & 0
\end{pmatrix}\begin{pmatrix}
    w  \\
    s
\end{pmatrix} = \begin{pmatrix}
    v \\
    0
\end{pmatrix},
\end{align}
for unknown $(w,s) \in \C^n \times \R$. If $\bar{p}^Tv = 0$, then $s = 0$ and $w = A^{INV}_{i\omega_0}v$ will be a solution to the system $(i\omega_0I_n - A)w = v$ satisfying $\bar{p}^Tw = 0$. 

For the second Lyapunov coefficient $l_2 = \frac{1}{\omega_0}\Re{(c_2(0))}$ we need the $w^3$, $w^3\bar{w}$, $w^2\bar{w}^2$ and $w^3\bar{w}^2$ terms from the homological equation. The first three terms yield the equations
\begin{align}
    H_{3000} &= (3i\omega_0I_n - A)^{-1}[3B(q,H_{2000}) + C(q,q,q)],\\
    H_{3100} &= (2i\omega_0I_n - A)^{-1}[3B(q,H_{2100}) + B(\bar{q},H_{3000}) + 3B(H_{1100},H_{2000})\nonumber\\
    &+ 3C(q,q,H_{1100}) +  3C(q,\bar{q},H_{2000})+ D(q,q,q,\bar{q}) -6c_1(0)H_{2000}  ], \\
    H_{2200} &= -A^{-1}[2B(q,\overline{H}_{2100}) + 2B(\bar{q},H_{2100}) + B(H_{2000},\overline{H}_{2000}) \nonumber\\
    &+ 2B(H_{1100},H_{1100}) + C(q,q,\overline{H}_{2000}) + 4C(q,\bar{q},H_{1100})\nonumber\\
    &+ C(\bar{q},\bar{q},H_{2000}) + D(q,q,\bar{q},\bar{q}) ].
\end{align}
Note that we used that $c_1(0) +\bar{c}_1(0) = 2\omega_0l_1 = 0$ to simplify the expression for $H_{2200}$. Collecting the $w^3\bar{w}^2$ terms yields the equation
\begin{align}\label{eq:AH3200}
    (i\omega_0I_n - A)H_{3200} &= M_{3200}-[12c_2(0)q + 6i\Im{(c_1(0))}H_{2100}],
\end{align}
where
\begin{align*}
    M_{3200} &=3B(q,H_{2200}) + 2B(\bar{q},H_{3100})+ B(\overline{H}_{2000},H_{3000}) + 6B(H_{1100},H_{2100}) \nonumber\\
    &+ 3B(\overline{H}_{2100},H_{2000}) + 3C(q,q,\overline{H}_{2100})+ 6C(q,\bar{q},H_{2100})  \nonumber\\
    &+ 3C(q,H_{2000},\overline{H}_{2000})+ 6C(q,H_{1100},H_{1100}) + C(\bar{q},\bar{q},H_{3000})   \nonumber\\
    &+ 6C(\bar{q},H_{1100},H_{2000})+ D(q,q,q,\overline{H}_{2000})+ 6D(q,q,\bar{q},H_{1100})   \nonumber\\
    &+ 3D(q,\bar{q},\bar{q},H_{2000})+ E(q,q,q,\bar{q},\bar{q}).
\end{align*}
By applying the Fredholm solvability condition to equation \cref{eq:AH3200}, we find that
\begin{align}
    c_2(0) &= \frac{1}{12}\bar{p}^TM_{3200}.
\end{align}
For the higher order coefficients we will also need $H_{3200}$. The unique solution to equation $\cref{eq:AH3200}$ satisfying $\bar{p}^TH_{3200} = 0$ can again be found by solving the corresponding bordered system:
\begin{equation}
    H_{3200} = A_{i\omega_0}^{INV}(M_{3200}-[12c_2(0)q + 6i\Im{(c_1(0))}H_{2100}]).
\end{equation}
An expression of the seventh order coefficient $c_3(0)$ has been derived in \cite{sotomayor2007bifurcation}, using the same normalization method. For this we need the coefficients of the terms $w^4, w^4\bar{w}, w^4\bar{w}^2$ and $w^3\bar{w}^3$ from the homological equation. These yield regular systems and the solutions are respectively
\begin{align}
    H_{4000} &= (4i\omega_0I_n - A)^{-1}[4B(q,H_{3000}) + 3B(H_{2000},H_{2000}) \nonumber\\
    &+ 6C(q,q,H_{2000}) + D(q,q,q,q)] ,\\
    H_{4100} &= (3i\omega_0I_n - A)^{-1}[4B(q,H_{3100}) + B(\bar{q},H_{4000}) + 4B(H_{1100},H_{3000}) \nonumber\\
    &+ 6B(H_{2000},H_{2100}) + 6C(q,q,H_{2100}) + 4C(q,\bar{q},H_{3000}) \nonumber\\
    &+ 12C(q,H_{1100},H_{2000}) + 3C(\bar{q},H_{2000},H_{2000}) + 4D(q,q,q,H_{1100}) \nonumber\\
    &+ 6D(q,q,\bar{q},H_{2000}) + E(q,q,q,q,\bar{q}) - 12c_1(0)H_{3000}] ,\\
    H_{4200} &= (2i\omega_0I_n - A)^{-1}[4B(q,H_{3200}) + 2B(\bar{q},H_{4100}) \nonumber\\
    &+ B(\overline{H}_{2000},H_{4000}) + 8B(H_{1100},H_{3100}) + 4B(\overline{H}_{2100},H_{3000}) \nonumber\\
    &+ 6B(H_{2000},H_{2200}) + 6B(H_{2100},H_{2100}) + 6C(q,q,H_{2200}) \nonumber\\
    &+ 8C(q,\bar{q},H_{3100}) + 4C(q,\overline{H}_{2000},H_{3000}) + 24C(q,H_{1100},H_{2100}) \nonumber\\
    &+ 12C(q,\overline{H}_{2100},H_{2000}) + C(\bar{q},\bar{q},H_{4000}) + 8C(\bar{q},H_{1100},H_{3000})\nonumber\\
    &+ 12C(\bar{q},H_{2000},H_{2100}) + 3C(\overline{H}_{2000},H_{2000},H_{2000}) + 12C(H_{1100},H_{1100},H_{2000}) \nonumber\\
    &+ 4D(q,q,q,\overline{H}_{2100}) + 12D(q,q,\bar{q},H_{2100}) + 6D(q,q,\overline{H}_{2000},H_{2000}) \nonumber\\
    &+ 12D(q,q,H_{1100},H_{1100}) + 4D(q,\bar{q},\bar{q},H_{3000}) + 24D(q,\bar{q},H_{1100},H_{2000}) \nonumber\\
    &+ 3D(\bar{q},\bar{q},H_{2000},H_{2000}) + E(q,q,q,q,\overline{H}_{2000}) + 8E(q,q,q,\bar{q},H_{1100})\nonumber\\
    &+ 6E(q,q,\bar{q},\bar{q},H_{2000}) + K(q,q,q,q,\bar{q},\bar{q}) -8(6c_2(0)H_{2000} + 2c_1(0)H_{3100} )],\\
    H_{3300} &= -A^{-1}[3B(q,\overline{H}_{3200}) + 3B(\bar{q},H_{3200})+ 3B(\overline{H}_{2000},H_{3100}) + B(\overline{H}_{3000},H_{3000}) \nonumber\\
    &+ 9B(H_{1100},H_{2200}) + 9B(H_{2100},\overline{H}_{2100}) + 3B(\overline{H}_{3100},H_{2000}) \nonumber\\
    &+ 3C(q,q,\overline{H}_{3100}) + 9C(q,\bar{q},H_{2200}) + 9C(q,\overline{H}_{2000},H_{2100}) \nonumber\\
    &+ 3C(q,\overline{H}_{3000},H_{2000}) + 18C(q,H_{1100},\overline{H}_{2100}) + 3C(\bar{q},\bar{q},H_{3100}) \nonumber\\
    &+ 3C(\bar{q},\overline{H}_{2000},H_{3000}) + 18C(\bar{q},H_{1100},H_{2100}) + 9C(\bar{q},\overline{H}_{2100},H_{2000}) \nonumber\\
    &+ 9C(\overline{H}_{2000},H_{1100},H_{2000}) + 6C(H_{1100},H_{1100},H_{1100}) + D(q,q,q,\overline{H}_{3000}) \nonumber\\
    &+ 9D(q,q,\bar{q},\overline{H}_{2100}) + 9D(q,q,\overline{H}_{2000},H_{1100}) + 9D(q,\bar{q},\bar{q},H_{2100}) \nonumber\\
    &+ 9D(q,\bar{q},\overline{H}_{2000},H_{2000}) + 18D(q,\bar{q},H_{1100},H_{1100}) + D(\bar{q},\bar{q},\bar{q},H_{3000}) \nonumber\\
    &+ 9D(\bar{q},\bar{q},H_{1100},H_{2000}) + 3E(q,q,q,\bar{q},\overline{H}_{2000}) + 9E(q,q,\bar{q},\bar{q},H_{1100}) \nonumber\\
    &+ 3E(q,\bar{q},\bar{q},\bar{q},H_{2000}) + K(q,q,q,\bar{q},\bar{q},\bar{q})  -72d_2H_{1100} ].
\end{align}
Here we used that $c_1(0) + \overline{c}_1(0) = 0$ and $c_2(0) + \overline{c}_2(0) = \Re{}(c_2(0)) = d_2$ to simplify the last two expressions. Note that the term $D(\bar{q},\bar{q},\bar{q},H_{3000})$ in the expression for $H_{3300}$ is missing in \cite{sotomayor2007bifurcation}, although it was included in their calculations. Collecting the $w^4\bar{w}^3$ terms results in the equation
\begin{align}\label{eq:AH43}
    (i\omega_0I_n - A)H_{4300} &= M_{4300} - (144c_3(0)q + 72(2c_2(0) + \overline{c}_2(0))H_{2100} \nonumber\\
    &+ 12i\Im{}\{c_1(0) \}H_{3200}),
\end{align}
where we used that $3c_1(0) + 2\overline{c}_1(0) = i\Im{}\{c_1(0)\}$. The term $M_{4300}$ contains all the derivatives and is presented in the supplementary material (SM1.1).

To obtain the equation for $c_3(0)$ we apply the Fredholm solvability condition to the singular system $\cref{eq:AH43}$. This yields the equation
\begin{align*}
    c_3(0) &= \frac{1}{144}\bar{p}^TM_{4300},
\end{align*}
where we used that $\bar{p}^T H_{2100} = 0$ and $\bar{p}^TH_{3200}  = 0$. The unique solution to equation $\cref{eq:AH43}$ satisfying $\bar{p}^T H_{4300} = 0$ can again be found by solving the corresponding bordered system:
\begin{equation}
    H_{4300} = A_{i\omega_0}^{INV}(M_{4300}-[144c_3(0)q + 72(2c_2(0) + \overline{c}_2(0))H_{2100} + 12i\Im{}\{c_1(0) \}H_{3200}]).
\end{equation}

\subsection{Parameter-related coefficients}\label{sec:parameter_coef_ODE}

\paragraph{Linear coefficients $K_{10}$,$K_{01}$}
We first determine the coefficients for the linear approximation of the parameter transformation $K$. Collecting the $\beta_1$ and $\beta_2$ terms in $\cref{eq:HOM_GH_ODE}$ yields for $\mu = (10),(01)$ the systems
\begin{align}
    AH_{00\mu} = -J_1K_{\mu}.
\end{align}
Let $e_1,e_2 \in \R^2$ be the standard basis vectors. Then we can write
\begin{equation}\label{eq:Kmu}
    K_{\mu} = \gamma_{1,\mu}e_1 + \gamma_{2,\mu}e_2,
\end{equation}
for some yet unknown $\gamma_{1,\mu},\gamma_{2,\mu} \in \R$. Since $A$ is regular, we have
\begin{equation}\label{eq:H00mu}
    H_{00\mu} = -\gamma_{1,\mu}A^{-1}J_1e_1 - \gamma_{2,\mu}A^{-1}J_1e_2.
\end{equation}
In the equations that follow we will occasionally use the following notation 
\begin{equation}\label{eq:delta_mu}
    \delta_{\mu}^{ij} =  
    \begin{cases}
       1, \quad \textrm{if } \mu = (ij), \\
       0, \quad \textrm{if } \mu = (ji)
    \end{cases}, \quad \textrm{ for } i,j \in \N_0.
\end{equation}
The $\beta_1w$ and $\beta_2w$ terms yield the systems
\begin{align}\label{eq:AH10mu}
    (i\omega_0I_n - A)H_{10\mu} =  A_1(q,K_{\mu}) +B(q,H_{00\mu}) -(\delta_{\mu}^{10} + ib_{1,\mu})q.
\end{align}
To reduce the length of the equations, it is convenient to define 
\begin{equation}
    \Gamma_i(q) = A_1(q,e_i) + B(q,-A^{-1}J_1e_i).
\end{equation}
Substituting equations $\cref{eq:Kmu}$ and $\cref{eq:H00mu}$ into equation $\cref{eq:AH10mu}$ results in 
\begin{align}\label{eq:AH10mu2}
    (i\omega_0I_n - A)H_{10\mu} =  \gamma_{1,\mu}\Gamma_1(q) + \gamma_{2,\mu}\Gamma_2(q) - (\delta_{\mu}^{10} + ib_{1,\mu})q,
\end{align}
Applying the Fredholm solvability condition to equation $\cref{eq:AH10mu2}$ results in
\begin{equation}
    \delta_{\mu}^{10} + ib_{1,\mu} = \bar{p}^T[\gamma_{1,\mu}\Gamma_1(q) + \gamma_{2,\mu}\Gamma_2(q)].
\end{equation}
If we take the real and imaginary parts of the above equation, we find that
\begin{equation}\label{eq:deltamu1}
    \delta_{\mu}^{10} = \gamma_{1,\mu}\Re{[\bar{p}^T\Gamma_1(q)]} + \gamma_{2,\mu}\Re{[\bar{p}^T\Gamma_2(q)]},
\end{equation}
and
\begin{equation}\label{eq:b1mu1}
    b_{1,\mu} = \gamma_{1,\mu}\Im{[\bar{p}^T\Gamma_1(q)]} + \gamma_{2,\mu}\Im{[\bar{p}^T\Gamma_2(q)]}.
\end{equation}
 A solution $H_{10\mu}$ of equation $\cref{eq:AH10mu2}$ satisfying $\bar{p}^TH_{10\mu}= 0$ can be obtained by solving the bordered system
\begin{equation*}
    \begin{pmatrix}
    i\omega_0I_n - A & q \\
    \bar{p}^T & 0
\end{pmatrix}\begin{pmatrix}
    H_{10\mu}  \\
    s
\end{pmatrix} = \begin{pmatrix}
     \gamma_{1,\mu}\Gamma_1(q) + \gamma_{2,\mu}\Gamma_2(q) - (\delta_{\mu}^{10} + ib_{1,\mu})q \\
    0
\end{pmatrix}.
\end{equation*}
Since the coefficients $b_{1,\mu}$, $\gamma_{1,\mu}$ and $\gamma_{2,\mu}$ are still unknown, this system cannot be solved directly. Instead, we derive another equation similar to $\cref{eq:deltamu1}$, where the only unknown are $\gamma_{1,\mu}$ and $\gamma_{2,\mu}$. This will then allow us to set up a linear system to solve for these coefficients. To achieve this, we would like to apply the bordered inverse linearly to the separate terms. However, applied to the separate terms, the solution to the bordered system \cref{eq:borderedsys} may not satisfy $s = 0$, losing its original interpretation. We write $w = A_{i\omega_0}^{Inv}v$ when $w$ is obtained by solving system \cref{eq:borderedsys} without ensuring $s = 0$, i.e. $w = A_{i\omega_0}^{Inv}v$  does not necessarily solve $(i\omega_0I_n - A)w = v$. The solution of equation $\cref{eq:AH10mu2}$ can then be written as
\begin{equation}\label{eq:H10mu}
    H_{10\mu} =  \gamma_{1,\mu}A^{Inv}_{i\omega_0}\Gamma_1(q) +  \gamma_{2,\mu}A^{Inv}_{i\omega_0}\Gamma_2(q) -(\delta_{\mu}^{10} + ib_{1,\mu})A^{Inv}_{i\omega_0}q.
\end{equation}
Collecting the $w^2\beta_i$ and $w\bar{w}\beta_i$ terms respectively yield the systems
\begin{align*}
     (2i\omega_0I_n - A)H_{20\mu} &= A_1(H_{2000},K_{\mu}) +2B(q,H_{10\mu})+ B(H_{00\mu},H_{2000})  \\
    &+ B_1(q,q,K_{\mu}) + C(q,q,H_{00\mu}) - 2(\delta_{\mu}^{10} + ib_{1,\mu})H_{2000} , \\
    -AH_{11\mu} &=  A_1(H_{1100},K_{\mu}) + 2\Re{\{B(\bar{q},H_{10\mu})\}} + B(H_{00\mu},H_{1100}) \\
    &+ B_1(q,\bar{q},K_{\mu}) + C(q,\bar{q},H_{00\mu})
    -2\delta_{\mu}^{10}H_{1100}.
\end{align*}
Both these systems are regular with solutions
\begin{align}
     H_{20\mu} &=  A_{2i\omega_0}^{-1}[A_1(H_{2000},K_{\mu}) +2B(q,H_{10\mu})+ B(H_{00\mu},H_{2000}) \nonumber\\
    &+ B_1(q,q,K_{\mu}) + C(q,q,H_{00\mu})- 2(\delta_{\mu}^{10} + ib_{1,\mu})H_{2000}] , \label{eq:H20mu1}\\
    H_{11\mu} &= -A^{-1}[A_1(H_{1100},K_{\mu}) + 2\Re{(B(\bar{q},H_{10\mu}))} + B(H_{00\mu},H_{1100}) \nonumber\\
    &+B_1(q,\bar{q},K_{\mu}) + C(q,\bar{q},H_{00\mu}) -2\delta_{\mu}^{10}H_{1100}],\label{eq:H11mu1}
\end{align}
where we write $A_{2i\omega_0}^{-1} = (2i\omega_0I_n - A)^{-1}$. We will now substitute equations $\cref{eq:Kmu}$, $\cref{eq:H00mu}$ and $\cref{eq:H10mu}$ into the above expressions. For this, it is convenient to define for $i = 1,2$ the following functions:
\begin{align}
\Lambda_i(u,v,w) &= \Gamma_i(u) + 2B(v,A_{i\omega_0}^{Inv}\Gamma_i(w)) + B_1(v,w,e_i) + C(v,w,-A^{-1}J_1e_i), \\
\Pi_i(u,v,w) &= \Gamma_i(u) + 2\Re{\{B\left(v,A_{i\omega_0}^{Inv}\Gamma_i(w)\right)\}} + B_1(v,w,e_i) + C(v,w,-A^{-1}J_1e_i).
\end{align}
Then the equations for $H_{20\mu}$ and $H_{11\mu}$ become
\begin{align*}
     H_{20\mu} &= \gamma_{1,\mu}A_{2i\omega_0}^{-1}\Lambda_1(H_{2000},q,q) + \gamma_{2,\mu}A_{2i\omega_0}^{-1}\Lambda_2(H_{2000},q,q) \\
     &-2(\delta_{\mu}^{10} + ib_{1,\mu})A_{2i\omega_0}^{-1}\left(H_{2000} + B(q, A_{i\omega_0}^{Inv}q)\right), \\
    H_{11\mu} &=  -\gamma_{1,\mu}A^{-1}\Pi_1(H_{1100},\bar{q},q)  -\gamma_{2,\mu}A^{-1}\Pi_2(H_{1100},\bar{q},q)  \\
    &+ 2\delta_{\mu}^{10}A^{-1}\left(H_{1100} + \Re{\left\{B(\bar{q}, A_{i\omega_0}^{Inv}q)\right\}}\right) + 2b_{1,\mu}A^{-1}\Re{\left\{iB(\bar{q}, A_{i\omega_0}^{Inv}q)\right\}}.
\end{align*}
The $w^2\bar{w}\beta_i$ terms yield the systems
\begin{align}\label{eq:AH21mu}
    (i\omega_0I_n - A)H_{21\mu} &=  A_1(H_{2100},K_{\mu}) + 2B(q,H_{11\mu}) + B(\bar{q},H_{20\mu}) + B(H_{00\mu},H_{2100}) \nonumber\\
    &+ B(\overline{H}_{10\mu},H_{2000}) + 2B(H_{10\mu},H_{1100}) + 2B_1(q,H_{1100},K_{\mu})  \nonumber \\
    &+ B_1(\bar{q},H_{2000},K_{\mu})+ C(q,q,\overline{H}_{10\mu}) + 2C(q,\bar{q}.H_{10\mu}) + 2C(q,H_{00\mu},H_{1100})  \nonumber\\
    &+ C(\bar{q},H_{00\mu},H_{2000})+ C_1(q,q,\bar{q},K_{\mu}) + D(q,q,\bar{q},H_{00\mu}) -[2(\delta_{\mu}^{01} + ib_{2,\mu})q \nonumber\\
    &+ (3\delta_{\mu}^{10} + ib_{1,\mu})H_{2100} + 2c_1(0)H_{10\mu}].
\end{align}
Applying the Fredholm solvability equation to equation $\cref{eq:AH21mu}$ results in
\begin{align}\label{eq:H21mu_Fred}
    \delta_{\mu}^{01} + ib_{2,\mu} &= \frac{1}{2}\bar{p}^T [A_1(H_{2100},K_{\mu}) + 2B(q,H_{11\mu}) + B(\bar{q},H_{20\mu}) + B(H_{00\mu},H_{2100}) \nonumber\\
    &+ B(\overline{H}_{10\mu},H_{2000}) + 2B(H_{10\mu},H_{1100}) + 2B_1(q,H_{1100},K_{\mu}) + B_1(\bar{q},H_{2000},K_{\mu}) \nonumber\\
    &+ C(q,q,\overline{H}_{10\mu}) + 2C(q,\bar{q}.H_{10\mu}) + 2C(q,H_{00\mu},H_{1100}) + C(\bar{q},H_{00\mu},H_{2000}) \nonumber\\
    &+ C_1(q,q,\bar{q},K_{\mu}) + D(q,q,\bar{q},H_{00\mu})],
\end{align}
where we used that $\bar{p}^TH_{2100}  = 0$ and $\bar{p}^TH_{10\mu} = 0$. If we substitute the expressions for $H_{00\mu}, H_{10\mu},H_{20\mu}, H_{11\mu}$, $K_{\mu}$ and $b_{1,\mu}$ into equation $\cref{eq:H21mu_Fred}$, we can finally solve for $\gamma_{1,\mu}$ and $\gamma_{2,\mu}$. After careful rewriting, this results in the following linear system of equations:
\begin{equation}\label{eq:PQ1}
    P\begin{pmatrix}
        \gamma_{1,\mu} \\
        \gamma_{2,\mu}
    \end{pmatrix} = Q_{\mu},
\end{equation}
where $P \in \R^{2\times2}$ is given by
\begin{align*}
    P_{1k} &= \Re{[\bar{p}^T\Gamma_k(q)]}, \\
    P_{2k} &= \frac{1}{2}\Re{}\Big\{ \bar{p}^T\Bigl[ \Gamma_k(H_{2100}) + 2B(q,-A^{-1}\Pi_k(H_{1100},\bar{q},q)) + B(\bar{q},A^{-1}_{2i\omega_0}\Lambda_k(H_{2000},q,q)) \nonumber\\
    &  +B(H_{2000},\overline{A_{i\omega_0}^{Inv}\Gamma_k(q)}) +2B(H_{1100},A_{i\omega_0}^{Inv}\Gamma_k(q)) + 2B_1(q,H_{1100},e_k)+ B_1(\bar{q},H_{2000},e_k) \nonumber\\
    &  +C(q,q,\overline{A_{i\omega_0}^{Inv}\Gamma_k(q)}) + 2C(q,\bar{q},A_{i\omega_0}^{Inv}\Gamma_k(q))+2C(q,H_{1100},-A^{-1}J_1e_k) \nonumber\\
    &  +C(\bar{q},H_{2000},-A^{-1}J_1e_k) + C_1(q,q,\bar{q},e_k) + D(q,q,\bar{q},-A^{-1}J_1e_k)  \nonumber\\
    &+\Im{[\bar{p}^T\Gamma_k(q)]} \Big(-4B\left(q, -A^{-1}\Re{\left\{iB(\bar{q}, A_{i\omega_0}^{Inv}q)\right\}}\right) - 2i B\left(\bar{q},A^{-1}_{2i\omega_0}\left(H_{2000} +B(q, A_{i\omega_0}^{Inv}q) \right)\right)\nonumber\\
    & + iB\left(H_{2000}, \overline{A_{i\omega_0}^{Inv}q}\right) -2iB\left(H_{1100},A_{i\omega_0}^{Inv}q\right) +iC\left(q,q, \overline{A_{i\omega_0}^{Inv}q}\right) - 2iC\left(q,\bar{q}, A_{i\omega_0}^{Inv}q \right)
    \Big)\Big] \Big\},  
\end{align*}
for $k = 1,2$ and $Q_{\mu} \in \R^2$, $\mu = (10),(01)$ is given by
\begin{align*}
   Q_{1,\mu} &= \delta_{\mu}^{10}, \\
    Q_{2,\mu} &= \delta_{\mu}^{01} +\frac{1}{2}\delta_{\mu}^{10} \Re{}\Big\{ \bar{p}^T\Bigl[ 4B\left(q,-A^{-1}\left(H_{1100} + \Re{\left\{B(\bar{q}, A_{i\omega_0}^{Inv}q)\right\}}\right) \right) \nonumber\\ 
    &+2B\left(\bar{q},A^{-1}_{2i\omega_0}\left(H_{2000} +B(q, A_{i\omega_0}^{Inv}q) \right)\right) + B\left(H_{2000}, \overline{A_{i\omega_0}^{Inv}q}\right) \nonumber\\
    & +2B\left(H_{1100},A_{i\omega_0}^{Inv}q\right) +C\left(q,q, \overline{A_{i\omega_0}^{Inv}q}\right) + 2C\left(q,\bar{q}, A_{i\omega_0}^{Inv}q \right)\Bigr]\Big\}.
\end{align*}
Thus, one first computes the linear coefficients $K_{10}$ and $K_{01}$ by solving the system $\cref{eq:PQ1}$. This system should be solvable if the transversality condition \cref{eq:trans_cond} holds. Once $K_{10}, K_{01}$ are known, we can determine the coefficients $H_{00\mu}$ from equation $\cref{eq:H00mu}$ and $b_{1,\mu}$ from equation $\cref{eq:b1mu1}$ for $\mu = (10),(01)$. With these coefficients we can determine $H_{10\mu}$ from equation $\cref{eq:H10mu}$, and $H_{20\mu}$ and $H_{11\mu}$ with $\mu = (10),(01)$ from respectively equations $\cref{eq:H20mu1}$ and $\cref{eq:H11mu1}$. The coefficient $b_{2,01}$ is then given by the imaginary part of the right-hand side of equation $\cref{eq:H21mu_Fred}$. Finally, we can use equation $\cref{eq:AH21mu}$ to solve for $H_{2101}$ satisfying $\bar{p}^TH_{2101} = 0$ by solving the corresponding bordered system. 

As we continue to higher-order coefficients, the part containing all the multilinear forms in the equations will increase in length. To keep things readable, we will not write all of these terms here. Instead, we use the notation $M_{ijkl}$ for parts containing only a linear combination of multilinear forms depending on known coefficients. The corresponding terms are presented in the supplementary material (S1.1). 

\paragraph{The coefficient $a_{3201}$}
 Note that $a_{3201} = \Re{}(g_{3201})$. To determine this coefficient we also need the coefficients $H_{3001}, H_{3101}$ and $H_{2201}$. These can be found by collecting the $w^3\beta_2$, $w^3\bar{w}\beta_2$ and $w^2\bar{w}^2\beta_2$ terms from the homological equation $\cref{eq:HOM_GH_ODE}$. This results in the following equations
 \begin{align*}
     H_{3001} &=  (3i\omega_0I_n -A)^{-1} [M_{3001} - 3ib_{1,01}H_{3000}], \nonumber\\
    H_{3101} &= (2i\omega_0I_n - A)^{-1}[M_{3101}-6(1 + ib_{2,01})H_{2000} - 6c_1(0)H_{2001} - 2ib_{1,01}H_{3100} ], \\
    H_{2201} &= -A^{-1}[M_{2201}  - 8 H_{1100} ].
 \end{align*}
 The coefficient $g_{3201}$ can now be found by applying the Fredholm alternative to the equation that follows from collecting the $w^3\bar{w}^2\beta_2$ terms. This yields the following equation
 \begin{align}\label{eq:AH3201}
     (i\omega_0I_n - A)H_{3201} &= M_{3201} - [12g_{3201}q + 12c_2(0)H_{1001}  +(18 + 6ib_{2,01})H_{2100} \nonumber\\
     &+ 6i\Im{}\{c_1(0)\}H_{2101} + ib_{1,01}H_{3200}].
 \end{align}
 Applying the Fredholm alternative to $\cref{eq:AH3201}$ now yields
 \begin{align}
     g_{3201} &=  \frac{1}{12}\bar{p}^TM_{3201}.
 \end{align}
 Now the unique solution for $H_{3201}$ satisfying $\bar{p}^TH_{3201}  = 0$ can be solved from equation $\cref{eq:AH3201}$ using the corresponding bordered system, i.e.
 \begin{align}
     H_{3201} &= A_{i\omega_0}^{INV}\left(M_{3201} - [12g_{3201}q + 12c_2(0)H_{1001}  +(18 + 6ib_{2,01})H_{2100} \right.\nonumber\\
     &\left.+ 6i\Im{}\{c_1(0)\}H_{2101} + ib_{1,01}H_{3200}] \right).
 \end{align}
 
\paragraph{The quadratic coefficients $K_{02},K_{11}$}
The quadratic coefficients are determined similarly. Collecting the $\beta_2^2$ and $\beta_1\beta_2$ terms from $\cref{eq:HOM_GH_ODE}$ yields the equations
\begin{align}\label{eq:AH0002}
    AH_{00\mu} &= -J_1K_{\mu} - M_{00\mu} \quad \mu = (02), (11), 
\end{align}
Let $e_1,e_2 \in \R^2$ be the standard basis vectors and write
\begin{align}\label{eq:K02}
    K_{\mu} = \gamma_{1,\mu}e_1 + \gamma_{2,\mu}e_2,
\end{align}
where $\gamma_{1,\mu}, \gamma_{1,\mu} \in \R$ are unknown constants that need to be determined. From equation $\cref{eq:AH0002}$ it follows that
\begin{equation}\label{eq:H0002}
    H_{00\mu} = -\gamma_{1,\mu}A^{-1}J_1e_1 - \gamma_{2,\mu}A^{-1}J_1e_2 - A^{-1}M_{00\mu}, 
\end{equation}
The $w\beta_2^2$ and $w\beta_1\beta_2$ terms yield the equations
\begin{align}\label{eq:AH1002}
    (i\omega_0I_n- A)H_{10\mu} &= A_1(q,K_{\mu}) +  B(q,H_{00\mu}) + r_{10\mu} - ib_{1,\mu}q  ,
\end{align}
where we define
\begin{align*}
    r_{1002} &=  M_{1002} - 2ib_{1,01}H_{1001},\\
    r_{1011} &=  M_{1011} - [(1+ib_{1,10})H_{1001} + ib_{1,01}H_{1010}].
\end{align*}
Applying the Fredholm alternative to equation $\cref{eq:AH1002}$ yields the equations
\begin{align}
    b_{1,\mu}i &= \bar{p}^T[A_1(q,K_{\mu}) +  B(q,H_{00\mu}) + M_{10\mu}].
\end{align}
Here we used that $\bar{p}^T H_{1010} = \bar{p}^TH_{1001} = 0$. Substituting expressions $\cref{eq:K02}$ and $\cref{eq:H0002}$ into the above equation results in 
\begin{align}\label{eq:ib102}
    b_{1,\mu}i &= \bar{p}^T[\gamma_{1,\mu}\Gamma_1(q) + \gamma_{2,\mu}\Gamma_2(q)  +  \tilde{M}_{10\mu}],
\end{align}
where 
\begin{align}\label{eq:Mtilde10mu}
    \tilde{M}_{10\mu} &= M_{10\mu} +  B(q,-A^{-1}M_{00\mu}) .
\end{align}
From equation $\cref{eq:ib102}$, it follows that 
\begin{align}
     \gamma_{1,\mu}\Re{[\bar{p}^T\Gamma_1(q)]} + \gamma_{2,\mu}\Re{[\bar{p}^T\Gamma_2(q)]} = -\Re{}\{\bar{p}^T \tilde{M}_{10\mu}\},
\end{align}
and
\begin{align}\label{eq:b102}
    b_{1,\mu} = \gamma_{1,\mu}\Im{[\bar{p}^T\Gamma_1(q)]} + \gamma_{2,\mu}\Im{[\bar{p}^T\Gamma_2(q)]} + \Im{}\{\bar{p}^T \tilde{M}_{10\mu}\}.
\end{align}
Furthermore, $H_{10\mu}$ is given by 
\begin{align}\label{eq:H1002}
    H_{10\mu} =  \gamma_{1,\mu}A_{i\omega_0}^{Inv}\Gamma_1(q)  + \gamma_{2,\mu}A_{i\omega_0}^{Inv}\Gamma_2(q) + A_{i\omega_0}^{Inv}\tilde{r}_{10\mu}  - ib_{1,\mu}A_{i\omega_0}^{Inv}q ,
\end{align}
where
\begin{align}\label{eq:rtilde10mu}
    \tilde{r}_{10\mu} &= r_{10\mu} + B(q, -A^{-1}M_{00\mu}).
\end{align}
Now collect the $w^2\beta^{\mu}$ terms and $w\bar{w}\beta^{\mu}$ terms for $\mu = (02),(11)$. These terms yield the equations
\begin{align}
    (2i\omega_0I_n- A)H_{20\mu} &= A_1(H_{2000},K_{\mu}) + 2B(q,H_{10\mu})+ B(H_{00\mu},H_{2000}) \nonumber \\\label{eq:AH2002}
    &+  B_1(q,q,K_{\mu})+ C(q,q,H_{00\mu})  +r_{20\mu}- 2ib_{1,\mu}H_{2000}, \\
    -AH_{11\mu} &= A_1(H_{1100},K_{\mu}) + 2\Re{}\left(B(\bar{q},H_{10\mu})\right)+ B(H_{00\mu},H_{1100}) \nonumber\\
    &+ B_1(q,\bar{q},K_{\mu})+ C(q,\bar{q},H_{00\mu})+ r_{11\mu}, \label{eq:AH1102}
\end{align}
where
\begin{align*}
    r_{2002} &= M_{2002}- 4ib_{1,01}H_{2001}, \\
    r_{2011} &= M_{2011}- [2(1 + ib_{1,10})H_{2001} + 2ib_{1,01}H_{2010}], \\
    r_{1102} &= M_{1102} ,\quad \textrm{ and } \quad r_{1111} = M_{1111}- 2H_{1101}.
\end{align*}
Substituting equations $\cref{eq:K02}$, $\cref{eq:H0002}$ and $\cref{eq:H1002}$ into the expressions for $H_{20\mu}$ and $H_{11\mu}$ and solving for the coefficients yields the equations
\begin{align}
    H_{20\mu} &= \gamma_{1,\mu}A_{2i\omega_0}^{-1}\Lambda_1(H_{2000},q,q) + \gamma_{2,\mu}A_{2i\omega_0}^{-1}\Lambda_2(H_{2000},q,q) \nonumber\\
     &+ A^{-1}_{2i\omega_0}\tilde{r}_{20\mu}-2ib_{1,\mu}A_{2i\omega_0}^{-1}\left(H_{2000} + B(q, A_{i\omega_0}^{Inv}q)\right), \label{eq:H20mu2}\\
    H_{11\mu} &=  -\gamma_{1,\mu}A^{-1}\Pi_1(H_{1100},\bar{q},q)  -\gamma_{2,\mu}A^{-1}\Pi_2(H_{1100},\bar{q},q)  \nonumber\\
    &-A^{-1}\tilde{r}_{11\mu}+ 2b_{1,\mu}A^{-1}\Re{\left\{iB(\bar{q}, A_{i\omega_0}^{Inv}q)\right\}}, \label{eq:H11mu2}
\end{align}
where we define
\begin{align*}
    \tilde{r}_{20\mu} &= r_{20\mu} + 2B\left(q,A_{i\omega_0}^{Inv}\tilde{r}_{10\mu}\right) + B(H_{2000},-A^{-1}M_{00\mu}) + C(q,q,-A^{-1}M_{00\mu}),\\
    \tilde{r}_{11\mu} &=r_{11\mu} + 2\Re{}\left(B\left(\bar{q},A_{i\omega_0}^{Inv}\tilde{r}_{10\mu}\right)\right) + B(H_{1100},-A^{-1}M_{00\mu}) + C(q,\bar{q},-A^{-1}M_{00\mu}).
\end{align*}
The $w^2\bar{w}\beta^{\mu}$ terms yield for $\mu = (02),(11)$ the equations
\begin{align}\label{eq:H2102}
    (i\omega_0I_n- A)H_{21\mu} &= A_1(H_{2100},K_{\mu}) + 2B(q,H_{11\mu})+ B(\bar{q},H_{20\mu})+ B(H_{00\mu},H_{2100}) \nonumber\\
    &+ B(\overline{H}_{10\mu},H_{2000})+ 2B(H_{10\mu},H_{1100})+ 2B_1(q,H_{1100},K_{\mu}) \nonumber\\
    &+ B_1(\bar{q},H_{2000},K_{\mu})+ C(q,q,\overline{H}_{10\mu})+ 2C(q,\bar{q},H_{10\mu})+ 2C(q,H_{00\mu},H_{1100}) \nonumber\\
    &+ C(\bar{q},H_{00\mu},H_{2000}) + C_1(q,q,\bar{q},K_{\mu})+ D(q,q,\bar{q},H_{00\mu}) + r_{21\mu} \nonumber\\
    &- (2ib_{2,\mu}q +ib_{1,\mu}H_{2100} + 2c_1(0)H_{10\mu}) ,
\end{align}
where
\begin{align*}
    r_{2102} &= M_{2102}- [4(1 + ib_{2,01})H_{1001}  + 2ib_{1,01}H_{2101}], \\
    r_{2111} &= M_{2111}- [2ib_{2,10}H_{1001}+2 (1+ib_{2,01})H_{1010}+ \left(3 + ib_{1,10}\right) H_{2101}+ib_{1,01}H_{2110}].
\end{align*}
Applying the Fredholm solvability condition to $\cref{eq:H2102}$ yields the equations
\begin{align}\label{eq:b2mu2}
    b_{2,\mu}i &= \frac{1}{2}\bar{p}^T[A_1(H_{2100},K_{\mu}) + 2B(q,H_{11\mu})+ B(\bar{q},H_{20\mu})+ B(H_{00\mu},H_{2100}) \nonumber\\
    &+ B(\overline{H}_{10\mu},H_{2000})+ 2B(H_{10\mu},H_{1100})+ 2B_1(q,H_{1100},K_{\mu}) \nonumber\\
    &+ B_1(\bar{q},H_{2000},K_{\mu})+ C(q,q,\overline{H}_{10\mu})+ 2C(q,\bar{q},H_{10\mu})+ 2C(q,H_{00\mu},H_{1100}) \nonumber\\
    &+ C(\bar{q},H_{00\mu},H_{2000}) + C_1(q,q,\bar{q},K_{\mu})+ D(q,q,\bar{q},H_{00\mu}) + M_{21\mu} ].
\end{align}
Now substitute the expressions for $H_{00\mu}, H_{10\mu},H_{20\mu}, H_{11\mu}$, $K_{\mu}$ and $b_{1,\mu}$ into the above expression to solve for $\gamma_{1,\mu}$ and $\gamma_{2,\mu}$. The resulting system of equations for $\gamma_{1,\mu}$ and $\gamma_{2,\mu}$ with $\mu = (02),(11)$ is given by \cref{eq:PQ1}, where $Q_{\mu} \in \R^2$ is given by
 \begin{align}
     Q_{\mu, 1} &= -\Re{}\{\bar{p}^T \tilde{M}_{10\mu}\}, \label{eq:Qmu21}\\
     Q_{\mu, 2} &= -\frac{1}{2}\Re{}\Big\{ \bar{p}^T \Bigl[  2B\left(q,-A^{-1}\tilde{r}_{11\mu}\right) + B\left(\bar{q},A_{2i\omega_0}^{-1}\tilde{r}_{20\mu}\right)+ B\left(H_{2100},-A^{-1}M_{00\mu}\right) \nonumber\\
     &+ B\left(H_{2000},\overline{A_{i\omega_0}^{Inv}\tilde{r}_{10\mu}}\right) +2B\left(H_{1100},A_{i\omega_0}^{Inv}\tilde{r}_{10\mu}\right) + C\left(q,q, \overline{A_{i\omega_0}^{Inv}\tilde{r}_{10\mu}}\right)\nonumber\\
     & + 2C\left(q,\bar{q}, A_{i\omega_0}^{Inv}\tilde{r}_{10\mu}\right) + 2C\left(q,H_{1100},-A^{-1}M_{00\mu} \right)+ C\left(\bar{q},H_{2000},-A^{-1}M_{00\mu}\right) \nonumber\\
     & +D(q,q,\bar{q},-A^{-1}M_{00\mu}) + M_{21\mu} +\Im{}\{\bar{p}^T \tilde{M}_{10\mu}\}\Big(-4B\left(q, -A^{-1}\Re{\left\{iB(\bar{q}, A_{i\omega_0}^{Inv}q)\right\}}\right)\nonumber\\
     & - 2i B\left(\bar{q},A^{-1}_{2i\omega_0}\left(H_{2000} +B(q, A_{i\omega_0}^{Inv}q) \right)\right)+ iB\left(H_{2000}, \overline{A_{i\omega_0}^{Inv}q}\right) -2iB\left(H_{1100},A_{i\omega_0}^{Inv}q\right)\nonumber\nonumber\\
    &  +iC\left(q,q, \overline{A_{i\omega_0}^{Inv}q}\right) - 2iC\left(q,\bar{q}, A_{i\omega_0}^{Inv}q \right)
    \Big)\Bigr]\Bigr\}. \label{eq:Qmu22}
 \end{align}
 Thus, to obtain the quadratic coefficients $K_{02}$ and $K_{11}$ we solve system $\cref{eq:PQ1}$ for $\mu = (02),(11)$. Once these coefficients are known, one can calculate the coefficients $H_{00\mu}$,$b_{1,02}$ and $H_{10\mu}$ from respectively the equations $\cref{eq:H0002}$,$\cref{eq:b102}$ and $\cref{eq:H1002}$. Then, the coefficients $H_{2002}$ and $H_{1102}$ are calculated from equations $\cref{eq:AH2002},\cref{eq:AH1102}$. Finally,  the coefficient $H_{2102}$ satisfying $\bar{p}^T H_{2102}  = 0$, can be determined from system $\cref{eq:H2102}$ using the corresponding bordered system.
 \paragraph{The cubic coefficient $K_{03}$} When collecting the equations to solve for the cubic coefficient $K_{03}$ one will notice that they can be written in the same way as we did for the quadratic coefficients. Only the terms $M_{ijkl}$ and $r_{ijkl}$ will change. Define
 \begin{align*}
    r_{1003} &= M_{1003}- (3ib_{1,02}H_{1001} + 3ib_{1,01}H_{1002}), \\
    r_{2003} &=  M_{2003}  -(6ib_{1,02}H_{2001} + 6ib_{1,01}H_{2002}), \\
    r_{1103} &= M_{1103}.
 \end{align*}
 With these definitions, the equations become the same as equations $\cref{eq:AH0002}$,$\cref{eq:AH1002}$,$\cref{eq:AH2002}$, $\cref{eq:AH1102}$ and $\cref{eq:b2mu2}$ for $\mu = (03)$. As a result, we can solve for the coefficients $\gamma_{1,03}$ and $\gamma_{2,03}$ by using system $\cref{eq:PQ1}$ where the components of $Q_{\mu}$ are given by \cref{eq:Qmu21} and \cref{eq:Qmu22} for $\mu = (03)$. Then $H_{0003}$ and $H_{1003}$ can be calculated from respectively the equations $\cref{eq:AH0002}$ and $\cref{eq:AH1002}$.

 \paragraph{Remark} The equations for the parameter-dependent coefficients can be simplified by noting that $A^{Inv}_{i\omega_0}q = 0$. To see this, write the corresponding bordered system as the following system of equations
\begin{align*}
    (i\omega_0I_n- A)w + qs &= q, \\
    \bar{p}^T w = 0.
\end{align*}
Taking the inner product with $p$ in the first equation yields $s = 1$. Then the first equation becomes $(i\omega_0I_n- A)w = 0$, with solution $w = \lambda q$. Filling this into the second equation yields $\bar{p}^T w  = \lambda \bar{p}^T q  = 0$ from which it follows that $\lambda = 0$. Thus, the unique solution to this system is $(w, s) = (0,1)$. As a result, all terms containing $A^{Inv}_{i\omega_0}q$ vanish. Nevertheless, we included them here as this will make the conversion to DDEs in the following section easier.

 \section{Coefficients of the parameter-dependent normal form and the predictor in DDEs}\label{sec:coeff_DDE}
Assume that system $\cref{eq:Classical_DDE}$ has an equilibrium at the origin at $\alpha = (0,0) \in \R^2$ with only one pair of purely imaginary simple eigenvalues
\[
\lambda_{1,2} = \pm i\omega_0, \quad \omega_0 > 0,
\]
and no other eigenvalues on the imaginary axis. Furthermore, we assume that the first Lyapunov coefficient $l_1(0) = 0$ and the second Lyapunov coefficient $l_2(0) \neq 0$.
We have eigenfunctions $\varphi$ and $\varphi^{\odot}$ satisfying
\begin{equation*}
    A\varphi = i\omega_0\varphi, \quad A^{\star}\varphi^{\odot} = i\omega_0\varphi^{\odot}, \quad \langle \varphi^{\odot}, \varphi \rangle = 1,
\end{equation*}
where $\langle \varphi^{\odot}, \varphi \rangle$ is given by the pairing $\cref{eq:xsun_dual}$.
Furthermore, introduce $q, p \in \C^n$ such that
\[
\Delta(i\omega_0)q = 0, \quad p^T\Delta(i\omega_0) = 0, \quad p^T\Delta'(i\omega_0)q = 1.
\]
With $q$ and $p$ as above, explicit expressions for the eigenfunctions $\varphi$ and $\varphi^{\odot}$ are respectively given by equations $\cref{eq:eigen_phi}$ and $\cref{eq:eigen_phi_sun}$. Points $y \in X_0$ of the real critical eigenspace can be represented in terms of the complex coordinate $z = \langle \varphi^{\odot}, y \rangle$ as,
\[
y = z \varphi + \bar{z}\bar{\varphi}.
\]
The homological equation $\cref{eq:HOM_DDE}$ becomes
\begin{align}\label{eq:HOM_GH_DDE}
    A^{\odot\star}jH(z,\bar{z},\beta) + J_1K(\beta)r^{\odot\star} + &R(H(z,\bar{z},\beta),K(\beta))\nonumber\\
    &= j(D_zH(z,\bar{z},\beta)\dot{z} + D_{\bar{z}}H(z,\bar{z},\beta)\dot{\bar{z}}).
\end{align}
Restricted to the two-dimensional center manifold, the system can be transformed into the same normalform $\cref{eq:wdot}$ with new unfolding parameters $\beta = (\beta_1,\beta_2)$ if the transversality condition $\cref{eq:trans_cond}$ holds. Thus, we take the same truncated normal form:
\begin{align}\label{eq:zdot}
    \dot{z} &= (i\omega_0 + \beta_1 +ib_1(\beta))z + (\beta_2 +  ib_2(\beta))z|z|^2 + (c_2(0) + g_{3201}\beta_2)z|z|^4 \nonumber\\
    &+ c_3(0)z|z|^6 , 
\end{align}
where $b_1(\beta)$ and $\beta_2(\beta)$ are expanded as $\cref{eq:b1_expansion}$ and $\cref{eq:b2_expansion}$.
The nonlinearity in the homological equation $\cref{eq:HOM_GH_DDE}$ is expanded as 
\begin{align*}
    R(u,\alpha) &=  [A_1(u,\alpha) + \frac{1}{2}B(u,u) + \frac{1}{2}J_2(\alpha,\alpha) + \frac{1}{6}C(u,u,u) + \frac{1}{2}B_1(u,u,\alpha) \\
    &+ \frac{1}{2}A_2(u,\alpha,\alpha) + \frac{1}{6}J_3(\alpha,\alpha,\alpha) +  \frac{1}{24}D(u,u,u,u) + \frac{1}{6}C_1(u,u,u,\alpha) + \frac{1}{4}B_2(u,u,\alpha,\alpha) \\
    &+ \frac{1}{6}A_3(u,\alpha,\alpha,\alpha) + \frac{1}{120}E(u,u,u,u,u) + \frac{1}{24}D_1(u,u,u,u,\alpha) + \frac{1}{12}C_2(u,u,u,\alpha,\alpha) \\
    &+ \frac{1}{12}B_3(u,u,\alpha,\alpha,\alpha)+\frac{1}{720}K(u,u,u,u,u,u) +  \frac{1}{120}E_1(u,u,u,u,u,\alpha) \\
    &+ \frac{1}{36}C_3(u,u,u,\alpha,\alpha,\alpha) +\frac{1}{5040}L(u,u,u,u,u,u,u) ]r^{\odot\star}.
\end{align*}
Finally, $H$ is expanded as
\begin{align}\label{eq:H_GH_DDE}
H(z,\bar{z},\beta)& = z\varphi + \bar{z}\bar{\varphi} + \sum_{n + m = 2}^7 \frac{1}{n!m!}H_{nm00}z^n\bar{z}^m + \sum_{n + m = 0}^5 H_{nm01}\frac{1}{n!m!}z^n\bar{z}^m\beta_2 \nonumber \\
&+ \sum_{n + m=0}^3 \frac{1}{n!m!}H_{nm10}z^n\bar{z}^m\beta_1 + \sum_{n + m=0}^3 \frac{1}{2n!m!}H_{nm02}z^n\bar{z}^m\beta_2^2 \nonumber\\
&+ \sum_{n + m=0}^1 \frac{1}{n!m!}H_{nm11}z^n\bar{z}^m\beta_1 \beta_2 + \sum_{n + m=0}^1 \frac{1}{6n!m!}H_{nm03}z^n\bar{z}^m\beta_2^3 \nonumber\\
&\textcolor{gray}{+ \frac{1}{6}H_{1103}z\bar{z}\beta_2^3 + \frac{1}{12}H_{2003}z^2\beta_2^3 + \frac{1}{12}H_{2103}z^2\bar{z}
\beta_2^3 }.
\end{align}
Note that $H$ is now a mapping into $X = C([-h,0],\R^n)$. Since $X$ is real, the property $H_{ijkl} = \overline{H}_{jikl}$ still holds. Additionally, the parameter relation $K$ remains the same as in the ODE case, given by $\cref{eq:K}$. 

The equations collected from the homological equation $\cref{eq:HOM_GH_DDE}$ have a similar structure as those in the ODE case. To obtain the corresponding equations for DDEs, one substitutes $A$ with $A^{\odot\star}$, $q$ with $\varphi$, appends $r^{\odot\star}$ after every derivative and adds a $j$ in front of the remaining terms. Despite these similarities, the solutions differ significantly as we now solve linear operator equations on $X^{\odot\star}$. 

In the following subsection, we will illustrate how the equations for the coefficients change and how each case can be solved using the results presented in \cref{sec:DDE_eqsolve}. We focus on presenting the critical normal form coefficients, as the parameter-dependent equations are solved similarly. Full derivations of the parameter-dependent equations are presented in the supplementary material (S1.2). 

\subsection{Critical normal form coefficients}
For the computation of the critical normal form coefficients up to $c_2(0)$ we follow \cite{janssens2010normalization}. Collecting the $z^2, z\bar{z}$ and $z^3$ terms from the homological equation $\cref{eq:HOM_GH_DDE}$ results in the following linear systems
\begin{align*}
    (2i\omega_0I_n- A^{\odot\star})jH_{2000} &= B(\varphi,\varphi)r^{\odot\star}, \\
    -A^{\odot\star}jH_{1100} &= B(\varphi,\bar{\varphi})r^{\odot\star}, \\
    (3i\omega_0I_n- A^{\odot\star})jH_{3000} &= [3B(\varphi,H_{2000}) + C(\varphi,\varphi,\varphi)]r^{\odot\star}.
\end{align*}
All three equations are regular with a right-hand side of the form $(w_0,w) = (w_0,0)$. Thus, applying  \cref{cor:inverse_cases}.1 with $\eta = 0$ results in the following solutions
\begin{align}
    H_{2000}(\theta) &= e^{2i\omega_0\theta}\Delta^{-1}(2i\omega_0)B(\varphi, \varphi), \\
    H_{1100}(\theta) &= \Delta^{-1}(0)B(\varphi, \bar{\varphi}), \\
    H_{3000}(\theta) &= e^{3i\omega_0\theta}\Delta^{-1}(3i\omega_0)[3B(\varphi,H_{2000}) + C(\varphi,\varphi,\varphi)].
\end{align}
Collecting the $z^2\bar{z}$ terms yields the following singular system
\begin{align*}
    (i\omega_0I_n- A^{\odot\star})jH_{2100} &= [2B(\varphi, H_{1100}) + B(\bar{\varphi},H_{2000}) + C(\varphi,\varphi,\bar{\varphi})]r^{\odot\star}- 2c_1(0)j\varphi.
\end{align*}
Applying the Fredhom solvability condition to the above system results in 
\begin{align*}
    c_1(0) = \frac{1}{2}\langle [2B(\varphi, H_{1100}) + B(\bar{\varphi},H_{2000}) + C(\varphi,\varphi,\bar{\varphi})]r^{\odot\star}, \varphi^{\odot} \rangle.
\end{align*}
Here we used that $\langle j\varphi, \varphi^{\odot} \rangle = \langle \varphi^{\odot},\varphi \rangle = 1$ which follows from the pairings $\cref{eq:xsun_dual}$ and $\cref{eq:xsunstar_dual}$. We can evaluate the above pairing using equation $\cref{eq:xsunstar_dual}$ and the expression $\cref{eq:eigen_phi_sun}$ for $\varphi^{\odot}$. This results in the following formula 
\begin{equation*}
    c_1(0) = \frac{1}{2} p^T[2B(\varphi, H_{1100}) + B(\bar{\varphi},H_{2000}) + C(\varphi,\varphi,\bar{\varphi})].
\end{equation*}
Finally, to obtain the expression for $H_{2100}(\theta)$ we apply \cref{cor:bordered_case} to obtain the unique solution
\begin{equation*}
    H_{2100}(\theta) = B_{i\omega_0}^{INV}(2B(\varphi, H_{1100}) + B(\bar{\varphi},H_{2000}) + C(\varphi,\varphi,\bar{\varphi}), -2c_1(0)q)(\theta)
\end{equation*}
satisfying $\langle \varphi^{\odot}, H_{2100} \rangle = 0.$ Henceforth, we denote the linear combination of multilinear forms as $M_{ijkl}$, where $M_{ijkl}$ is defined as in the ODE case but with all vectors $q$ replaced by functions $\varphi$. Similarly to the first three systems we find for $H_{2200}$ the equation
\begin{align}
    H_{2200}(\theta) &= \Delta^{-1}(0)M_{2200}.
\end{align}
The $z^3\bar{z}$ terms yield the linear system
\begin{align*}
    (2i\omega_0I- A^{\odot\star})jH_{3100} &= M_{3100}r^{\odot\star} -6c_1(0)jH_{2000}.
\end{align*}
Notice that the right hand side is of the form $(w_0,w)$ with $w_0 = M_{3100} - 6c_1(0)H_{2000}(0)$ and 
\[
w(\theta) = -6c_1(0)H_{2000}(\theta) = e^{2i\omega_0\theta}\Delta^{-1}(2i\omega_0)(-6c_1(0)B(\varphi, \varphi)).
\]
Thus, we can apply \cref{cor:inverse_case1} with $\eta = -6c_1(0)B(\varphi,\varphi)$ and $\xi_1 = \xi_2 = 0$, which yields 
\begin{align}
    H_{3100}(\theta) &= e^{2i\omega_0\theta} \Delta^{-1}(2i\omega_0)M_{3100}-6c_1(0)\Delta^{-1}(2i\omega_0)[\Delta'(2i\omega_0)  - \theta \Delta(2i\omega_0)]H_{2000}(\theta).
\end{align}
Collecting the $z^3\bar{z}^2$ terms results in the equation
\begin{align}\label{eq:AjH3200}
(i\omega_0I_n- A^{\odot\star})jH_{3200 } &= M_{3200}r^{\odot\star}-12c_2(0)j\varphi - 6i\Im{\{c_1(0)\}}jH_{2100}.
\end{align}
Applying the Fredholm solvability condition gives
\begin{align*}
    c_2(0) &= \frac{1}{12}p^T M_{3200}.
\end{align*}
We proceed by deriving the coefficients needed to compute the seventh-order coefficient $c_3(0)$.  We can now apply \cref{cor:bordered_case} to find the solution
\begin{align*}
    H_{3200}(\theta) = B_{i\omega_0}^{INV}\left(M_{3200}, -12c_2(0)q - 6i\Im\{c_1(0)\}H_{2100}(0), -12i\Im\{c_1(0)\}c_1(0)q\right)(\theta),
\end{align*}
satisfying $\langle \varphi^{\odot}, H_{3200} \rangle = 0.$ Solving for the expressions of the coefficients $H_{4000}, H_{4100}$ and $H_{3300}$ is similar to how we solved for $H_{3100}$ with the application of \cref{cor:inverse_case1}. This results in the following expressions
\begin{align}
    H_{4000}(\theta) &=e^{4i\omega_0\theta}\Delta^{-1}(4i\omega_0)M_{4000} ,\\
    H_{4100}(\theta) &= e^{3i\omega_0\theta}\Delta^{-1}(3i\omega_0)M_{4100}- 12c_1(0)\Delta^{-1}(3i\omega_0)[\Delta'(3i\omega_0)  - \theta \Delta(3i\omega_0)]H_{3000}(\theta) ,
\end{align}
and
\begin{align}
    H_{3300}(\theta) &= \Delta^{-1}(0)M_{3300}-72d_2\Delta^{-1}(0)[\Delta'(0)  - \theta \Delta(0)]H_{1100}(\theta) .
\end{align}
The equation for $H_{4200}$ is of the following form
\begin{align}
    (2i\omega_0I_n - A^{\odot\star})jH_{4200} &= M_{4200}r^{\odot\star} -48c_2(0)jH_{2000} - 16c_1(0)jH_{3100} ,
\end{align}
 Using the linearity of the inverse operator, we can apply \cref{cor:inverse_case1} for the first two terms. Then to determine the inverse of the last term, we apply \cref{cor:inverse_case2} with $M = M_{3100}$, $\hat{\eta} = -6c_1(0)H_{2000}(0)$ and $\hat{\xi} = 0$. This results in the following expression
\begin{align}
    H_{4200}(\theta) &= e^{2i\omega_0\theta}\Delta^{-1}(2i\omega_0)M_{4200}- 48c_2(0)\Delta^{-1}(2i\omega_0)[\Delta'(2i\omega_0) - \theta\Delta(2i\omega_0)]H_{2000}(\theta) \nonumber \\
    &- 16c_1(0)e^{2i\omega_0\theta}\Delta^{-1}(2i\omega_0)\Bigl( [\Delta'(2i\omega_0) - \theta\Delta(2i\omega_0)]H_{3100}(0) \nonumber\\
    &+ 3c_1(0)[\Delta''(2i\omega_0) - \theta^2\Delta(2i\omega_0)]H_{2000}(0)\Bigr).
\end{align}

The third Lyapunov coefficient is derived from the system 
\begin{align}\label{eq:AjH43}
    (i\omega_0I_n - A^{\odot\star})jH_{4300} &= M_{4300}r^{\odot\star} - (144c_3(0)j\varphi + 72(2c_2(0) + \overline{c_2}(0))jH_{2100} \nonumber\\
    &+ 12i\Im{}\{c_1(0) \}jH_{3200}),
\end{align}
Applying the Fredholm solvability condition to equation $\cref{eq:AjH43}$ will then result in 
\begin{equation*}
    c_3(0) = \frac{1}{144} p^T M_{4300},
\end{equation*}
where we used that $\langle \varphi^{\odot}, H_{2100} \rangle = 0$ and $\langle \varphi^{\odot}, H_{3200} \rangle = 0$.
\subsection{Parameter-dependent coefficients}
To arrive at the equations for the unknown coefficients in the parameter transformation $K$, we can follow the same procedure as in \cref{sec:parameter_coef_ODE} with the new equations. The complete derivation is shown in the supplementary material (SM2). The resulting systems remain very similar to the systems $\cref{eq:PQ1}$. After changing all $q$ to $\varphi$ and $\bar{p}$ to $p$ one essentially only needs to determine the correct representations for all terms with an inverse matrix operator in front of it. Nevertheless, one needs to be careful. For instance, after defining $\Gamma_i(\varphi) = A_1(\varphi,e_i) + B(\varphi,\Delta^{-1}(0)J_1e_i)$, equation \cref{eq:AH10mu2} for the coefficients $H_{10\mu}$ will become
\begin{align}\label{eq:AH10mu2DDE}
    (i\omega_0I_n - A^{\odot\star})jH_{10\mu} =  [\gamma_{1,\mu}\Gamma_1(\varphi) + \gamma_{2,\mu}\Gamma_2(\varphi)]r^{\odot\star}- (\delta_{\mu}^{10} + ib_{1,\mu})j\varphi.
\end{align}
A solution $H_{10\mu}$ of equation \cref{eq:AH10mu2DDE} satisfying $\langle \varphi^{\odot}, H_{10\mu} \rangle = 0$ can be obtained with \cref{cor:bordered_case}:
\begin{align}\label{eq:H10muDDE}
    H_{10\mu}(\theta)=  B_{i\omega_0}^{INV}(\gamma_{1,\mu}\Gamma_1(\varphi) + \gamma_{2,\mu}\Gamma_2(\varphi), -(\delta_{\mu}^{10} + ib_{1,\mu})q)(\theta)
\end{align}
To free the unknown coefficients $\gamma_{1,\mu}$ and $\gamma_{2,\mu}$ we need to rewrite the equation. As with ODEs, we can apply the bordered inverse to the separate terms but this loses the original interpretation. Thus, we write $B^{Inv}_{i\omega_0}$ whenever only the notation from \cref{cor:bordered_case} is used, i.e. $v = B^{Inv}_{i\omega_0}(\eta, \xi_0, \xi_1, \ldots, \xi_m)$ is defined by equation \cref{eq:bordered_case_rep} from \cref{cor:bordered_case} but does not nevessarily solve $(i\omega_0I - A^{\odot\star})jv = jw$. With this we can rewrite \eqref{eq:H10muDDE} as
\begin{equation}\label{eq:H10muDDE2}
    H_{10\mu}=  \gamma_{1,\mu}B^{Inv}_{i\omega_0}(\Gamma_1(\varphi)) +  \gamma_{2,\mu}B^{Inv}_{i\omega_0}(\Gamma_2(\varphi)) -(\delta_{\mu}^{10} + ib_{1,\mu})B^{Inv}_{i\omega_0}(0,q).
\end{equation}
The important thing is that with all three terms combined, \cref{eq:H10muDDE2} gives us the same expression as \cref{eq:H10muDDE} which solves equation \cref{eq:AH10mu2DDE}. 
As an example, the system that solves for the coefficients $K_{10}$ and $K_{01}$ will have the following form:
\begin{equation}\label{eq:PQ1_DDE}
    P\begin{pmatrix}
        \gamma_{1,\mu} \\
        \gamma_{2,\mu}
    \end{pmatrix} = Q_{\mu},
\end{equation}
where $P \in \R^{2\times2}$ has for $k = 1,2$ the compenents
\begin{align}
    P_{1k} &= \Re{[p^T\Gamma_k(\varphi)]}, \label{eq:P1k_DDE}\\
    P_{2k} &= \frac{1}{2}\Re{}\Bigl\{ p^T\Bigl[ \Gamma_k(H_{2100}) + 2B\left(\varphi,\Delta^{-1}(0)\Pi_k(H_{1100},\bar{\varphi},\varphi)\right) \nonumber\\
    &+ B\left(\bar{\varphi},e^{2i\omega_0\theta}\Delta^{-1}(2i\omega_0)\Lambda_k(H_{2000},\varphi,\varphi)\right) +B\left(H_{2000},\overline{B_{i\omega_0}^{Inv}(\Gamma_k(\varphi))}\right) \nonumber\\
    &+2B\left(H_{1100},B_{i\omega_0}^{Inv}(\Gamma_k(\varphi))\right) + 2B_1(\varphi,H_{1100},e_k)   \nonumber\\
    &+ B_1(\bar{\varphi},H_{2000},e_k) +C\left(\varphi,\varphi,\overline{B_{i\omega_0}^{Inv}(\Gamma_k(\varphi))}\right) + 2C\left(\varphi,\bar{\varphi},B_{i\omega_0}^{Inv}(\Gamma_k(\varphi))\right) \nonumber\\
    &+2C\left(\varphi,H_{1100},\Delta^{-1}(0)J_1e_k\right) +C(\bar{\varphi},H_{2000},\Delta^{-1}(0)J_1e_k) + C_1(\varphi,\varphi,\bar{\varphi},e_k)   \nonumber\\
    &+ D(\varphi,\varphi,\bar{\varphi},\Delta^{-1}(0)J_1e_k) + \Im{[p^T\Gamma_k(\varphi)]}\Bigl(-4B\left(\varphi,\Delta^{-1}(0)\Re{[iB\left(\bar{\varphi},B_{i\omega_0}^{Inv}(0,q)\right)]} \right) \nonumber\\
    &-2iB\left(\bar{\varphi},e^{2i\omega_0\theta}\Delta^{-1}(2i\omega_0)\left([\Delta'(2i\omega_0) - \theta\Delta(2i\omega_0)]H_{2000}(0) + B\left(\varphi, B_{i\omega_0}^{Inv}(0,q) \right)  \right)\right)  \nonumber\\
    &\left.\left. \left.+iB\left(H_{2000},\overline{B_{i\omega_0}^{Inv}(0,q)}\right) - 2iB\left(H_{1100},B_{i\omega_0}^{Inv}(0,q)\right) + iC\left(\varphi,\varphi,\overline{B_{i\omega_0}^{Inv}(0,q)}\right)\right.\right.\right. \nonumber\\
    &- 2iC\left(\varphi,\bar{\varphi},B_{i\omega_0}^{Inv}(0,q)\right) \Bigr)\Bigr] \Bigr\},\label{eq:P2k_DDE}
\end{align}
 and $Q_{\mu} \in \R^2$ with $\mu = (10),(01)$ is given by
\begin{align}
   Q_{1,\mu} &= \delta_{\mu}^{10}, \label{eq:Q1mu1_DDE}\\
    Q_{2,\mu} &= \frac{1}{2}\delta_{\mu}^{10} \Re{}\Bigl\{ p^T\Bigl[ 4B\left(\varphi,\Delta^{-1}(0)\left([\Delta'(0) - \theta\Delta(0)]H_{1100}(\theta) + \Re{}\left\{ B\left(\bar{\varphi}, B_{i\omega_0}^{Inv}(0,q)\right)\right\}\right) \right) \nonumber\\
    &+2B\left(\bar{\varphi},e^{2i\omega_0\theta}\Delta^{-1}(2i\omega_0)\left([\Delta'(2i\omega_0) - \theta\Delta(2i\omega_0)]H_{2000}(0) + B\left(\varphi, B_{i\omega_0}^{Inv}(0,q) \right) \right) \right)  \nonumber\\
    &+B\left(H_{2000},\overline{B_{i\omega_0}^{Inv}(0,q)}\right) + 2B\left(H_{1100},B_{i\omega_0}^{Inv}(0,q)\right) + C\left(\varphi,\varphi,\overline{B_{i\omega_0}^{Inv}(0,q)}\right) \nonumber\\
    &+ 2C\left(\varphi,\bar{\varphi},B_{i\omega_0}^{Inv}(0,q)\right)  \Bigr]\Bigr\} + \delta_{\mu}^{01} . \label{eq:Q2mu1_DDE}
\end{align}

\section{The higher order \texorpdfstring{$LPC$}{LPC} predictors}
\label{sec:predictors}

\subsection{The predictor for ODEs}
The parameter values on the $LPC$ curve are approximated as
\begin{align*}
    \beta_1 = d_2\varepsilon^4 + 2(d_3 - a_{3201}d_2)\varepsilon^6, \quad \beta_2 = -2d_2\varepsilon^2 + (4a_{3201}d_2 - 3d_3)\varepsilon^4 ,
\end{align*}
for $\varepsilon > 0$. 
 Taking the expansion  $\cref{eq:K}$ for $K$, we have the following parameter predictor
 \begin{align}\label{eq:K_eps}
     \alpha &= \alpha_0 + K(\beta_1,  \beta_2).
 \end{align}
To approximate the solution in phase space we substitute $w = \varepsilon e^{i\psi}$ into the expansion \cref{eq:H_GH} of $H$ together with the approximations in the $\beta$ parameters, where it is not necessary to include the last three terms that are marked grey.  Thus, for $\psi \in [0,2\pi]$ the periodic orbit can be approximated as
\begin{align}
    x &= x_0 + H(\varepsilon e^{i\psi},\varepsilon e^{-i\psi}, \beta_1,\beta_2).
\end{align}
 Finally, the period is approximated by equation $\cref{eq:period}$,  
 \begin{align*}
    T &= 2\pi/(\omega_0 + (\Im{(c_1(0))} -2d_2b_{1,01})\varepsilon^2 + [d_2b_{1,10} + (4a_{3201}d_2 - 3d_3)b_{1,01} + 2d_2^2b_{1,02} \\
    &- 2d_2b_{2,01} + \Im{(c_2(0))}]\varepsilon^4).
\end{align*}
\subsection{The predictor for DDEs}
The parameter approximations for the $LPC$ curve for DDEs remain the same as for ODEs, where the expression for $K$ is given by $\cref{eq:K}$ and the period is given by $\cref{eq:period}$. Similarly, the periodic orbit is approximated for $\psi \in [0,2\pi]$ as
\begin{align}
    x &= x_0 + H(\varepsilon e^{i\psi}, \varepsilon e^{-i\psi},\beta_1,\beta_2).
\end{align}
However, now this defines a function in $C([-h,0],\R^n)$. Since, we know that the solution that we are looking for is periodic, we have that $x(\theta) = x(\theta + T)$ for all $\theta \in [-h,0]$. Thus, it is enough to ensure that $x(0) = x(T)$ when solving the system for the localisation of the periodic cycle. Therefore, we take $H_{ijkl} = H_{ijkl}(0)$ for the coefficients in the approximation of $H$ when we approximate the cycle.

\section{Examples}
\label{sec:examples}
In this section, we demonstrate the performance of the higher-order predictor. We begin with two ODE systems: the Bazykin and Khibnik prey-predator model and the extended Lorenz-84 model. For these examples, we visually compare the numerically computed Limit Point of Cycles curve with both the first-order and higher-order predictors in the parameter space to illustrate the significant improvement in accuracy.

We then conclude with a more challenging DDE system, the FitzHugh-Nagumo neural model with delay. For this case, we also provide a detailed convergence plot, which quantitatively shows that the new method achieves a higher order of convergence, providing robust evidence of its superior performance.

\subsection{ODE Examples} 
We have implemented all the equations for the higher-order predictor in the programming language Julia. For the numerical computation of the generalized Hopf points and the continuation of the $LPC$ curves, we used the existing Julia package BifurcationKit \cite{veltz2020bifurcationkitjl}. The option to compute the higher-order predictor has been added to the BifurcationKit package. In the following two examples, we will compare the new higher-order predictor to the first-order predictor from \cite{kuznetsov2008switching}. Note that in both examples, all derivatives have been computed symbolically using the Symbolic.jl package. We have included a section in the BifurcationKit documentation (cite) that walks through the Lorenz-84 example.
\subsubsection{Bazykin and Khibnik prey-predator model}
\begin{figure}[!ht]
    \centering    
    \includegraphics{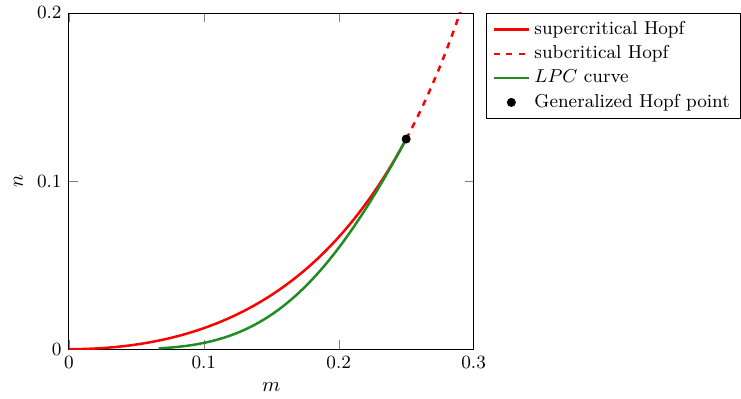}
    \caption{The bifurcation structure of the Bazykin-Khibnik model \cref{eq:BazKhib} in the $(m, n)$-parameter plane. The analytical curve of Hopf bifurcations changes stability at the generalized Hopf (GH) point, transitioning between subcritical and supercritical regimes. The GH point also serves as the origin for a curve of Limit Point of Cycles (LPC) bifurcations, which has been continued numerically.}
    \label{fig:BazKhib_bif_diagram}
\end{figure}

As a first example, we consider a version of a prey-predator system by Bazykin and Khibnik \cite{bazykin1981strict}. The model consists of the following two equations
\begin{equation}
\begin{cases} \label{eq:BazKhib}
    \dot{x} &= {\displaystyle \frac{x^2(1 - x)}{n + x} - xy},\\
    \dot{y} &= -y(m - x),
\end{cases}
\end{equation}
where $x,y \geq 0$ and $0 < m < 1$. One can show analytically that an Andronov-Hopf bifurcation occurs along the curve $n =m^2/(1 - 2m)$. The first Lyapunov coefficient along this curve is positive for $0 < m < \frac{1}{4}$ and negative for $\frac{1}{4} < m < \frac{1}{2}$, while it vanishes at $(m,n)=(\frac{1}{4},\frac{1}{8})$. Thus, a generalized Hopf bifurcation occurs in this system. The bifurcation diagram of $\cref{eq:BazKhib}$ is shown in \cref{fig:BazKhib_bif_diagram}.

\begin{figure}[!ht]
    \centering
    \includegraphics{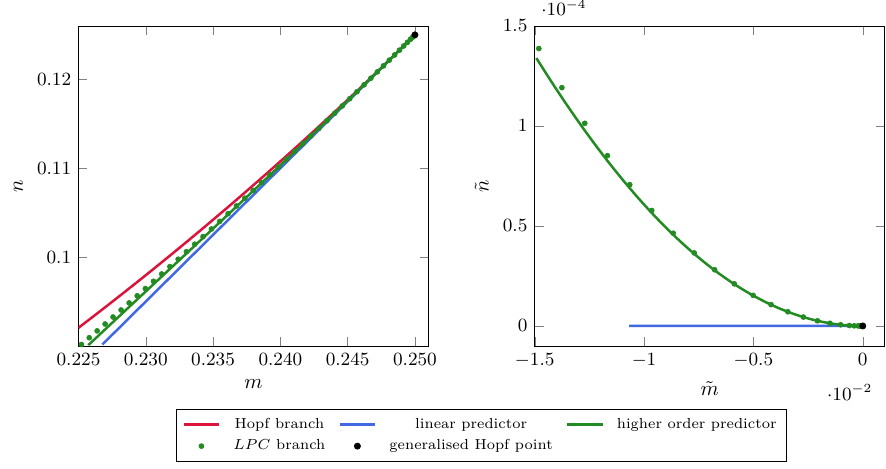}
    \caption{On the left, the bifurcation diagram near the generalized Hopf bifurcation in the  Bazykin-Khibnik Model $\cref{eq:BazKhib}$ is shown together with the first-order and the higher-order $LPC$ predictors. The figure on the right shows the $LPC$ curve and the predictors after a translation to the origin and a rotation.}
    \label{fig:BazKhib_LPC_curve}
\end{figure}

At the generalized Hopf bifurcation we have the eigenvalues $\lambda_{1,2} = \pm i\omega_0$, with $\omega_0 = \sqrt{2}/4$. We take the following eigenvectors
\begin{align*}
    q = \begin{pmatrix}
        \frac{1}{3}\sqrt{3}\\
        -\frac{1}{3}\sqrt{6}i
    \end{pmatrix}, \quad p = \begin{pmatrix}
        \frac{1}{2}\sqrt{3} \\
        -\frac{1}{4}\sqrt{6}i
    \end{pmatrix},
\end{align*}
satisfying $\bar{q}^Tq = \bar{p}^Tq = 1$. Using these eigenvectors, the second Lyapunov coefficient has the exact value\footnote{This value was calculated in MATLAB using the Symbolic Math Toolbox} $l_2 = -\frac{1024}{729}\sqrt{2}$. In our implementation in Julia, this value is numerically approximated by $l_2=-1.986494770740791$. \Cref{fig:BazKhib_LPC_curve} shows a close-up of the bifurcation diagram, including the first-order and the higher-order predictors in parameter space.

\subsubsection{The extended Lorenz-84 model}
\begin{figure}[!ht]
    \centering
    \includegraphics{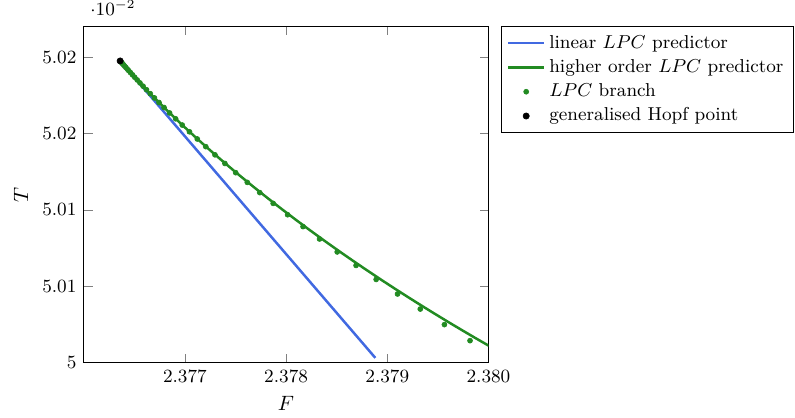}
    \caption{The numerically computed $LPC$ curve emanating from the generalized Hopf point in the extended Lorenz-84 model \cref{eq:Lorenz84} together with the first-order and higher-order predictors.}
    \label{fig:Lorenz84_LPC_curve}
\end{figure}

Our second example is an extended version of the Lorenz-84, which approximates the dynamics of an atmospheric flow model \cite{van2003baroclinic}. In this model, $X$ represents the strength of the jet stream, while $Y$ and $Z$ model the sine and cosine coefficients of the baroclinic wave. In \cite{kuznetsov2004fold}, the model was extended to include the variable $U$ to study the effect of external parameters, such as temperature, on the jet stream and baroclinic waves. The same model was used as an example in \cite{kuznetsov2008switching}. The extended model consists of four equations:
\begin{equation}\label{eq:Lorenz84}
    \begin{cases}
        \dot{X} &= -Y^2 - Z^2 - \alpha X + \alpha F - \gamma U^2,\\
        \dot{Y} &= XY -\beta XZ - Y + G,\\
        \dot{Z} &= \beta XY + XZ - Z,\\
        \dot{U} &= -\delta U + \gamma UX + T.
    \end{cases}
\end{equation}
We take $F$ and $T$ as the bifurcation parameters and fix the parameters
\[
\alpha = 0.25, \quad \beta = 1, \quad G = 0.25, \quad \delta = 1.04, \quad \gamma = 0.987.
\]
Then, the system has a generalized Hopf bifurcation at $(F,T) \approx (2.3763, 0.05019)$.  At the generalized Hopf bifurcation, we have the purely imaginary eigenvalues $\lambda = \pm i\omega_0$, with $\omega_0 = 0.690367$ and we take the eigenvectors such that $\bar{q}^Tq = \bar{p}^Tq = 1$. The second Lyapunov coefficient is computed as $l_2 = 0.22567$. \cref{fig:Lorenz84_LPC_curve} shows the first-order and the higher-order predictors in parameter space next to the numerically continued $LPC$ curve.

\subsection{DDE example: Coupled FHN neural system with delay}
The numerical computation of the generalized Hopf points and the continuation of the $LPC$ curves and the calculation of the higher-order coefficients for the higher-order predictor have been implemented in Julia. The code to produce the numerical results of the example below can be obtained upon request. We compare the higher-order predictor to the first-order predictor from \cite{bosschaert2020switching}.

We consider the coupled FitzHugh-Nagumo model from \cite{zhen2010bautin} which was also used to test the predictor in \cite{bosschaert2020switching}. This system consists of the following set of equations
\begin{equation}\label{eq:FHN}
    \begin{cases}
        \dot{u}_1(t) = -\frac{1}{3}u_1^3(t) + (c + \alpha)u_1^2(t) + du_1(t) - u_2(t) + 2\beta f(u_1(t - \tau)), \\
        \dot{u}_2(t) = \varepsilon(u_1(t) - bu_2(t)).
    \end{cases}
\end{equation}
In this model, $\alpha$ and $\beta$ measure the synaptic strength in self-connection and neighbourhood interaction, respectively, and $\tau > 0$ represents the time delay in signal transmission. The function $f$ is a sufficiently smooth sigmoidal amplification function. The parameters $b$ and $\varepsilon$ are assumed to satisfy $0 < b < 1$ and $0 < \varepsilon \ll 1$. 

\begin{figure}[!ht]
    \centering
    \includegraphics{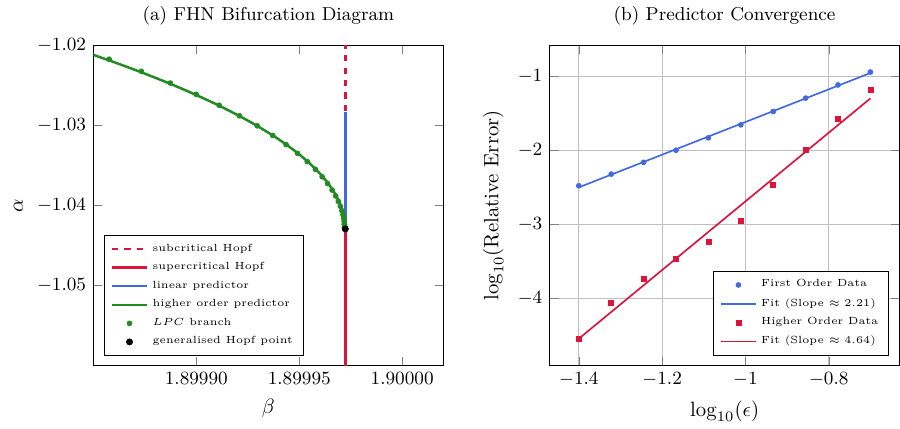}
    \caption{Comparison of predictor performance and convergence. (a) The bifurcation diagram near the generalized Hopf bifurcation in the coupled FitzHugh-Nagumo model $\cref{eq:FHN}$ together with first-order and higher-order predictors. (b) Log-log convergence plot showing the relative error of the predictors.}
    \label{fig:FHN_LPC.png}
\end{figure}

As in \cite{zhen2010bautin}, we take $\beta$ and $\alpha$ as the bifurcation parameters and fix the parameters
\[
b = 0.9, \quad \varepsilon = 0.08, \quad c = 2.0528, \quad d = -3.2135, \quad \tau = 1.7722
\]
Furthermore, we use $f(u) = \tanh(u)$ for the sigmoid amplification function. There is a generalized Hopf point at $(\beta, \alpha) = (1.9, -1.0429)$. The left and right vectors null vectors of the characteristic matrix function evaluated at the critical eigenvalues are taken such that $\bar{q}^T q = 1$ and $p^T \Delta'(i\omega_0)q = 1$. The second Lyapunov coefficient is negative and has the value $l_2 = -15.6733$. \Cref{fig:FHN_LPC.png} (a) shows that while the linear predictor is only tangent to the Hopf curve, our higher-order predictor provides a visually excellent approximation, being nearly identical to the computed LPC curve. \Cref{fig:FHN_LPC.png} (b) quantitatively evaluates the performance of the predictors by plotting their relative error against the step size $\epsilon$ on a log-log scale. The relative error is calculated between the approximate solution, as generated by the first- and higher-order predictors, and the ``true" solution, which is taken as the result from a Newton-Raphson method converged to a tolerance of 1e-8. The slope of the fitted lines reveals the order of convergence for each method. The analysis demonstrates that the higher-order predictor (red squares) achieves a convergence rate with a slope of approximately $4.64$, more than double the $\approx 2.21$ slope of the first-order predictor (blue circles). This confirms that the higher-order method provides a significantly more accurate initial guess, from which the convergence to the true solution branch is much faster.

\section{Conclusions}
\label{sec:conclusions}
In this work, we derived higher-order approximations to the $LPC$ curve emanating from the generalized Hopf bifurcation in ordinary and delay differential equations through parameter-dependent center manifold reduction and normalization. This is the only generic codimension-two equilibrium bifurcation where non-hyperbolic limit cycles emanate tangentially to the codimension-one equilibrium bifurcation curve from the bifurcation point. By doing so, we have completed a comprehensive framework for analyzing such bifurcations, providing deeper insights into their dynamics and establishing a solid foundation for future advances in the study of complex dynamical systems.

In addition to presenting a detailed derivation of the parameter-dependent normal form coefficients and transformations, we demonstrated the practical applications and robustness of our method. This was validated through several numerical examples, culminating in a convergence plot for the challenging FitzHugh-Nagumo model with delay, which provides quantitative proof of the predictor's higher-order accuracy.

The higher-order predictors developed in this work have been implemented and are publicly available on GitHub. Specifically, the method for ODEs is integrated into the Julia package BifurcationKit.jl, while the DDE implementation resides in a custom package. However, the underlying mathematical framework is general and not specific to these Julia implementations. Therefore, it could be readily integrated into other established bifurcation analysis toolboxes to benefit a wider user base. For instance, the ODE predictor could be incorporated into the widely-used MATLAB package MatCont, and similarly, the DDE predictor could be added to the well-known software DDE-BIFTOOL.

Since our method applies in all situations where the center manifold theorem holds, future research could focus on extending these techniques to a wider range of systems, including partial differential equations and stochastic models, thereby enhancing their versatility and potential for modelling complex real-world phenomena.

We conclude by noting that research on initiating the continuation of curves emanating from bifurcation points remains an active and dynamic area. In particular, starting the continuation of Shilnikov homoclinic bifurcation curves near zero-Hopf bifurcations remains a challenging and unresolved problem. Addressing this challenge will provide another powerful numerical tool for demonstrating the existence of chaos in dynamical systems, as was recently done in \cite{bosschaert2022choas}, and will complement the research reported in \cite{BALDOMA2020105162}.

%\appendix
%\section{Derivation of the parameter-dependent coefficients for DDEs} %%%\label{ap:coef_DDE}

% \section*{Acknowledgments}
% We would like to acknowledge the assistance of volunteers in putting
% together this example manuscript and supplement.

\ifarxiv
\clearpage
\hypertarget{mysupplement}{}
% === Decide between SIAM's journal class or myarticle class ===========
%
% When using the SIAM class, produce a separate supplement that
% refers to the main document using the xr-hyper package that is
% loaded in the class file.
%
% For a pre-print, skip directly to the contents of the supplement
% and include them as a separate second part of the main file.

\newif\ifsiam

% --- Uncomment precisely one of the following two. --------------------

%\siamtrue
\siamfalse

\newif\ifarxiv
\ifsiam \arxivfalse \else \arxivtrue \fi

% --- For the SIAM journal class ---------------------------------------

\ifsiam
\documentclass[review,onefignum,onetabnum]{siamonline250211}

% Optional PDF information
\ifpdf
\hypersetup{
  pdftitle={Supplementary Materials: \TheTitle},
  pdfauthor={\TheAuthors}
}
\fi

% Needed by the xr-hyper package
\externaldocument[]{article}

% Start of the document. Set title and page headings.
\begin{document}
\maketitle
\pagestyle{myheadings}
\thispagestyle{plain}
\markboth{\TheAuthors}{SWITCHING TO NONHYPERBOLIC CYCLES IN DDES}
\phantomsection\label{mysupplement}
\fi

% --- For the myarticle class ------------------------------------------
\ifarxiv
\begin{mytitle}
  \title{\textsc{supplementary materials for:}\\ \TheTitle}
  \maketitle
\end{mytitle}
\thispagestyle{plain}

\ResetCounters
\pdfbookmark[0]{Supplement}{supplement}
\fi

% === End of document class selection, beginning of contents ===========
\section{Details on the equations}
In this supplement, we present all the details of the equations not given in the paper but needed for implementation.

\subsection{Multilinear parts}\label{sec:M_ODE}
In this section, we give the full expressions of the terms containing all the multilinear forms which we did not write out in the derivation for the coefficients. The expressions follow the same order of appearance as the corresponding equations in \cref{sec:coeffs_ODE} of the main article. To change to the DDE setting, one only needs to change all vectors $q$ to the function $\varphi$.

\paragraph{The coefficient $H_{4300}$}
\begin{align*}
    M_{4300} &= 4B(q,H_{3300}) + 3B(\bar{q},H_{4200}) + 3B(\overline{H}_{2000},H_{4100}) + B(\overline{H}_{3000},H_{4000})\\
    &+ 12B(H_{1100},H_{3200}) + 12B(\overline{H}_{2100},H_{3100}) + 4B(\overline{H}_{3100},H_{3000}) \\
    &+ 6B(H_{2000},\overline{H}_{3200}) + 18B(H_{2100},H_{2200}) + 6C(q,q,\overline{H}_{3200}) \\
    &+ 12C(q,\bar{q},H_{3200}) + 12C(q,\overline{H}_{2000},H_{3100}) + 4C(q,\overline{H}_{3000},H_{3000}) \\
    &+ 36C(q,H_{1100},H_{2200}) + 36C(q,\overline{H}_{2100},H_{2100}) + 12C(q,\overline{H}_{3100},H_{2000})\\
    &+ 3C(\bar{q},\bar{q},H_{4100}) + 3C(\bar{q},\overline{H}_{2000},H_{4000}) + 24C(\bar{q},H_{1100},H_{3100}) \\
    &+ 12C(\bar{q},\overline{H}_{2100},H_{3000}) + 18C(\bar{q},H_{2000},H_{2200}) + 18C(\bar{q},H_{2100},H_{2100}) \\
    &+ 12C(\overline{H}_{2000},H_{1100},H_{3000}) + 18C(\overline{H}_{2000},H_{2000},H_{2100}) \\
    &+ 3C(\overline{H}_{3000},H_{2000},H_{2000}) + 36C(H_{1100},H_{1100},H_{2100}) \\
    &+ 36C(H_{1100},\overline{H}_{2100},H_{2000}) + 4D(q,q,q,\overline{H}_{3100}) + 18D(q,q,\bar{q},H_{2200}) \\
    &+ 18D(q,q,\overline{H}_{2000},H_{2100}) + 6D(q,q,\overline{H}_{3000},H_{2000}) + 36D(q,q,H_{1100},\overline{H}_{2100}) \\
    &+ 12D(q,\bar{q},\bar{q},H_{3100}) + 12D(q,\bar{q},\overline{H}_{2000},H_{3000}) + 72D(q,\bar{q},H_{1100},H_{2100}) \\
    &+ 36D(q,\bar{q},\overline{H}_{2100},H_{2000}) + 36D(q,\overline{H}_{2000},H_{1100},H_{2000}) \\
    &+ 24D(q,H_{1100},H_{1100},H_{1100}) + D(\bar{q},\bar{q},\bar{q},H_{4000}) + 12D(\bar{q},\bar{q},H_{1100},H_{3000}) \\
    &+ 18D(\bar{q},\bar{q},H_{2000},H_{2100}) + 9D(\bar{q},\overline{H}_{2000},H_{2000},H_{2000}) \\
    &+ 36D(\bar{q},H_{1100},H_{1100},H_{2000}) + E(q,q,q,q,\overline{H}_{3000}) + 12E(q,q,q,\bar{q},\overline{H}_{2100})\\
    &+12E(q,q,q,\overline{H}_{2000},H_{1100}) + 18E(q,q,\bar{q},\bar{q},H_{2100}) + 18E(q,q,\bar{q},\overline{H}_{2000},H_{2000}) \\
    &+ 36E(q,q,\bar{q},H_{1100},H_{1100}) + 4E(q,\bar{q},\bar{q},\bar{q},H_{3000}) + 36E(q,\bar{q},\bar{q},H_{1100},H_{2000}) \\
    &+ 3E(\bar{q},\bar{q},\bar{q},H_{2000},H_{2000}) + 3K(q,q,q,q,\bar{q},\overline{H}_{2000}) + 12K(q,q,q,\bar{q},\bar{q},H_{1100}) \\
    &+ 6K(q,q,\bar{q},\bar{q},\bar{q},H_{2000}) + L(q,q,q,q,\bar{q},\bar{q},\bar{q}).
\end{align*}
\paragraph{The coefficient $a_{3201}$}
\begin{align*}
    M_{3001} &= A_1(H_{3000},K_{01}) + 3B(q,H_{2001}) + B(H_{0001},H_{3000}) \nonumber\\
    &+ 3B(H_{1001},H_{2000}) + 3B_1(q,H_{2000},K_{01}) + 3C(q,q,H_{1001}) \nonumber\\
    &+ 3C(q,H_{0001},H_{2000}) + C_1(q,q,q,K_{01}) + D(q,q,q,H_{0001}), \\
    M_{3101} &= A_1(H_{3100},K_{01}) + 3B(q,H_{2101}) + B(\bar{q},H_{3001}) \nonumber\\
    &+ B(H_{0001},H_{3100}) + B(\overline{H}_{1001},H_{3000}) + 3B(H_{1001},H_{2100}) \nonumber\\
    &+ 3B(H_{1100},H_{2001}) + 3B(H_{1101},H_{2000}) + 3B_1(q,H_{2100},K_{01}) \nonumber\\
    &+ B_1(\bar{q},H_{3000},K_{01}) + 3B_1(H_{1100},H_{2000},K_{01}) + 3C(q,q,H_{1101}) \nonumber \\
    &+ 3C(q,\bar{q},H_{2001})+ 3C(q,H_{0001},H_{2100}) + 3C(q,\overline{H}_{1001},H_{2000}) \nonumber\\
    &+ 6C(q,H_{1001},H_{1100})+ C(\bar{q},H_{0001},H_{3000}) + 3C(\bar{q},H_{1001},H_{2000})\nonumber\\
    &+ 3C(H_{0001},H_{1100},H_{2000}) + 3C_1(q,q,H_{1100},K_{01}) + 3C_1(q,\bar{q},H_{2000},K_{01}) \nonumber\\
    &+ D(q,q,q,\overline{H}_{1001})+ 3D(q,q,\bar{q},H_{1001}) + 3D(q,q,H_{0001},H_{1100})\nonumber\\
    &+ 3D(q,\bar{q},H_{0001},H_{2000})+ D_1(q,q,q,\bar{q},K_{01}) + E(q,q,q,\bar{q},H_{0001}), \\
    M_{2201} &= A_1(H_{2200},K_{01}) + 2B(q, \overline{H}_{2101}) + 2B(\bar{q},H_{2101}) + B(H_{0001},H_{2200}) \\
    &+ 2B(\overline{H}_{1001},H_{2100}) +B(\overline{H}_{2000},H_{2001}) + B(\overline{H}_{2001},H_{2000})\\
    &+ 2B(H_{1001},\overline{H}_{2100}) + 4B(H_{1100},H_{1101}) + 2B_1(q, \overline{H}_{2100}, K_{01}) \\
    &+ B_1(\overline{H}_{2000},H_{2000},K_{01}) + 2B_1(H_{1100},H_{1100},K_{01})+ 2B_1(\bar{q}, H_{2100},K_{01}) \\
    &+ C(q,q, \overline{H}_{2001}) + 4C(q,\bar{q},H_{1101}) + 2C(q, H_{0001}, \overline{H}_{2100}) \\
    &+ 4C(q, \overline{H}_{1001}, H_{1100} ) +2C(q,\overline{H}_{2000},H_{1001}) + C(\bar{q},\bar{q},H_{2001}) \\
    &+ 2C(\bar{q},H_{0001},H_{2100}) + 2C(\bar{q},\overline{H}_{1001},H_{2000}) + 4C(\bar{q},H_{1001},H_{1100})\\
    &+C( H_{0001},\overline{H}_{2000},H_{2000}) + 2C(H_{0001},H_{1100},H_{1100} ) + C_1(q,q,\overline{H}_{2000},K_{01}) \\
    &+ 4C_1(q,\bar{q},H_{1100},K_{01}) + C_1(\bar{q},\bar{q},H_{2000},K_{01})
    + 2D(q,q,\bar{q},\overline{H}_{1001}) \\
    &+ D(q,q,H_{0001},\overline{H}_{2000}) + 2D(q,\bar{q},\bar{q},H_{1001}) + 4D(q,\bar{q},H_{0001},H_{1100})\\
    &+ D(\bar{q},\bar{q},H_{0001},H_{2000})+ D_1(q,q,\bar{q},\bar{q},K_{01})+ E(q,q,\bar{q},\bar{q},H_{0001}),\\
    M_{3201} &= A_1(H_{3200},K_{01}) + 3B(q,H_{2201}) + 2B(\bar{q},H_{3101}) + B(H_{0001},H_{3200})\nonumber\\
    &+ 2B(\overline{H}_{1001},H_{3100}) + B(\overline{H}_{2000},H_{3001}) + B(\overline{H}_{2001},H_{3000})\nonumber \\
    &+ 3B(H_{1001},H_{2200}) + 6B(H_{1100},H_{2101}) + 6B(H_{1101},H_{2100})\nonumber\\
    &+ 3B(\overline{H}_{2100},H_{2001}) + 3B(\overline{H}_{2101},H_{2000}) + 3B_1(q,H_{2200},K_{01}) \nonumber\\
    &+ 2B_1(\bar{q},H_{3100},K_{01}) + B_1(\overline{H}_{2000},H_{3000},K_{01}) + 6B_1(H_{1100},H_{2100}, K_{01})\nonumber\\
    &+ 3B_1(\overline{H}_{2100},H_{2000},K_{01}) + 3C(q,q,\overline{H}_{2101}) + 6C(q,\bar{q},H_{2101}) \nonumber\\
    &+ 3C(q,H_{0001},H_{2200}) + 6C(q,\overline{H}_{1001},H_{2100}) + 3C(q,\overline{H}_{2000},H_{2001})\nonumber\\
    &+ 3C(q,\overline{H}_{2001},H_{2000}) + 6C(q,H_{1001},\overline{H}_{2100}) + 12C(q,H_{1100},H_{1101}) \nonumber\\
    &+ C(\bar{q},\bar{q},H_{3001}) + 2C(\bar{q},H_{0001},H_{3100}) + 2C(\bar{q},\overline{H}_{1001},H_{3000}) \nonumber\\
    &+ 6C(\bar{q},H_{1001},H_{2100}) + 6C(\bar{q},H_{1100},H_{2001}) + 6C(\bar{q},H_{1101},H_{2000})\nonumber\\
    &+ C(H_{0001},\overline{H}_{2000},H_{3000}) + 6C(H_{0001},H_{1100},H_{2100}) + 3C(H_{0001},\overline{H}_{2100},H_{2000})\nonumber\\
    &+ 6C(\overline{H}_{1001}, H_{1100},H_{2000}) + 3C(\overline{H}_{2000},H_{1001},H_{2000}) \nonumber\\
    &+ 6C(H_{1001},H_{1100},H_{1100}) + 3C_1(q,q,\overline{H}_{2100},K_{01}) + 6C_1(q,\bar{q},H_{2100},K_{01})\nonumber\\
    &+ 3C_1(q,\overline{H}_{2000},H_{2000},K_{01}) + 6C_1(q,H_{1100},H_{1100},K_{01}) + C_1(\bar{q},\bar{q},H_{3000},K_{01}) \nonumber\\
    &+ 6C_1(\bar{q},H_{1100},H_{2000},K_{01}) + D(q,q,q,\overline{H}_{2001}) + 6D(q,q,\bar{q},H_{1101}) \nonumber\\
    &+ 3D(q,q,H_{0001},\overline{H}_{2100}) + 6D(q,q,\overline{H}_{1001},H_{1100}) + 3D(q,q,\overline{H}_{2000},H_{1001})\nonumber\\
    &+ 3D(q,\bar{q},\bar{q},H_{2001}) + 6D(q,\bar{q},H_{0001},H_{2100}) + 6D(q,\bar{q},\overline{H}_{1001},H_{2000})\nonumber\\
    &+ 12D(q,\bar{q},H_{1001},H_{1100}) + 3D(q,H_{0001},\overline{H}_{2000},H_{2000})\nonumber\\
    &+ 6D(q,H_{0001},H_{1100},H_{1100}) + D(\bar{q},\bar{q},H_{0001},H_{3000}) + 3D(\bar{q},\bar{q},H_{1001},H_{2000})\nonumber\\
    &+ 6D(\bar{q},H_{0001},H_{1100},H_{2000}) + D_1(q,q,q,\overline{H}_{2000},K_{01}) + 6D_1(q,q,\bar{q},H_{1100},K_{01}) \nonumber\\
    &+ 3D_1(q,\bar{q},\bar{q},H_{2000},K_{01}) + 2E(q,q,q,\bar{q},\overline{H}_{1001}) + E(q,q,q,H_{0001},\overline{H}_{2000})\nonumber \\
    &+ 3E(q,q,\bar{q},\bar{q},H_{1001}) + 6E(q,q,\bar{q},H_{0001},H_{1100}) + 3E(q,\bar{q},\bar{q},H_{0001},H_{2000}) \nonumber\\
    &+ E_1(q,q,q,\bar{q},\bar{q},K_{01}) + K(q,q,q,\bar{q},\bar{q},H_{0001}).
\end{align*}

 \paragraph{The coefficient $K_{02}$} 
 \begin{align*}
    M_{0002} &= 2A_1(H_{0001},K_{01}) + B(H_{0001},H_{0001})+J_2(K_{01},K_{01}), \\
    M_{1002} &= 2A_1(H_{1001},K_{01}) + 2B(H_{0001},H_{1001}) + A_2(q,K_{01},K_{01})\\
    &+ 2B_1(q,H_{0001},K_{01}) + C(q,H_{0001},H_{0001}), \\
    M_{2002} &= 2A_1(H_{2001},K_{01})  + 2B(H_{0001},H_{2001})  + 2B(H_{1001},H_{1001}) \\
    &+ A_2(H_{2000}, K_{01},K_{01}) + 4B_1(q,H_{1001},K_{01}) + 2B_1(H_{0001},H_{2000},K_{01})\\
    &+ 4C(q,H_{0001},H_{1001}) + C(H_{0001},H_{0001},H_{2000}) + B_2(q,q,K_{01},K_{01})\\
    &+ 2C_1(q,q,H_{0001},K_{01}) + D(q,q,H_{0001},H_{0001}), \\
    M_{1102} &= 2A_1(H_{1101},K_{01})  + 2B(H_{0001},H_{1101})  + 2B(\overline{H}_{1001},H_{1001})\\
    &+ A_2(H_{1100},K_{01},K_{01}) + 4\Re{}\left(B_1(\bar{q},H_{1001},K_{01})\right) + 2B_1(H_{0001},H_{1100},K_{01})\\
    &+ 4\Re{}\left(C(\bar{q},H_{0001},H_{1001})\right) + C(H_{0001},H_{0001},H_{1100}) + B_2(q,\bar{q},K_{01},K_{01}) \\
    &+ 2C_1(q,\bar{q},H_{0001},K_{01}) + D(q,\bar{q},H_{0001},H_{0001}), \\
    M_{2102} &= 2A_1(H_{2101},K_{01})  + 2B(H_{0001},H_{2101})+ 2B(\overline{H}_{1001},H_{2001})  + 4B(H_{1001},H_{1101})\\
    &+ A_2(H_{2100},K_{01},K_{01}) + 4B_1(q,H_{1101},K_{01}) + 2B_1(\bar{q},H_{2001},K_{01}) \\
    &+ 2B_1(H_{0001},H_{2100},K_{01}) + 2B_1(\overline{H}_{1001},H_{2000},K_{01})+ 4B_1(H_{1001},H_{1100},K_{01})\\
    &+ 4C(q,H_{0001},H_{1101})  + 4C(q,\overline{H}_{1001},H_{1001}) + 2C(\bar{q},H_{0001},H_{2001}) \\
    &+ 2C(\bar{q},H_{1001},H_{1001}) + C(H_{0001},H_{0001},H_{2100}) + 2C(H_{0001},\overline{H}_{1001},H_{2000}) \\
    &+ 4C(H_{0001},H_{1001},H_{1100}) + 2B_2(q,H_{1100},K_{01}, K_{01}) + B_2(\bar{q},H_{2000},K_{01},K_{01})  \\
    &+ 2C_1(q,q,\overline{H}_{1001},K_{01})+ 4C_1(q,\bar{q},H_{1001},K_{01}) + 4C_1(q,H_{0001},H_{1100},K_{01}) \\
    &+ 2C_1(\bar{q},H_{0001},H_{2000},K_{01})  +2D(q,q,H_{0001},\overline{H}_{1001}) + 4D(q,\bar{q},H_{0001},H_{1001}) \\
    &+ 2D(q,H_{0001},H_{0001},H_{1100}) + D(\bar{q},H_{2000},H_{0001},H_{0001}) + C_2(q,q,\bar{q},K_{01},K_{01}) \\
    &+ 2D_1(q,q,\bar{q},H_{0001},K_{01}) + E(q,q,\bar{q},H_{0001},H_{0001}).
\end{align*}

\paragraph{The coefficient $K_{11}$} 
\begin{align*}
    M_{0011} &= A_1(H_{0010},K_{01}) + A_1(H_{0001},K_{10}) + B(H_{0001},H_{0010}) + J_2(K_{01},K_{10}), \\
    M_{1011} &= A_1(H_{1010},K_{01}) + A_1(H_{1001},K_{10})  + B(H_{0001},H_{1010}) \\
    &+B(H_{0010},H_{1001}) + A_2(q,K_{01},K_{10}) + B_1(q,H_{0010},K_{01}) \\
    &+ B_1(q,H_{0001},K_{10}) + C(q,H_{0001},H_{0010}),\\
    M_{2011} &= A_1(H_{2010},K_{01}) + A_1(H_{2001},K_{10})+ B(H_{0001},H_{2010}) \\
    &+ B(H_{0010},H_{2001}) +  2B(H_{1001},H_{1010}) + A_2(H_{2000},K_{01},K_{10})  \\
    &+ 2B_1(q,H_{1010},K_{01}) + 2B_1(q,H_{1001},K_{10}) + B_1(H_{0010},H_{2000},K_{01}) \\
    &+ B_1(H_{0001},H_{2000},K_{10})  + 2C(q,H_{0001},H_{1010}) + 2C(q,H_{0010},H_{1001}) \\
    &+ C(H_{0001},H_{0010},H_{2000}) + B_2(q,q,K_{10},K_{01}) + C_1(q,q,H_{0010},K_{01}) \\
    &+ C_1(q,q,H_{0001},K_{10}) + D(q,q,H_{0001},H_{0010}), \\
    M_{1111} &= A_1(H_{1110},K_{01}) + A_1(H_{1101},K_{10})   + B(H_{0001},H_{1110}) \\
    &+ B(H_{0010},H_{1101})  + 2\Re{}\left(B(\overline{H}_{1001},H_{1010})\right)+ A_2(H_{1100},K_{01},K_{10})  \\
    &+ 2\Re{}\left(B_1(\bar{q},H_{1010},K_{01})\right) + 2\Re{}\left(B_1(\bar{q},H_{1001},K_{10})\right)  + B_1(H_{0010},H_{1100},K_{01})\\
    &+ B_1(H_{0001},H_{1100},K_{10})  +  2\Re{}\left(C(\bar{q},H_{0001},H_{1010})\right)+ 2\Re{}\left(C(\bar{q},H_{0010},H_{1001})\right)\\
    &+ C(H_{0001},H_{0010},H_{1100}) + B_2(q,\bar{q},K_{01},K_{10}) + C_1(q,\bar{q},H_{0010},K_{01}) \\
    &+ C_1(q,\bar{q},H_{0001},K_{10}) + D(q,\bar{q},H_{0001},H_{0010}), \\
    M_{2111} &= A_1(H_{2110},K_{01})+A_1(H_{2101},K_{10})+B(H_{0001},H_{2110})+B(H_{0010},H_{2101})\\
    &+B(\overline{H}_{1001},H_{2010})+B(\overline{H}_{1010},H_{2001})+2 B(H_{1001},H_{1110})+2 B(H_{1010},H_{1101})\\
    &+A_2(H_{2100},K_{01},K_{10})+2 B_1(q,H_{1110},K_{01})+2 B_1(q,H_{1101},K_{10})+B_1(\bar{q},H_{2010},K_{01})\\
    &+B_1(\bar{q},H_{2001},K_{10})+B_1(H_{0010},H_{2100},K_{01})+B_1(\overline{H}_{1010},H_{2000},K_{01}) \\
    &+2B_1(H_{1010},H_{1100},K_{01})+B_1(H_{0001},H_{2100},K_{10})+B_1(\overline{H}_{1001},H_{2000},K_{10}) \\
    &+2B_1(H_{1001},H_{1100},K_{10})+2 C(q,H_{0001},H_{1110})+2 C(q,H_{0010},H_{1101})\\
    &+2 C(q,\overline{H}_{1001},H_{1010})+2 C(q,\overline{H}_{1010},H_{1001})+C(\bar{q},H_{0001},H_{2010})\\
    &+C(\bar{q},H_{0010},H_{2001})+2C(\bar{q},H_{1001},H_{1010})+C(H_{0001},H_{0010},H_{2100})\\
    &+C(H_{0001},\overline{H}_{1010},H_{2000})+2 C(H_{0001},H_{1010},H_{1100})+C(H_{0010},\overline{H}_{1001},H_{2000})\\
    &+2C(H_{0010},H_{1001},H_{1100})+2B_2(q,H_{1100},K_{01},K_{10})+B_2(\bar{q},H_{2000},K_{01},K_{10})\\
    &+C_1(q,q,\overline{H}_{1010},K_{01})+C_1(q,q,\overline{H}_{1001},K_{10})+2C_1(q,\bar{q},H_{1010},K_{01})\\
    &+2C_1(q,\bar{q},H_{1001},K_{10})+2C_1(q,H_{0010},H_{1100},K_{01})+2C_1(q,H_{0001},H_{1100},K_{10})\\
    &+C_1(\bar{q},H_{0010},H_{2000},K_{01})+C_1(\bar{q},H_{0001},H_{2000},K_{10})+D(q,q,H_{0001},\overline{H}_{1010})\\
    &+D(q,q,H_{0010},\overline{H}_{1001})+2 D(q,\bar{q},H_{0001},H_{1010})+2 D(q,\bar{q},H_{0010},H_{1001})\\
    &+2 D(q,H_{0001},H_{0010},H_{1100})+D(\bar{q},H_{0001},H_{0010},H_{2000})+C_2(q,q,\bar{q},K_{01},K_{10})\\
    &+D_1(q,q,\bar{q},H_{0010},K_{01})+D_1(q,q,\bar{q},H_{0001},K_{10}) +E(q,q,\bar{q},H_{0001},H_{0010}).
\end{align*}

\paragraph{The coefficient $K_{03}$} 
 \begin{align*}
     M_{0003} &=  3A_1(H_{0002},K_{01}) + 3A_1(H_{0001},K_{02}) + 3B(H_{0001}, H_{0002}) \\
     &+ 3J_2(K_{01}, K_{02}) + 3A_2(H_{0001},K_{01}, K_{01}) + 3B_1(H_{0001}, H_{0001},K_{01}) \\
     &+ J_3(K_{01}, K_{01}, K_{01}) + C(H_{0001}, H_{0001}, H_{0001}),\\
     M_{1003} &= 3A_1(H_{1002},K_{01}) + 3A_1(H_{1001},K_{02})+ 3B(H_{0001}, H_{1002})\\
     &+ 3B(H_{0002}, H_{1001})+ 3A_2(q, K_{01}, K_{02})+ 3A_2(H_{1001},K_{01}, K_{01})\\
     &+ 3B_1(q, H_{0002}, K_{01})+ 3B_1(q, H_{0001}, K_{02})+ 6B_1(H_{0001}, H_{1001},K_{01})\\
     &+ 3C(q, H_{0001}, H_{0002})+ 3C(H_{0001}, H_{0001}, H_{1001})+ A_3(q, K_{01}, K_{01}, K_{01})\\
     &+ 3B_2(q, H_{0001}, K_{01}, K_{01}) + 3C_1(q, H_{0001}, H_{0001}, K_{01})  \\
     &+ D(q, H_{0001}, H_{0001}, H_{0001}), \\
     M_{2003} &= 3A_1(H_{2002},K_{01}) + 3A_1( H_{2001},K_{02}) +  3B(H_{0001}, H_{2002}) + 3B(H_{0002}, H_{2001})\\
     &+ 6B(H_{1001}, H_{1002}) +3A_2(H_{2001},K_{01}, K_{01}) + 3A_2(H_{2000},K_{01}, K_{02})\\
     &+ 6B_1(q, H_{1002}, K_{01}) + 6B_1(q, H_{1001}, K_{02}) + 6B_1(H_{0001}, H_{2001},K_{01}) \\
     &+ 3B_1(H_{0002}, H_{2000},K_{01}) + 6B_1(H_{1001}, H_{1001},K_{01})  +3B_1(H_{0001}, H_{2000},K_{02})\\
     &+ 6C(q, H_{0001}, H_{1002}) + 6C(q, H_{0002}, H_{1001}) + 3C(H_{0001}, H_{0001}, H_{2001}) \\
     &+ 3C(H_{0001}, H_{0002}, H_{2000})+ 6C(H_{0001}, H_{1001}, H_{1001}) + A_3(H_{2000},K_{01}, K_{01}, K_{01}) \\
     &+ 3B_2(q, q, K_{01}, K_{02}) + 6B_2(q, H_{1001}, K_{01}, K_{01}) + 3B_2( H_{0001}, H_{2000},K_{01}, K_{01})\\
     &+ 3C_1(q, q, H_{0002}, K_{01}) + 3C_1(q, q, H_{0001}, K_{02}) + 12C_1(q, H_{0001}, H_{1001}, K_{01}) \\
     &+ 3C_1(H_{0001}, H_{0001}, H_{2000},K_{01}) + 3D(q, q, H_{0001}, H_{0002})+ 6D(q, H_{0001}, H_{0001}, H_{1001}) \\
     &+ D(H_{0001}, H_{0001}, H_{0001}, H_{2000}) + B_3(q, q, K_{01}, K_{01}, K_{01}) + 3C_2(q, q, H_{0001}, K_{01}, K_{01})\\
     &+3D_1(q, q, H_{0001}, H_{0001}, K_{01}) + E(q, q, H_{0001}, H_{0001}, H_{0001}),\\
     M_{1103} &= 3A_1(H_{1102},K_{01}) + 3A_1(H_{1101},K_{02}) + 3B(H_{0001}, H_{1102}) + 3B(H_{0002}, H_{1101})\\
     &+ 3B(\overline{H}_{1001}, H_{1002}) + 3B(\overline{H}_{1002}, H_{1001}) +3A_2(H_{1101},K_{01}, K_{01}) \\
     &+ 3A_2(H_{1100},K_{01}, K_{02})  + 6\Re{}(B_1(\bar{q}, H_{1002}, K_{01}))+ 6\Re{}(B_1(\bar{q}, H_{1001}, K_{02})) \\
     &+ 6B_1(H_{0001}, H_{1101},K_{01})+ 3B_1( H_{0002}, H_{1100},K_{01}) + 6B_1(\overline{H}_{1001}, H_{1001},K_{01}) \\
     &+ 3B_1(H_{0001},H_{1100},K_{02}) + 6\Re{}(C(\bar{q}, H_{0001}, H_{1002}))+ 6\Re{}(C(\bar{q}, H_{0002}, H_{1001})) \\
     &+ 3C(H_{0001}, H_{0001}, H_{1101})+ 3C(H_{0001}, H_{0002}, H_{1100}) + 6C(H_{0001}, \overline{H}_{1001}, H_{1001}) \\
     &+ A_3(H_{1100},K_{01}, K_{01}, K_{01}) + 3B_2(q, \bar{q}, K_{01}, K_{02}) + 6\Re{}(B_2(\bar{q}, H_{1001}, K_{01}, K_{01})) \\
     & + 3B_2(H_{0001}, H_{1100},K_{01}, K_{01}) + 3C_1(q, \bar{q}, H_{0002}, K_{01}) + 3C_1(q, \bar{q}, H_{0001}, K_{02}) \\
     &+ 12\Re{}(C_1(\bar{q}, H_{0001}, H_{1001}, K_{01})) + 3C_1(H_{0001}, H_{0001}, H_{1100},K_{01}) \\
     &+ 3D(q, \bar{q}, H_{0001}, H_{0002}) + 6\Re{}(D(\bar{q}, H_{0001}, H_{0001}, H_{1001}) )\\
     &+ D(H_{0001}, H_{0001}, H_{0001}, H_{1100})+ B_3(q, \bar{q}, K_{01}, K_{01}, K_{01}) \\
     &+ 3C_2(q, \bar{q}, H_{0001}, K_{01}, K_{01}) + 3D_1(q, \bar{q}, H_{0001}, H_{0001}, K_{01}) \\
     &+E(q, \bar{q}, H_{0001}, H_{0001}, H_{0001}), \\
     M_{2103} &=  3A_1(H_{2102},K_{01}) + 3A_1(H_{2101},K_{02})+ 3B(H_{0001}, H_{2102}) + 3B(H_{0002}, H_{2101})\\
     &+ 3B(\overline{H}_{1001}, H_{2002})+ 3B(\overline{H}_{1002}, H_{2001})+ 6B(H_{1001}, H_{1102}) + 6B(H_{1002}, H_{1101})\\
     &+ 3A_2(H_{2101},K_{01}, K_{01}) + 3A_2(H_{2100},K_{01}, K_{02}) + 6B_1(q, H_{1102}, K_{01}) \\
     &+ 6B_1(q, H_{1101}, K_{02}) + 3B_1(\bar{q}, H_{2002}, K_{01}) + 3B_1(\bar{q}, H_{2001}, K_{02}) \\
     &+ 6B_1(H_{0001}, H_{2101},K_{01}) + 3B_1( H_{0002}, H_{2100},K_{01}) + 6B_1(\overline{H}_{1001}, H_{2001},K_{01}) \\
     &+ 3B_1(\overline{H}_{1002}, H_{2000},K_{01}) + 12B_1(H_{1001}, H_{1101},K_{01}) + 6B_1(H_{1002}, H_{1100},K_{01})\\
     &+ 3B_1(H_{0001}, H_{2100},K_{02}) + 3B_1(\overline{H}_{1001}, H_{2000},K_{02}) + 6B_1(H_{1001}, H_{1100},K_{02})\\
     &+ 6C(q,H_{0001},H_{1102}) + 6C(q,H_{0002},H_{1101}) + 6C(q,\overline{H}_{1001},H_{1002}) \\
     &+ 6C(q,\overline{H}_{1002},H_{1001}) + 3C(\bar{q},H_{0001},H_{2002})+ 3C(\bar{q},H_{0002},H_{2001})\\
     &+ 6C(\bar{q},H_{1001},H_{1002})+3C(H_{0001},H_{0001},H_{2101}) + 3C(H_{0001},H_{0002},H_{2100})\\
     &+6C(H_{0001},\overline{H}_{1001},H_{2001}) + 3C(H_{0001},\overline{H}_{1002},H_{2000}) +12 C(H_{0001},H_{1001},H_{1101}) \\
     &+ 6 C(H_{0001},H_{1002},H_{1100}) +3C(H_{0002},\overline{H}_{1001},H_{2000}) + 6C(H_{0002},H_{1001},H_{1100}) \\
     &+ 6C(\overline{H}_{1001},H_{1001},H_{1001})+ A_3(H_{2100},K_{01},K_{01},K_{01}) + 6B_2(q,H_{1101},K_{01},K_{01}) \\
     &+ 6B_2(q,H_{1100},K_{01},K_{02}) + 3B_2(\bar{q},H_{2001},K_{01},K_{01}) + 3B_2(\bar{q},H_{2000},K_{01},K_{02}) \\
     &+3B_2(H_{0001},H_{2100},K_{01},K_{01}) + 3B_2(\overline{H}_{1001},H_{2000},K_{01},K_{01})\\
     &+ 6B_2(H_{1001},H_{1100},K_{01},K_{01}) +3C_1(q,q,\overline{H}_{1002},K_{01}) + 3C_1(q,q,\overline{H}_{1001},K_{02}) \\
     &+ 6C_1(q,\bar{q},H_{1002},K_{01})+ 6C_1(q,\bar{q},H_{1001},K_{02}) + 12C_1(q,H_{0001},H_{1101},K_{01}) \\
     &+ 6C_1(q,H_{0002},H_{1100},K_{01}) + 12C_1(q,\overline{H}_{1001},H_{1001},,K_{01}) \\
     &+ 6C_1(q,H_{0001},H_{1100},K_{02}) + 6C_1(\bar{q},H_{0001},H_{2001},K_{01}) \\
     &+ 3C_1(\bar{q},H_{0002},H_{2000},K_{01}) + 6 C_1(\bar{q},H_{1001},H_{1001},K_{01}) \\
     &+ 3C_1(\bar{q},H_{0001},H_{2000},K_{02}) + 3C_1(H_{0001},H_{0001},H_{2100},K_{01}) \\
     &+ 6C_1(H_{0001},\overline{H}_{1001},H_{2000},K_{01}) + 12 C_1(H_{0001},H_{1001},H_{1100},K_{01})\\
     &+ 3 D(q,q,H_{0001},\overline{H}_{1002})+ 3D(q,q,H_{0002},\overline{H}_{1001}) + 6 D(q,\bar{q},H_{0001},H_{1002})\\
     &+ 6 D(q,\bar{q},H_{0002},H_{1001}) + 6 D(q,H_{0001},H_{0001},H_{1101}) + 6 D(q,H_{0001},H_{0002},H_{1100}) \\
     &+ 12 D(q,H_{0001},\overline{H}_{1001},H_{1001}) + 3 D(\bar{q},H_{0001},H_{0001},H_{2001}) \\
     &+ 3 D(\bar{q},H_{0001},H_{0002},H_{2000}) + 6 D(\bar{q},H_{0001},H_{1001},H_{1001}) \\
     &+ D(H_{0001},H_{0001},H_{0001},H_{2100}) + 3D(H_{0001},H_{0001},\overline{H}_{1001},H_{2000})\\
     &+ 6D(H_{0001},H_{0001},H_{1001},H_{1100})+2B_3(q,H_{1100},K_{01},K_{01},K_{01}) \\
     &+ B_3(\bar{q},H_{2000},K_{01},K_{01},K_{01}) +3C_2(q,q,\bar{q},K_{01},K_{02}) + 3C_2(q,q,\overline{H}_{1001},K_{01},K_{01}) \\
     &+6C_2(q,\bar{q},H_{1001},K_{01},K_{01}) + 6C_2(q,H_{0001},H_{1100},K_{01},K_{01}) \\
     &+3C_2(\bar{q},H_{0001},H_{2000},K_{01},K_{01})+3D_1(q,q,\bar{q},H_{0002},K_{01}) \\
     &+ 3D_1(q,q,\bar{q},H_{0001},K_{02}) +6D_1(q,q,H_{0001},\overline{H}_{1001},K_{01}) \\
     &+12D_1(q,\bar{q},H_{0001},H_{1001},K_{01}) +6D_1(q,H_{0001},H_{0001},H_{1100},K_{01}) \\
     &+3D_1(\bar{q},H_{0001},H_{0001},H_{2000},K_{01}) +  3E(q,q,\bar{q},H_{0001},H_{0002}) \\
     &+ 3E(q,q,H_{0001},H_{0001},\overline{H}_{1001}) + 6E(q,\bar{q},H_{0001},H_{0001},H_{1001}) \\
     &+ 2E(q,H_{0001},H_{0001},H_{0001},H_{1100}) + E(\bar{q},H_{0001},H_{0001},H_{0001},H_{2000}) \\
     &+C_3(q,q,\bar{q},K_{01},K_{01},K_{01}) + 3E_1(q,q,\bar{q},H_{0001},H_{0001},K_{01}) \\
     &+ K(q,q,\bar{q},H_{0001},H_{0001},H_{0001}).
 \end{align*}
 
\newpage
\subsection{Derivation of the parameter-dependent coefficients for DDEs}
\paragraph{Linear coefficients $K_{10}$,$K_{01}$}
To find the equations from which we can solve for $K_{10},K_{01}$, we follow the same steps as in the ODE case. Collecting the $\beta_1$ and $\beta_2$ terms in \cref{eq:HOM_GH_DDE} yields for $\mu = (10),(01)$ the systems
\begin{align}\label{eq:AjH00mu1}
    -A^{\odot\star}jH_{00\mu} = J_1K_{\mu}r^{\odot\star}.
\end{align}
Write
\begin{equation}\label{eq:Kmu1DDE}
    K_{\mu} = \gamma_{1,\mu}e_1 + \gamma_{2,\mu}e_2,
\end{equation}
for unknown $\gamma_{1,\mu},\gamma_{2,\mu} \in \R$ and $e_1,e_2 \in \R^2$ the standard basis vectors. Substituting equation \cref{eq:Kmu1DDE} into \cref{eq:AjH00mu1} and applying  \cref{cor:inverse_cases}.1 results in
\begin{equation}\label{eq:H00muDDE}
    H_{00\mu}(\theta) = \gamma_{1,\mu}\Delta^{-1}(0)J_1e_1 + \gamma_{2,\mu}\Delta(0)^{-1}J_1e_2.
\end{equation}
The $\beta_1z$ and $\beta_2z$ terms yield the systems
\begin{align}\label{eq:AH10muDDE}
    (i\omega_0I_n-A^{\odot\star})jH_{10\mu} =  [A_1(\varphi,K_{\mu}) +B(\varphi,H_{00\mu})]r^{\odot\star} -(\delta_{\mu}^{10} + ib_{1,\mu})j\varphi,
\end{align}
where $\delta_{\mu}^{10}$ is defined as in \cref{eq:delta_mu}. Define $\Gamma_i(\varphi) = A_1(\varphi,e_i) + B(\varphi,\Delta^{-1}(0)J_1e_i)$.
Substituting equations \cref{eq:Kmu1DDE} and \cref{eq:H00muDDE} into equation \cref{eq:AH10muDDE} then yields
\begin{align}\label{S-eq:AH10mu2DDE}
    (i\omega_0I_n - A^{\odot\star})jH_{10\mu} =  [\gamma_{1,\mu}\Gamma_1(\varphi) + \gamma_{2,\mu}\Gamma_2(\varphi)]r^{\odot\star}- (\delta_{\mu}^{10} + ib_{1,\mu})j\varphi.
\end{align}

Applying the Fredholm solvability condition to equation \cref{S-eq:AH10mu2DDE} results in
\begin{equation}
    \delta_{\mu}^{10} + ib_{1,\mu} = p^T[\gamma_{1,\mu}\Gamma_1(\varphi) + \gamma_{2,\mu}\Gamma_2(\varphi)].
\end{equation}
Taking the real and imaginary parts yield
\begin{equation}
    \delta_{\mu}^{10} = \gamma_{1,\mu}\Re{[p^T\Gamma_1(\varphi)]} + \gamma_{2,\mu}\Re{[p^T\Gamma_2(\varphi)]},
\end{equation}
and
\begin{equation}\label{eq:b1mu1DDE}
    b_{1,\mu} = \gamma_{1,\mu}\Im{[p^T\Gamma_1(\varphi)]} + \gamma_{2,\mu}\Im{[p^T\Gamma_2(\varphi)]}.
\end{equation}
 A solution $H_{10\mu}$ of equation \cref{S-eq:AH10mu2DDE} satisfying $\langle \varphi^{\odot}, H_{10\mu} \rangle = 0$ can be obtained with \cref{cor:bordered_case}:
\begin{align}\label{S-eq:H10muDDE}
    H_{10\mu}(\theta)=  B_{i\omega_0}^{INV}(\gamma_{1,01}\Gamma_1(\varphi) + \gamma_{2,01}\Gamma_2(\varphi), -(\delta_{\mu}^{10} + ib_{1,\mu})q)(\theta).
\end{align}
To free the unknown coefficients $\gamma_{1,\mu}$ and $\gamma_{2,\mu}$ we need to rewrite the equation as we did for ODEs. As with ODEs, applying the bordered inverse to the separate terms loses the original interpretation. Thus, we write $B^{Inv}_{i\omega_0}$ whenever only the notation from \cref{cor:bordered_case} is used, i.e. $v = B^{Inv}_{i\omega_0}(\eta, \xi_0, \xi_1, \ldots, \xi_m)$ is defined by equation \cref{eq:bordered_case_rep} from \cref{cor:bordered_case} but does not solve $(i\omega_0I_n- A^{\odot\star})jv = jw$. With this we can rewrite \eqref{eq:H10muDDE} as
\begin{equation}\label{S-eq:H10muDDE2}
    H_{10\mu}(\theta) =  \gamma_{1,\mu}B^{Inv}_{i\omega_0}(\Gamma_1(\varphi))(\theta) +  \gamma_{2,\mu}B^{Inv}_{i\omega_0}(\Gamma_2(\varphi))(\theta) -(\delta_{\mu}^{10} + ib_{1,\mu})B^{Inv}_{i\omega_0}(0,q)(\theta).
\end{equation} The important thing is that with all three terms combined, \cref{S-eq:H10muDDE2} gives us the same expression as \cref{S-eq:H10muDDE} which solves equation \cref{S-eq:AH10mu2DDE}. Collecting the $z^2\beta_i$ and $z\bar{z}\beta_i$ terms respectively yield the systems
\begin{align*}
     (2i\omega_0I_n- A^{\odot\star})jH_{20\mu} &= [A_1(H_{2000},K_{\mu}) +2B(\varphi,H_{10\mu})+ B(H_{00\mu},H_{2000})  \\
    &+ B_1(\varphi,\varphi,K_{\mu}) + C(\varphi,\varphi,H_{00\mu})]r^{\odot\star} - 2(\delta_{\mu}^{10} + ib_{1,\mu})jH_{2000} , \\
    -A^{\odot\star}jH_{11\mu} &=  [A_1(H_{1100},K_{\mu}) + 2\Re{\{B(\bar{\varphi},H_{10\mu})\}} + B(H_{00\mu},H_{1100}) \\
    &+ B_1(\varphi,\bar{\varphi},K_{\mu}) + C(\varphi,\bar{\varphi},H_{00\mu})]r^{\odot\star}
    -2\delta_{\mu}^{10}jH_{1100}.
\end{align*}
Both these systems are regular and their solutions follow from \cref{cor:inverse_case1}. The resulting equations are
\begin{align}
     H_{20\mu}(\theta) &=  e^{2i\omega_0\theta}\Delta^{-1}(2i\omega_0)[A_1(H_{2000},K_{\mu}) +2B(\varphi,H_{10\mu})+ B(H_{00\mu},H_{2000}) \nonumber\\
    &+ B_1(\varphi,\varphi,K_{\mu}) + C(\varphi,\varphi,H_{00\mu})] \nonumber\\
    &- 2(\delta_{\mu}^{10} + ib_{1,\mu})\Delta^{-1}(2i\omega_0)[\Delta'(2i\omega_0) - \theta\Delta(2i\omega_0)]H_{2000}(\theta) , \label{eq:H20mu1_DDE}\\
    H_{11\mu}(\theta) &= \Delta^{-1}(0)[A_1(H_{1100},K_{\mu}) + 2\Re{(B(\bar{\varphi},H_{10\mu}))} + B(H_{00\mu},H_{1100}) \nonumber\\
    &+B_1(\varphi,\bar{\varphi},K_{\mu}) + C(\varphi,\bar{\varphi},H_{00\mu})] -2\delta_{\mu}^{10}\Delta^{-1}(0)[\Delta'(0) - \theta\Delta(0)]H_{1100}(\theta), \label{eq:H11mu1_DDE}
\end{align}
We now substitute equations \cref{eq:Kmu1DDE}, \cref{eq:H00muDDE} and \cref{S-eq:H10muDDE} into the above expressions. For this, it is convenient to define for $i = 1,2$ the following functions 
\begin{align*}
\Lambda_i(u,v,w) &= \Gamma_i(u) + 2B\left(v,B_{i\omega_0}^{Inv}(\Gamma_i(w))\right) + B_1(v,w,e_i) + C(v,w,\Delta^{-1}(0)J_1e_i), \\
\Pi_i(u,v,w) &= \Gamma_i(u) + 2\Re{\{B\left(v,B_{i\omega_0}^{Inv}(\Gamma_i(w))\right)\}} + B_1(v,w,e_i) + C(v,w,\Delta^{-1}(0)J_1e_i).
\end{align*}
Then the equations for $H_{20\mu}$ and $H_{11\mu}$ become
\begin{align*}
     H_{20\mu}(\theta) &= \gamma_{1,\mu} e^{2i\omega_0\theta}\Delta^{-1}(2i\omega_0)\Lambda_1(H_{2000},\varphi,\varphi) + \gamma_{2,\mu} e^{2i\omega_0\theta}\Delta^{-1}(2i\omega_0)\Lambda_2(H_{2000},\varphi,\varphi) \\
     &- 2\delta_{\mu}^{10}\Delta^{-1}(2i\omega_0)\left([\Delta'(2i\omega_0) - \theta\Delta(2i\omega_0)]H_{2000}(\theta) + e^{2i\omega_0\theta}B\left(\varphi, B_{i\omega_0}^{Inv}(0,q) \right)\right),\\ 
     &-2ib_{1,\mu}\Delta^{-1}(2i\omega_0)\left([\Delta'(2i\omega_0) - \theta\Delta(2i\omega_0)]H_{2000}(\theta) + e^{2i\omega_0\theta}B\left(\varphi, B_{i\omega_0}^{Inv}(0,q) \right)\right),\\
    H_{11\mu}(\theta) &=  \gamma_{1,\mu}\Delta^{-1}(0)\Pi_1(H_{1100},\bar{\varphi},\varphi)  +\gamma_{2,\mu}\Delta^{-1}(0)\Pi_2(H_{1100},\bar{\varphi},\varphi) \\
    &-2\delta_{\mu}^{10}\Delta^{-1}(0)\left([\Delta'(0) - \theta\Delta(0)]H_{1100}(\theta) +\Re{}\left\{ B\left(\bar{\varphi}, B_{i\omega_0}^{Inv}(0,q)\right)\right\} \right)\\
    &-2b_{1,\mu}\Delta^{-1}(0)\Re{}\left\{ iB\left(\bar{\varphi}, B_{i\omega_0}^{Inv}(0,q)\right)\right\},
\end{align*}
The $z^2\bar{z}\beta_i$ terms yield the systems
\begin{align}\label{eq:AH21mu_DDE}
    (i\omega_0I_n- A^{\odot\star})jH_{21\mu} &=  [A_1(H_{2100},K_{\mu}) + 2B(\varphi,H_{11\mu}) + B(\bar{\varphi},H_{20\mu}) + B(H_{00\mu},H_{2100}) \nonumber\\
    &+ B(\overline{H}_{10\mu},H_{2000}) + 2B(H_{10\mu},H_{1100}) + 2B_1(\varphi,H_{1100},K_{\mu})  \nonumber \\
    &+ B_1(\bar{\varphi},H_{2000},K_{\mu})+ C(\varphi,\varphi,\overline{H}_{10\mu}) + 2C(\varphi,\bar{\varphi},H_{10\mu})  \nonumber\\
    &+ 2C(\varphi,H_{00\mu},H_{1100})+ C(\bar{\varphi},H_{00\mu},H_{2000}) + C_1(\varphi,\varphi,\bar{\varphi},K_{\mu})\nonumber\\
    & + D(\varphi,\varphi,\bar{\varphi},H_{00\mu})]r^{\odot\star}-[2(\delta_{\mu}^{01} + ib_{2,\mu})j\varphi + (3\delta_{\mu}^{10} + ib_{1,\mu})jH_{2100}  \nonumber\\&+2c_1(0)jH_{10\mu}].
\end{align}
Applying the Fredholm Alternative to equation \cref{eq:AH21mu_DDE} results in the equation
\begin{align}\label{eq:H21mu_Fred_DDE}
    \delta_{\mu}^{01} + ib_{2,\mu} &= \frac{1}{2}p^T [A_1(H_{2100},K_{\mu}) + 2B(\varphi,H_{11\mu}) + B(\bar{\varphi},H_{20\mu}) + B(H_{00\mu},H_{2100}) \nonumber\\
    &+ B(\overline{H}_{10\mu},H_{2000}) + 2B(H_{10\mu},H_{1100}) + 2B_1(\varphi,H_{1100},K_{\mu}) + B_1(\bar{\varphi},H_{2000},K_{\mu}) \nonumber\\
    &+ C(\varphi,\varphi,\overline{H}_{10\mu}) + 2C(\varphi,\bar{\varphi}.H_{10\mu}) + 2C(\varphi,H_{00\mu},H_{1100}) + C(\bar{\varphi},H_{00\mu},H_{2000}) \nonumber\\
    &+ C_1(\varphi,\varphi,\bar{\varphi},K_{\mu}) + D(\varphi,\varphi,\bar{\varphi},H_{00\mu})],
\end{align}
where we used that $\langle \varphi^{\odot}, H_{2100} \rangle = 0$ and $\langle \varphi^{\odot}, H_{10\mu} \rangle = 0$. If we substitute the expressions for $H_{00\mu}, H_{10\mu},H_{20\mu}, H_{11\mu}$, $K_{\mu}$ and $b_{1,\mu}$ into equation \cref{eq:H21mu_Fred_DDE}, we can solve for $\gamma_{1,\mu}$ and $\gamma_{2,\mu}$.

After substitution and some rewriting we arrive at the following system 
\begin{equation}\label{eq:SPQ1_DDE}
    P\begin{pmatrix}
        \gamma_{1,\mu} \\
        \gamma_{2,\mu}
    \end{pmatrix} = Q_{\mu},
\end{equation}
where $P \in \R^{2\times2}$ has for $k = 1,2$ the compenents
\begin{align}
    P_{1k} &= \Re{[p^T\Gamma_k(\varphi)]}, \label{eq:SP1k_DDE}\\
    P_{2k} &= \frac{1}{2}\Re{}\Bigl\{ p^T\Bigl[ \Gamma_k(H_{2100}) + 2B\left(\varphi,\Delta^{-1}(0)\Pi_k(H_{1100},\bar{\varphi},\varphi)\right) \nonumber\\
    &+ B\left(\bar{\varphi},e^{2i\omega_0\theta}\Delta^{-1}(2i\omega_0)\Lambda_k(H_{2000},\varphi,\varphi)\right) +B\left(H_{2000},\overline{B_{i\omega_0}^{Inv}(\Gamma_k(\varphi))}\right) \nonumber\\
    &+2B\left(H_{1100},B_{i\omega_0}^{Inv}(\Gamma_k(\varphi))\right) + 2B_1(\varphi,H_{1100},e_k)   \nonumber\\
    &+ B_1(\bar{\varphi},H_{2000},e_k) +C\left(\varphi,\varphi,\overline{B_{i\omega_0}^{Inv}(\Gamma_k(\varphi))}\right) + 2C\left(\varphi,\bar{\varphi},B_{i\omega_0}^{Inv}(\Gamma_k(\varphi))\right) \nonumber\\
    &+2C\left(\varphi,H_{1100},\Delta^{-1}(0)J_1e_k\right) +C(\bar{\varphi},H_{2000},\Delta^{-1}(0)J_1e_k) + C_1(\varphi,\varphi,\bar{\varphi},e_k)   \nonumber\\
    &+ D(\varphi,\varphi,\bar{\varphi},\Delta^{-1}(0)J_1e_k) + \Im{[p^T\Gamma_k(\varphi)]}\Bigl(-4B\left(\varphi,\Delta^{-1}(0)\Re{[iB\left(\bar{\varphi},B_{i\omega_0}^{Inv}(0,q)\right)]} \right) \nonumber\\
    &-2iB\left(\bar{\varphi},e^{2i\omega_0\theta}\Delta^{-1}(2i\omega_0)\left([\Delta'(2i\omega_0) - \theta\Delta(2i\omega_0)]H_{2000}(0) + B\left(\varphi, B_{i\omega_0}^{Inv}(0,q) \right)  \right)\right)  \nonumber\\
    &\left.\left. \left.+iB\left(H_{2000},\overline{B_{i\omega_0}^{Inv}(0,q)}\right) - 2iB\left(H_{1100},B_{i\omega_0}^{Inv}(0,q)\right) + iC\left(\varphi,\varphi,\overline{B_{i\omega_0}^{Inv}(0,q)}\right)\right.\right.\right. \nonumber\\
    &- 2iC\left(\varphi,\bar{\varphi},B_{i\omega_0}^{Inv}(0,q)\right) \Bigr)\Bigr] \Bigr\},\label{eq:SP2k_DDE}
\end{align}
 and $Q_{\mu} \in \R^2$ with $\mu = (10),(01)$ is given by
\begin{align}
   Q_{1,\mu} &= \delta_{\mu}^{10}, \label{eq:SQ1mu1_DDE}\\
    Q_{2,\mu} &= \frac{1}{2}\delta_{\mu}^{10} \Re{}\Bigl\{ p^T\Bigl[ 4B\left(\varphi,\Delta^{-1}(0)\left([\Delta'(0) - \theta\Delta(0)]H_{1100}(\theta) + \Re{}\left\{ B\left(\bar{\varphi}, B_{i\omega_0}^{Inv}(0,q)\right)\right\}\right) \right) \nonumber\\
    &+2B\left(\bar{\varphi},e^{2i\omega_0\theta}\Delta^{-1}(2i\omega_0)\left([\Delta'(2i\omega_0) - \theta\Delta(2i\omega_0)]H_{2000}(0) + B\left(\varphi, B_{i\omega_0}^{Inv}(0,q) \right) \right) \right)  \nonumber\\
    &+B\left(H_{2000},\overline{B_{i\omega_0}^{Inv}(0,q)}\right) + 2B\left(H_{1100},B_{i\omega_0}^{Inv}(0,q)\right) + C\left(\varphi,\varphi,\overline{B_{i\omega_0}^{Inv}(0,q)}\right) \nonumber\\
    &+ 2C\left(\varphi,\bar{\varphi},B_{i\omega_0}^{Inv}(0,q)\right)  \Bigr]\Bigr\} +\delta_{\mu}^{01} . \label{eq:SQ2mu1_DDE}
\end{align}

To summarise, we first compute the linear coefficients $K_{10}$ and $K_{01}$ by solving system \cref{eq:SPQ1_DDE}. Then we can determine the coefficients $H_{00\mu}$ from equation \cref{eq:H00muDDE}, $b_{1,\mu}$ from equation \cref{eq:b1mu1DDE} and $H_{10\mu}$ from equation \cref{S-eq:H10muDDE}. Once we have these, the coefficients $H_{20\mu}$ and $H_{11\mu}$ can respectively be calculated from equations \cref{eq:H20mu1_DDE} and \cref{eq:H11mu1_DDE}. The coefficient $b_{2,01}$ is determined from the imaginary part of equation \cref{eq:H21mu_Fred_DDE}. Finally, we will also need $H_{2101}$ and $H_{2110}$ for what follows. The solution to equation \cref{eq:AH21mu_DDE} can be found by applying \cref{cor:bordered_case}. Denote the part containing all the multilinear forms in equation \cref{eq:AH21mu_DDE} by $M_{21\mu}$. This results in the equation
\begin{align*}
     H_{21\mu}(\theta) &= B_{i\omega_0}^{INV}\Bigl(M_{21\mu}, [-2(\delta_{\mu}^{01} + ib_{2,\mu})q- (3\delta_{\mu}^{10} + ib_{1,\mu})H_{2100}(0)- 2c_1(0)H_{10\mu}(0)],\\
     &-4(2\delta_{\mu}^{10}+ib_{1,\mu})c_1(0)q\Bigr)(\theta), \quad \mu = (01),(10).
 \end{align*}

\paragraph{The coefficient $a_{3201}$}
As in the ODE case, to determine the coefficient $a_{3201}$ we first need the coefficients $H_{3001}, H_{3101}$ and $H_{2201}$. These can be found by collecting the $z^3\beta_2$, $z^3\bar{z}\beta_2$ and $z^2\bar{z}^2\beta_2$ terms from the homological equation \cref{eq:HOM_GH_DDE}. For convenience, we will write $M_{ijkl}$ for the term containing all the multilinear forms as these are the same as in the corresponding equations for ODEs after changing $q$ to $\varphi$. The equations for $H_{3001}$ and $H_{2201}$ follow from \cref{cor:inverse_case1}. This results in the expressions
 \begin{align*}
     H_{3001}(\theta) &=  e^{3i\omega_0\theta}\Delta^{-1}(3i\omega_0)M_{3001} - 3ib_{1,01}\Delta^{-1}(3i\omega_0)[\Delta'(3i\omega_0) -\theta\Delta(3i\omega_0)]H_{3000}(\theta), \nonumber\\
     H_{2201}(\theta) &= \Delta^{-1}(0)M_{2201} - 8\Delta^{-1}(0)[\Delta'(0) - \theta\Delta(0)]H_{1100}(\theta) .
 \end{align*}
 For $H_{3101}$, we need to be more careful. The $z^3\bar{z}\beta_2$ terms yield the equation
 \begin{align*}\label{eq:H3101}
    (2i\omega_0I_n- A^{\odot\star})jH_{3101} &= M_{3101}r^{\odot\star}-6(1 + ib_{2,01})jH_{2000} - 6c_1(0)jH_{2001} - 2ib_{1,01}jH_{3100}.
 \end{align*}
 To solve for $H_{3101}$, we can separately apply \cref{cor:inverse_case1} for the part $M_{3101}r^{\odot\star}-6(1 + ib_{2,01})jH_{2000}$ and Corollary \cref{cor:inverse_case2} for the last two terms. This yields the following solution
 \begin{align}
     H_{3101}(\theta) &= e^{2i\omega_0\theta}\Delta^{-1}(2i\omega_0)M_{3101} - 6(1 + ib_{2,01})\Delta^{-1}(2i\omega_0)[\Delta'(2i\omega_0) - \theta\Delta(2i\omega_0)]H_{2000}(\theta) \nonumber \\
     &- 6c_1(0)e^{2i\omega_0\theta}\Delta^{-1}(2i\omega_0) \Bigl([\Delta'(2i
    \omega_0) - \theta\Delta(2i\omega_0)]H_{2001}(0) \nonumber\\
    &+ ib_{1,01} [\Delta''(2i\omega_0)- \theta^2\Delta(2i\omega_0)]H_{2000}(0)\Bigr) \nonumber \\
    &-2ib_{1,01}e^{2i\omega_0\theta}\Delta^{-1}(2i\omega_0) \Bigl([\Delta'(2i
    \omega_0) - \theta\Delta(2i\omega_0)]H_{3100}(0) \nonumber\\
    &+ 3c_1(0)[\Delta''(2i\omega_0)- \theta^2\Delta(2i\omega_0)]H_{2000}(0)\Bigr).
 \end{align}
 Finally, collecting the $z^3\bar{z}^2\beta_2$ terms results in the following equation
 \begin{align}\label{eq:AH3201_DDE}
     (i\omega_0I_n- A^{\odot\star})jH_{3201} &= M_{3201}r^{\odot\star} - [12g_{3201}j\varphi + 12c_2(0)jH_{1001}  +(18 + 6ib_{2,01})jH_{2100} \nonumber\\
     &+ 6i\Im{}\{c_1(0)\}jH_{2101} + ib_{1,01}jH_{3200}].
 \end{align}
The coefficient $a_{3201}$ can now be found by applying the Fredholm solvability condition to equation \cref{eq:AH3201_DDE} and then taking the real part. Thus
\begin{align}
    a_{3201} = \frac{1}{12}\Re{}\left\{ p^T M_{3201} \right\}.
\end{align}
\paragraph{Quadratic coefficients $K_{02},K_{11},b_{1,02}$}
In the ODE case, we put all of the terms that contained known coefficients into a rest term $r_{ijkl}$. Since the multilinear forms in $r_{ijkl}$ will be the same for DDEs but with $q$ changed to $\varphi$, we will indicate those with $M_{ijkl}$. Collecting the $\beta_2^2$ and $\beta_1\beta_2$ terms from \cref{eq:HOM_GH_DDE} yields for $\mu = (02),(11)$ the equations
\begin{align}\label{eq:AH0002_DDE}
    -A^{\odot\star}jH_{00\mu} &= [J_1K_{\mu}+ M_{00\mu}]r^{\odot\star}.
\end{align}
Let $e_1,e_2 \in \R^2$ be the standard basis vectors and write
\begin{align}\label{eq:K02_DDE}
    K_{\mu} = \gamma_{1,\mu}e_1 + \gamma_{2,\mu}e_2,
\end{align}
where $\gamma_{1,\mu}, \gamma_{1,\mu} \in \R$ are unknown constants that need to be determined. From equation \cref{eq:AH0002_DDE} it follows after applying \cref{cor:inverse_cases}.1 that
\begin{align}\label{eq:H0002_DDE}
    H_{00\mu}(\theta) = \gamma_{1,\mu}\Delta^{-1}(0)J_1e_1 + \gamma_{2,\mu}\Delta^{-1}(0)J_1e_2 + \Delta^{-1}(0)M_{00\mu}.
\end{align}
The $z\beta_2^2$ and $z\beta_1\beta_2$ terms yield for $\mu = (02),(11)$ the equations
\begin{align}
    (i\omega_0I_n- A^{\odot\star})jH_{10\mu} &= [A_1(\varphi,K_{\mu}) +  B(\varphi,H_{00\mu})]r^{\odot\star} - ib_{1,\mu}j\varphi + r_{10\mu}, \label{eq:AH10mu_DDE} 
\end{align}
where
\begin{align*}
    r_{1002} &= M_{1002}r^{\odot\star} -2ib_{1,01}jH_{1001}, \\
    r_{1011} &= M_{1011}r^{\odot\star} -(1 + ib_{1,10})jH_{1001} - ib_{1,01}jH_{1010}.
\end{align*}
Applying the Fredholm alternative to equations \cref{eq:AH10mu_DDE}  yields the equations 
\begin{align}\label{eq:ibmu2_DDE}
    ib_{1,\mu} &= p^T[A_1(\varphi,K_{\mu}) +  B(\varphi,H_{00\mu}) + M_{10\mu}].
\end{align}
Here we used that $\langle \varphi^{\odot}, H_{1010} \rangle = \langle \varphi^{\odot}, H_{1001} \rangle = 0$. Substituting equations \cref{eq:K02_DDE} and \cref{eq:H0002_DDE} into equation \cref{eq:ibmu2_DDE} results in 
\begin{align}\label{eq:ib102_DDE}
    ib_{1,\mu} &= p^T[\gamma_{1,\mu}\Gamma_1(\varphi) + \gamma_{2,\mu}\Gamma_2(\varphi)  +  \tilde{M}_{10\mu}],
\end{align}
where 
\begin{align*}
    \tilde{M}_{10\mu} &= B(\varphi,\Delta^{-1}(0)M_{00\mu}) + M_{10\mu}.
\end{align*}
From equation \cref{eq:ib102_DDE}, it follows that 
\begin{align}
     \gamma_{1,\mu}\Re{[p^T\Gamma_1(\varphi)]} + \gamma_{2,\mu}\Re{[p^T\Gamma_2(\varphi)]} = -\Re{}\{p^T \tilde{M}_{10\mu}\},
\end{align}
and
\begin{align}\label{eq:b102_DDE}
    b_{1,\mu} = \gamma_{1,\mu}\Im{[p^T\Gamma_1(\varphi)]} + \gamma_{2,\mu}\Im{[p^T\Gamma_2(\varphi)]} + \Im{}\{p^T \tilde{M}_{10\mu}\}.
\end{align}
Furthermore, a piecewise application of \cref{cor:bordered_case} to equations \cref{eq:AH10mu_DDE} yields for $\mu = (02),(11)$ the solutions
\begin{align}\label{eq:H10mu2_DDE}
    H_{10\mu}(\theta) &=  \gamma_{1,\mu}B_{i\omega_0}^{Inv}(\Gamma_1(\varphi))(\theta)  + \gamma_{2,\mu}B_{i\omega_0}^{Inv}(\Gamma_2(\varphi))(\theta) \nonumber\\
    &- ib_{1,\mu}B_{i\omega_0}^{Inv}(0,q)(\theta)  + B_{10\mu}(\theta),
\end{align}
where
\begin{align*}
    B_{1002}(\theta) &= B_{i\omega_0}^{Inv}(\tilde{M}_{1002})(\theta) - 2ib_{1,01}B_{i\omega_0}^{Inv}(0, H_{1001}(0),ib_{1,01}q)(\theta), \\
    B_{1011}(\theta) &=  B_{i\omega_0}^{Inv}(\tilde{M}_{1011})(\theta) -(1 + ib_{1,10})B_{i\omega_0}^{Inv}(0, H_{1001}(0),ib_{1,01}q)(\theta) \\
    &-ib_{1,01}B_{i\omega_0}^{Inv}(0,H_{1010}(0),(1 + ib_{1,10})q)(\theta).\
\end{align*}
Now collect the $z^2\beta^{\mu}$  terms for $\mu = (02),(11)$. These terms yield the equations
\begin{align}
    (2i\omega_0I_n- A^{\odot\star})jH_{20\mu} &= [A_1(H_{2000},K_{\mu}) + 2B(\varphi,H_{10\mu})+ B(H_{00\mu},H_{2000}) \nonumber \\\label{eq:AH20mu2_DDE}
    &+  B_1(\varphi,\varphi,K_{\mu})+ C(\varphi,\varphi,H_{00\mu})]r^{\odot\star}  - 2ib_{1,\mu}jH_{2000} + r_{20\mu},
\end{align}
where
\begin{align*}
    r_{2002} &= M_{2002}r^{\odot\star} - 4ib_{1,01}jH_{2001}, \\
    r_{2011} &= M_{2011}r^{\odot\star} - 2(1 + ib_{1,10})jH_{2001} - 2ib_{1,01}jH_{2010}.
\end{align*}
Using \cref{cor:inverse_cases} and \cref{cor:inverse_case2}, we find that the solutions to equations \cref{eq:AH20mu2_DDE} are
\begin{align}
    H_{2002}(\theta) &= e^{2i\omega_0\theta}\Delta^{-1}(2i\omega_0)[A_1(H_{2000},K_{02}) + 2B(\varphi,H_{1002})+ B(H_{0002},H_{2000}) \nonumber \\
    &+  B_1(\varphi,\varphi,K_{02})+ C(\varphi,\varphi,H_{0002}) +M_{2002}] \nonumber\\
    &-2ib_{1,02}\Delta^{-1}(2i\omega_0)[\Delta'(2i\omega_0) - \theta\Delta(2i\omega_0)]H_{2000}(\theta) \nonumber \\
    &-4ib_{1,01}e^{2i\omega_0\theta}\Delta^{-1}(2i\omega_0) \Bigl([\Delta'(2i
    \omega_0) - \theta\Delta(2i\omega_0)]H_{2001}(0) \nonumber\\
    &+ ib_{1,01} [\Delta''(2i\omega_0)- \theta^2\Delta(2i\omega_0)]H_{2000}(0)\Bigr), \label{eq:H2002_DDE}\\
    H_{2011}(\theta) &= e^{2i\omega_0\theta}\Delta^{-1}(2i\omega_0)[A_1(H_{2000},K_{11}) + 2B(\varphi,H_{1011})+ B(H_{0011},H_{2000}) \nonumber \\
    &+  B_1(\varphi,\varphi,K_{11})+ C(\varphi,\varphi,H_{0011}) +M_{2011}] \nonumber\\
    &-2ib_{1,11}\Delta^{-1}(2i\omega_0)[\Delta'(2i\omega_0) - \theta\Delta(2i\omega_0)]H_{2000}(\theta) \nonumber \\
    &-2(1 + ib_{1,10})e^{2i\omega_0\theta}\Delta^{-1}(2i\omega_0) \Bigl([\Delta'(2i
    \omega_0) - \theta\Delta(2i\omega_0)]H_{2001}(0) \nonumber\\
    &+ ib_{1,01} [\Delta''(2i\omega_0)- \theta^2\Delta(2i\omega_0)]H_{2000}(0)\Bigr) \nonumber \\
    &- 2ib_{1,01}e^{2i\omega_0\theta}\Delta^{-1}(2i\omega_0)\Bigl( [\Delta'(2i
    \omega_0) - \theta\Delta(2i\omega_0)]H_{2010}(0) \nonumber\\
    &+ (1 + ib_{1,10}) [\Delta''(2i\omega_0)- \theta^2\Delta(2i\omega_0)]H_{2000}(0)\Bigr). \label{eq:H2011_DDE}
\end{align}
Substituting equations \cref{eq:K02_DDE},\cref{eq:H0002_DDE} and \cref{eq:H10mu2_DDE} into the above expressions result in 
\begin{align}\label{eq:H20mu2_DDE}
    H_{20\mu}(\theta) &= \gamma_{1,\mu}e^{2i\omega_0\theta}\Delta^{-1}(2i\omega_0)\Lambda_1(H_{2000}, \varphi, \varphi) + \gamma_{2,\mu}e^{2i\omega_0\theta}\Delta^{-1}(2i\omega_0)\Lambda_2(H_{2000}, \varphi, \varphi) \nonumber \\
    &-2ib_{1,\mu}\Delta^{-1}(2i\omega_0)\left([\Delta'(2i\omega_0) - \theta\Delta(2i\omega_0)]H_{2000}(\theta) + e^{2i\omega_0\theta}B\left(\varphi, B_{i\omega_0}^{Inv}(0,q)\right)\right) \nonumber \\
    &+ e^{2i\omega_0\theta}\Delta^{-1}(2i\omega_0)\tilde{r}_{20\mu}(\theta),
\end{align}
where we have for $\mu = (02),(11)$:
\begin{align*}
    \tilde{r}_{2002}(\theta) &= M_{2002} + 2B(\varphi, B_{1002})  + B(H_{2000}, \Delta^{-1}(0)M_{0002}) + C(\varphi,\varphi, \Delta^{-1}(0)M_{0002}) \\
    &-4ib_{1,01}\Bigl([\Delta'(2i
    \omega_0) - \theta\Delta(2i\omega_0)]H_{2001}(0) \nonumber\\
    &+ ib_{1,01}[\Delta''(2i\omega_0)- \theta^2\Delta(2i\omega_0)]H_{2000}(0)\Bigr), \\
    \tilde{r}_{2011}(\theta) &= M_{2011} + 2B(\varphi, B_{1011})  + B(H_{2000}, \Delta^{-1}(0)M_{0011}) + C(\varphi,\bar{\varphi}, \Delta^{-1}(0)M_{0011}) \\
    &-2(1 + ib_{1,10})\Bigl([\Delta'(2i
    \omega_0) - \theta\Delta(2i\omega_0)]H_{2001}(0) \nonumber\\
    &+ ib_{1,01} [\Delta''(2i\omega_0)- \theta^2\Delta(2i\omega_0)]H_{2000}(0)\Bigr) \nonumber \\
    &- 2ib_{1,01}\Bigl( [\Delta'(2i
    \omega_0) - \theta\Delta(2i\omega_0)]H_{2010}(0) \nonumber\\
    &+ (1 + ib_{1,10}) [\Delta''(2i\omega_0)- \theta^2\Delta(2i\omega_0)]H_{2000}(0)\Bigr).
\end{align*}
Similarly, collecting the $z\bar{z}\beta^{\mu}$ terms for $\mu = (02),(11)$ yield the equations
\begin{align}
    -A^{\odot\star}jH_{11\mu} &= [A_1(H_{1100},K_{\mu}) + 2\Re{}\left(B(\bar{\varphi},H_{10\mu})\right)+ B(H_{00\mu},H_{1100}) \nonumber\\
    &+ B_1(\varphi,\bar{\varphi},K_{\mu})+ C(\varphi,\bar{\varphi},H_{00\mu})]r^{\odot\star} + r_{11\mu}, \label{eq:HA1102_DDE} 
\end{align}
where
\begin{align*}
    r_{1102} = M_{1102}r^{\odot\star} \quad \textrm{ and } \quad r_{1111} = M_{1111}r^{\odot\star} - 2jH_{1101}.
\end{align*}
Using \cref{cor:inverse_case1}, we find that the solutions to equations \cref{eq:HA1102_DDE} are
\begin{align}
    H_{1102}(\theta) &= \Delta^{-1}(0)[A_1(H_{1100},K_{02}) + 2\Re{}\left(B(\bar{\varphi},H_{1002})\right)+ B(H_{0002},H_{1100}) \nonumber\\
    &+ B_1(\varphi,\bar{\varphi},K_{02})+ C(\varphi,\bar{\varphi},H_{0002})+ M_{1102}], \label{eq:H1102_DDE}\\
    H_{1111}(\theta) &= \Delta^{-1}(0)[A_1(H_{1100},K_{11}) + 2\Re{}\left(B(\bar{\varphi},H_{1011})\right)+ B(H_{0011},H_{1100}) \nonumber\\
    &+ B_1(\varphi,\bar{\varphi},K_{11})+ C(\varphi,\bar{\varphi},H_{0011})+ M_{1111}] \nonumber \\
    &-2\Delta^{-1}(0)[\Delta'(0) - \theta\Delta(0)]H_{1101}(\theta).\label{eq:H1111_DDE}
\end{align}
Substituting equations \cref{eq:K02_DDE},\cref{eq:H0002_DDE} and \cref{eq:H10mu2_DDE} into the above expressions result in 
\begin{align}\label{eq:H11mu_DDE}
    H_{11\mu}(\theta) &= \gamma_{1,
    \mu}\Delta^{-1}(0)\Pi_1(H_{1100},\bar{\varphi},\varphi) +  \gamma_{2,
    \mu}\Delta^{-1}(0)\Pi_1(H_{1100},\bar{\varphi},\varphi) \nonumber\\
    &-2b_{1,\mu}\Delta^{-1}(0)\Re{}\left\{ iB(\left(\bar{\varphi},B_{i\omega_0}^{Inv}(0,q)\right)\right\} + \Delta^{-1}(0)\tilde{r}_{11\mu}(\theta),
\end{align}
where $\tilde{r}_{11\mu}$ is given for $\mu = (02),(11)$ by 
\begin{align*}
    \tilde{r}_{1102}(\theta) &= M_{1102} + 2\Re{}\{ B(\bar{\varphi}, B_{1002}) \} + B(H_{1100}, \Delta^{-1}(0)M_{0002}) + C(\varphi,\bar{\varphi}, \Delta^{-1}(0)M_{0002}), \\
    \tilde{r}_{1111}(\theta) &= M_{1111} + 2\Re{}\{ B(\bar{\varphi}, B_{1011}) \} + B(H_{1100}, \Delta^{-1}(0)M_{0011}) + C(\varphi,\bar{\varphi}, \Delta^{-1}(0)M_{0011}) \\
    &-2[\Delta'(0) - \theta\Delta(0)]H_{1101}(\theta).
\end{align*}
The $z^2\bar{z}\beta^{\mu}$ terms yield for $\mu = (02),(11)$ the equations
\begin{align}\label{eq:AH21mu2_DDE}
    (i\omega_0I_n- A^{\odot\star})jH_{21\mu} &= [A_1(H_{2100},K_{\mu}) + 2B(\varphi,H_{11\mu})+ B(\bar{\varphi},H_{20\mu})+ B(H_{00\mu},H_{2100}) \nonumber\\
    &+ B(\overline{H}_{10\mu},H_{2000})+ 2B(H_{10\mu},H_{1100})+ 2B_1(\varphi,H_{1100},K_{\mu}) \nonumber\\
    &+ B_1(\bar{\varphi},H_{2000},K_{\mu})+ C(\varphi,\varphi,\overline{H}_{10\mu})+ 2C(\varphi,\bar{\varphi},H_{10\mu})\nonumber\\
    &+ 2C(\varphi,H_{00\mu},H_{1100}) + C(\bar{\varphi},H_{00\mu},H_{2000}) + C_1(\varphi,\varphi,\bar{\varphi},K_{\mu})\nonumber\\
    &+ D(\varphi,\varphi,\bar{\varphi},H_{00\mu}) ]r^{\odot\star} + r_{21\mu} \nonumber\\
    &- (2ib_{2,\mu}j\varphi +ib_{1,\mu}jH_{2100} + 2c_1(0)jH_{10\mu}) ,
\end{align}
where
\begin{align*}
    r_{2102} &= M_{2102}r^{\odot\star}- [4(1 + ib_{2,01})jH_{1001}  + 2ib_{1,01}jH_{2101}], \\
    r_{2111} &= M_{2111}r^{\odot\star}- [2ib_{2,10}jH_{1001}+2 (1+ib_{2,01})jH_{1010}+ \left(3 + ib_{1,10}\right) jH_{2101}+ib_{1,01}jH_{2110}].
\end{align*}
Applying the Fredholm solvability condition to the equation \cref{eq:AH21mu2_DDE} results for $\mu = (02),(11)$ in the equations
\begin{align}\label{eq:b2mu2_DDE}
    b_{2,\mu}i &= \frac{1}{2}p^T[A_1(H_{2100},K_{\mu}) + 2B(\varphi,H_{11\mu})+ B(\bar{\varphi},H_{20\mu})+ B(H_{00\mu},H_{2100}) \nonumber\\
    &+ B(\overline{H}_{10\mu},H_{2000})+ 2B(H_{10\mu},H_{1100})+ 2B_1(\varphi,H_{1100},K_{\mu}) \nonumber\\
    &+ B_1(\bar{\varphi},H_{2000},K_{\mu})+ C(\varphi,\varphi,\overline{H}_{10\mu})+ 2C(\varphi,\bar{\varphi},H_{10\mu})+ 2C(\varphi,H_{00\mu},H_{1100}) \nonumber\\
    &+ C(\bar{\varphi},H_{00\mu},H_{2000}) + C_1(\varphi,\varphi,\bar{\varphi},K_{\mu})+ D(\varphi,\varphi,\bar{\varphi},H_{00\mu}) + M_{21\mu}],
\end{align}
\\[3mm]Finally, substitute equations \cref{eq:K02_DDE},\cref{eq:H0002_DDE}, \cref{eq:H10mu2_DDE}, \cref{eq:H20mu2_DDE} and \cref{eq:H11mu_DDE} into the above expression and solve for the coefficients $\gamma_{1,\mu}$,$\gamma_{2,\mu}$. We arrive at the following system

\begin{equation}\label{eq:PQ2_DDE}
    P\begin{pmatrix}
        \gamma_{1,\mu} \\
        \gamma_{2,\mu}
    \end{pmatrix} = Q_{\mu},
\end{equation}
where the components of $P \in \R^{2\times2}$ are given by equations \cref{eq:SP1k_DDE} and \cref{eq:SP2k_DDE}. Meanwhile, the components of $Q_{\mu} \in \R^2$, where $\mu = (02),(11)$, are
\begin{align*}
   Q_{1,\mu} &=  -\Re{}\{p^T \tilde{M}_{10\mu}\}, \\
    Q_{2,\mu} &= -\frac{1}{2}\Re{}\Bigl\{ p^T\Bigl[ 2B\left(\varphi, \Delta^{-1}(0)\tilde{r}_{11\mu}\right) + B\left(\bar{\varphi},e^{2i\omega_0\theta}\Delta^{-1}(2i\omega_0)\tilde{r}_{20\mu}\right)\\
    &+B\left(H_{2100},\Delta^{-1}(0)M_{00\mu}\right) + B\left(H_{2000},\overline{B_{10\mu}}\right) + 2B\left(H_{1100},B_{10\mu}\right) \\
    &+C\left(\varphi,\varphi,\overline{B_{10\mu}}\right) + 2C\left(\varphi,\bar{\varphi},B_{10\mu}\right) + 2C\left(\varphi,H_{1100},\Delta^{-1}(0)M_{00\mu}\right) \\
    &+ C\left(\bar{\varphi},H_{2000},\Delta^{-1}(0)M_{00\mu}\right) + D\left(\varphi, \varphi, \bar{\varphi}, \Delta^{-1}(0)M_{00\mu}\right) + M_{21\mu} \\
    &+\Im{}\{p^T \tilde{M}_{10\mu}\}\Bigl(-4B\left(\varphi,\Delta^{-1}(0)\Re{[iB\left(\bar{\varphi},B_{i\omega_0}^{Inv}(0,1)\right)]} \right) \nonumber\\
    &-2iB\left(\bar{\varphi},e^{2i\omega_0\theta}\Delta^{-1}(2i\omega_0)\left([\Delta'(2i\omega_0) - \theta\Delta(2i\omega_0)]H_{2000}(0) + B\left(\varphi, B_{i\omega_0}^{Inv}(0,q) \right)  \right)\right)  \nonumber\\
    &\left.\left. \left.+iB\left(H_{2000},\overline{B_{i\omega_0}^{Inv}(0,q)}\right) - 2iB\left(H_{1100},B_{i\omega_0}^{Inv}(0,q)\right) + iC\left(\varphi,\varphi,\overline{B_{i\omega_0}^{INV}(0,q)}\right)\right.\right.\right. \nonumber\\
    &- 2iC\left(\varphi,\bar{\varphi},B_{i\omega_0}^{Inv}(0,q)\right)) \Bigr)\Bigr]\Bigr\}.
\end{align*}
 Thus, to obtain the quadratic coefficients $K_{02}$ and $K_{11}$ we solve system \cref{eq:PQ2_DDE}. Once these coefficients are known, one can calculate the coefficients $H_{00\mu}$,$b_{1,\mu}$ and $H_{10\mu}$ from respectively the equations \cref{eq:H0002_DDE},\cref{eq:b102_DDE} and \cref{eq:H10mu2_DDE}. Then, the coefficients $H_{2002}$ and $H_{1102}$ are calculated from respectively the equations \cref{eq:H2002_DDE} and \cref{eq:H1102_DDE}. Finally, the coefficient $H_{2102}$, needs to be determined from system \cref{eq:AH21mu2_DDE}. This can be achieved by a piecewise application of \cref{cor:bordered_case}.
  \begin{align*}
     H_{2102}&= B_{i\omega_0}^{Inv}\Bigl([A_1(H_{2100},K_{\mu}) + 2B(\varphi,H_{11\mu})+ B(\bar{\varphi},H_{20\mu})+ B(H_{00\mu},H_{2100}) \nonumber\\
    &+ B(\overline{H}_{10\mu},H_{2000})+ 2B(H_{10\mu},H_{1100})+ 2B_1(\varphi,H_{1100},K_{\mu}) \nonumber\\
    &+ B_1(\bar{\varphi},H_{2000},K_{\mu})+ C(\varphi,\varphi,\overline{H}_{10\mu})+ 2C(\varphi,\bar{\varphi},H_{10\mu})+ 2C(\varphi,H_{00\mu},H_{1100}) \nonumber\\
    &+ C(\bar{\varphi},H_{00\mu},H_{2000}) + C_1(\varphi,\varphi,\bar{\varphi},K_{\mu})+ D(\varphi,\varphi,\bar{\varphi},H_{00\mu}) + M_{2102}], -2ib_{2,02}q\Bigr) \\
    &- ib_{1,02}B_{i\omega_0}^{Inv}\Bigl(0,H_{2100}(0), 2c_1(0)q) -4(1 + ib_{2,01})B_{i\omega_0}^{Inv}(0,H_{1001}(0),ib_{1,01}q\Bigr) \\
     &-2c_1(0)B_{i\omega_0}^{Inv}\Bigl(0,H_{1002}(0), [ib_{1,02}q + 2ib_{1,01}H_{1001}(0)], -b_{1,01}^2q\Bigr) \\
     &-2ib_{1,01}B_{i\omega_0}^{Inv}\Bigl(0,H_{2101}(0), [2(1 + ib_{2,01})q + ib_{1,01}H_{2100}(0) \\
     &+ 2c_1(0)H_{1001}(0)], 2ib_{1,01}c_1(0)q\Bigr).
 \end{align*}

\paragraph{The coefficient $K_{03}$}
Collecting the $\beta_2^3$ terms yields the equation
 \begin{align}\label{eq:AH0003_DDE}
     -A^{\odot\star}jH_{0003}  &= [J_1 K_{03} + M_{0003}]r^{\odot\star}.
 \end{align}
 With $e_1,e_2 \in \R^2$ the standard basis vectors and $\gamma_{1,03},\gamma_{2,03} \in \R$ we write
 \begin{align}\label{eq:K03_DDE}
     K_{03} = \gamma_{1,03}e_1 + \gamma_{2,03}e_2.
 \end{align}
 The solution to equation \cref{eq:AH0003_DDE} follows from \cref{cor:inverse_cases}.1 and is given by
 \begin{align}\label{eq:H0003_DDE}
    H_{0003}(\theta) = \gamma_{1,03}\Delta^{-1}(0)J_1e_1 + \gamma_{2,03}\Delta^{-1}(0)J_1e_2 + \Delta^{-1}(0)M_{0003}.
\end{align}
 The $z\beta_2^3$ terms yield the equation
 \begin{align}\label{eq:AH1003_DDE}
     (i\omega_0 I - A^{\odot\star})jH_{1003} &= [A_1(\varphi,K_{03})  + B(\varphi, H_{0003})]r^{\odot\star} - ib_{1,03}j\varphi + r_{1003},
 \end{align}
 where 
 \begin{align*}
     r_{1003} &= M_{1003}r^{\odot\star} - 3ib_{1,02}jH_{1001} - 3ib_{1,01}jH_{1002}.
 \end{align*}
Applying the Fredholm solvability condition to equation \cref{eq:AH1003_DDE}, after substituting equations \cref{eq:K03_DDE} and \cref{eq:H0003_DDE}, yields the following expression:
\begin{align}\label{eq:ib103_DDE}
    ib_{1,03} = \gamma_{1,03}p^T\Gamma_1(\varphi) + \gamma_{2,03}p^T\Gamma_2(\varphi) + p^T \tilde{M}_{1003},
\end{align}
where
 \begin{align*}
    \tilde{M}_{1003} &= B(\varphi,\Delta^{-1}(0)M_{0003}) + M_{1003}.
\end{align*}
To solve for $H_{1003}$ from equation \cref{eq:AH1003_DDE}, we apply \cref{cor:bordered_case}. The result is:
 \begin{align}\label{eq:H1003_DDE}
    H_{1003}(\theta) &=  \gamma_{1,03}B_{i\omega_0}^{Inv}(\Gamma_1(\varphi))(\theta)  + \gamma_{2,03}B_{i\omega_0}^{Inv}(\Gamma_2(\varphi))(\theta) \nonumber\\
    &- ib_{1,03}B_{i\omega_0}^{Inv}(0,q)(\theta)  + B_{1003}(\theta),
\end{align}
where 
\begin{align*}
    B_{1003}(\theta) &= B_{i\omega_0}^{Inv}(\tilde{M}_{1003})(\theta) - 3ib_{1,02}B_{i\omega_0}^{Inv}(0,H_{1001}(0),ib_{1,01}q)(\theta)\\
    &- 3ib_{1,01}B_{i\omega_0}^{Inv}(0,H_{1002}(0), [ib_{1,02}q + 2ib_{1,01}H_{1001}(0)], -b_{1,01}^2q)(\theta). 
\end{align*}
 Collecting the $z^2\beta_2^3$ and $z\bar{z}\beta_2^3$ terms yields the equations:
 \begin{align}
     (2i\omega_0I_n- A^{\odot\star})jH_{2003} &= [A_1(H_{2000},K_{03}) + 2B(\varphi, H_{1003}) + B(H_{0003}, H_{2000}) \label{eq:AH2003_DDE}\\
     &+ B_1(\varphi, \varphi, K_{03}) + C(\varphi, \varphi, H_{0003})]r^{\odot\star} -2ib_{1,03}jH_{2000} + r_{2003},\nonumber\\
     - A^{\odot\star}jH_{1103} &= [A_1(H_{1100},K_{03}) + 2\Re{}(B(\bar{\varphi}, H_{1003})) + B(H_{0003}, H_{1100})  \nonumber\\
     &+ B_1(\varphi, \bar{\varphi}, K_{03})+ C(\varphi, \bar{\varphi}, H_{0003})]r^{\odot\star} + r_{1103},\label{eq:AH1103_DDE}
 \end{align}
 where
 \begin{align*}
     r_{2003} &= M_{2003}r^{\odot\star} -(6ib_{1,02}jH_{2001} + 6ib_{1,01}jH_{2002}) \quad \textrm{ and } \quad r_{1103} = M_{1103} r^{\odot\star}.
 \end{align*}
 To solve equation \cref{eq:AH2003_DDE}, we apply \cref{cor:inverse_cases}.1 for the terms containing the multilinear forms and \cref{cor:inverse_case2} to the remaining terms. This yields the following solution:
 \begin{align}
     H_{2003}(\theta) &= e^{2i\omega_0\theta}\Delta^{-1}(2i\omega_0)[A_1(H_{2000},K_{03}) + 2B(\varphi, H_{1003}) + B(H_{0003}, H_{2000}) \label{eq:H2003_DDE_2}\\
     &+ B_1(\varphi, \varphi, K_{03}) + C(\varphi, \varphi, H_{0003}) + M_{2003}] \nonumber\\
     &-2ib_{1,03}\Delta^{-1}(2i\omega_0)[\Delta'(2i\omega_0) - \theta\Delta(2i\omega_0)]H_{2000}(\theta) \nonumber \\
     &-6ib_{1,02}e^{2i\omega_0\theta}\Delta^{-1}(2i\omega_0) \Bigl([\Delta'(2i
    \omega_0) - \theta\Delta(2i\omega_0)]H_{2001}(0) \nonumber\\
    &+ ib_{1,01} [\Delta''(2i\omega_0)- \theta^2\Delta(2i\omega_0)]H_{2000}(0)\Bigr) \nonumber \\
    &- 6ib_{1,01}e^{2i\omega_0\theta}\Delta^{-1}(2i\omega_0)\Bigl( [\Delta'(2i\omega_0) - \theta\Delta(2i\omega_0)]H_{2002}(0)\nonumber\\
    &+ [\Delta''(2i\omega_0) - \theta^2\Delta(2i\omega_0)](ib_{1,02}H_{2000}(0) + 2ib_{1,01}H_{2001}(0)) \nonumber \\
    & - \frac{4}{3}b_{1,01}^2[\Delta'''(2i\omega_0) - \theta^3\Delta(2i\omega_0)]H_{2000}(0)\Bigr), 
\end{align}
Since the right-hand side of equation \cref{eq:AH1103_DDE} is of the form $(w_0,0)$, the solution for $H_{1103}$ is obtained by simply applying \cref{cor:inverse_cases}.1. This results in the equation:
\begin{align}
     H_{1103}(\theta) &= \Delta^{-1}(0)[A_1(H_{1100},K_{03}) + 2\Re{}(B(\bar{\varphi}, H_{1003})) + B(H_{0003}, H_{1100})  \nonumber\\
     &+ B_1(\varphi, \bar{\varphi}, K_{03})+ C(\varphi, \bar{\varphi}, H_{0003}) + M_{1103}]. \label{eq:H1103_DDE_2}
 \end{align}
Substituting equations \cref{eq:K03_DDE},\cref{eq:H0003_DDE}, and \cref{eq:H1003_DDE} into equations \cref{eq:H2003_DDE_2},\cref{eq:H1103_DDE_2} result in:
\begin{align}
    H_{2003}(\theta) &= \gamma_{1,03}e^{2i\omega_0\theta}\Delta^{-1}(2i\omega_0)\Lambda_1(H_{2000}, \varphi, \varphi) + \gamma_{2,03}e^{2i\omega_0\theta}\Delta^{-1}(2i\omega_0)\Lambda_2(H_{2000}, \varphi, \varphi) \nonumber \\
    &-2ib_{1,03}\Delta^{-1}(2i\omega_0)\left([\Delta'(2i\omega_0) - \theta\Delta(2i\omega_0)]H_{2000}(\theta) +e^{2i\omega_0\theta}B\left(\varphi, B_{i\omega_0}^{INV}(0,q)\right)\right) \nonumber \\
    &+ e^{2i\omega_0\theta}\Delta^{-1}(2i\omega_0)\tilde{r}_{2003}, \label{eq:H2003_DDE_3}\\
    H_{1103}(\theta) &= \gamma_{1,
    03}\Delta^{-1}(0)\Pi_1(H_{1100},\bar{\varphi},\varphi) +  \gamma_{2,
    03}\Delta^{-1}(0)\Pi_1(H_{1100},\bar{\varphi},\varphi) \nonumber\\
    &-2b_{1,03}\Delta^{-1}(0)\Re{}\left\{ iB(\left(\bar{\varphi},B_{i\omega_0}^{INV}(0,q)\right)\right\} + \Delta^{-1}(0)\tilde{r}_{1103},\label{eq:H1103_DDE_3}
\end{align}
where we have 
\begin{align*}
    \tilde{r}_{2003} &= M_{2003} + 2B(\varphi, B_{1003})  + B(H_{2000}, \Delta^{-1}(0)M_{0003}) + C(\varphi,\bar{\varphi}, \Delta^{-1}(0)M_{0003}) \\
    &-6ib_{1,02} \Bigl([\Delta'(2i
    \omega_0) - \theta\Delta(2i\omega_0)]H_{2001}(0) \nonumber\\
    &+ ib_{1,01} [\Delta''(2i\omega_0)- \theta^2\Delta(2i\omega_0)]H_{2000}(0)\Bigr) \nonumber \\
    &- 6ib_{1,01}\Bigl( [\Delta'(2i\omega_0) - \theta\Delta(2i\omega_0)]H_{2002}(0)\nonumber\\
    &+ [\Delta''(2i\omega_0) - \theta^2\Delta(2i\omega_0)](ib_{1,02}H_{2000}(0) + 2ib_{1,01}H_{2001}(0)) \nonumber \\
    & - \frac{4}{3}b_{1,01}^2[\Delta'''(2i\omega_0) - \theta^3\Delta(2i\omega_0)]H_{2000}(0)\Bigr), \\
    \tilde{r}_{1103} &= M_{1103} + 2\Re{}\{ B(\bar{\varphi}, B_{1003}) \} + B(H_{1100}, \Delta^{-1}(0)M_{0003}) + C(\varphi,\bar{\varphi}, \Delta^{-1}(0)M_{0003}).
\end{align*}
 Applying the Fredholm solvability condition to the $z^2\bar{z}\beta_2^3$ terms yields the equation:
 \begin{align}
     ib_{2,03} &=\frac{1}{2} p^T [A_1(H_{2100},K_{03}) + 2B(\varphi, H_{1103}) + B(\bar{\varphi}, H_{2003})  + B(H_{0003}, H_{2100}) \nonumber \\
     &+ B(\overline{H}_{1003}, H_{2000}) + 2B(H_{1003}, H_{1100})  + 2B_1(\varphi, H_{1100}, K_{03}) \nonumber\\
     &+ B_1(\bar{\varphi}, H_{2000}, K_{03}) + C(\varphi,\varphi,\overline{H}_{1003}) + 2C(\varphi,\bar{\varphi},H_{1003}) + 2C(\varphi,H_{0003},H_{1100})\nonumber\\
     &+ C(\bar{\varphi},H_{0003},H_{2000})  + C_1(\varphi,\varphi,\bar{\varphi},K_{03}) + D(\varphi,\varphi,\bar{\varphi},H_{0003}) + M_{2103}].\label{eq:AH2103_DDE}
 \end{align}
Observe that equations  \cref{eq:H0003_DDE},\cref{eq:H1003_DDE},\cref{eq:H2003_DDE_3}, \cref{eq:H1103_DDE_3}, and \cref{eq:AH2103_DDE} correspond to equations \cref{eq:H0002_DDE},\cref{eq:H10mu2_DDE},\cref{eq:H20mu2_DDE}, \cref{eq:H11mu_DDE}, and \cref{eq:b2mu2_DDE} respectively, for $\mu = (03)$. Consequently, we can solve for the coefficients $\gamma_{1,03}$ and $\gamma_{2,03}$ by utilizing system \cref{eq:PQ2_DDE} for $\mu = (03)$.

\subsection{Remaining center manifold coefficients}
\subsubsection{Parameter-independent coefficients ODE}
For our center manifold approximation, we also need the following parameter-independent coefficients whose expressions can be found by respectively collecting the $w^5$,$w^6$,$w^5\bar{w}$, $w^7$, $w^6\bar{w}$,$w^5\bar{w}^2$ terms. The resulting equations are
\begin{align*}
    H_{5000} &= (5i\omega_0I_n- A)^{-1}[5B(q,H_{4000}) + 10B(H_{2000},H_{3000}) + 10C(q,q,H_{3000}) \\
    &+ 15C(q,H_{2000},H_{2000}) + 10D(q,q,q,H_{2000}) + E(q,q,q,q,q)], \\
    H_{6000} &= (6i\omega_0I_n - A)^{-1}[6B(q,H_{5000}) + 15B(H_{2000},H_{4000}) + 10B(H_{3000},H_{3000}) \\
    &+ 15C(q,q,H_{4000}) + 60C(q,H_{2000},H_{3000}) + 15C(H_{2000},H_{2000},H_{2000})\\
    &+ 20D(q,q,q,H_{3000}) + 45D(q,q,H_{2000}.H_{2000}) + 15E(q,q,q,q,H_{2000}) \\
    &+ K(q,q,q,q,q,q)],\\
    H_{5100} &= (4i\omega_0I_n- A)^{-1}[5B(q,H_{4100}) + B(\bar{q},H_{5000}) + 5B(H_{1100},H_{4000}) \\
    &+ 10B(H_{2000},H_{3100}) + 10B(H_{2100},H_{3000}) + 10C(q,q,H_{3100}) + 5C(q,\bar{q},H_{4000}) \\
    &+ 20C(q,H_{1100},H_{3000}) + 30C(q,H_{2000},H_{2100}) + 10C(\bar{q},H_{2000},H_{3000}) \\
    &+ 15C(H_{1100},H_{2000},H_{2000}) + 10D(q,q,q,H_{2100}) + 10D(q,q,\bar{q},H_{3000}) \\
    &+ 30D(q,q,H_{1100},H_{2000})+ 15D(q,\bar{q},H_{2000},H_{2000}) + 5E(q,q,q,q,H_{1100}) \\
    &+ 10E(q,q,q,\bar{q},H_{2000}) + K(q,q,q,q,q,\bar{q}) - 20c_1(0)H_{4000}], \\
    H_{7000} &= (7i\omega_0I_n- A)^{-1}[7B(q,H_{6000}) + 21B(H_{2000},H_{5000}) + 35B(H_{3000},H_{4000}) \\
    &+ 21C(q,q,H_{5000}) + 105C(q,H_{2000},H_{4000}) + 70C(q,H_{3000},H_{3000}) \\
    &+ 105C(H_{2000},H_{2000},H_{3000}) + 35D(q,q,q,H_{4000}) \\
    &+ 210D(q,q,H_{2000},H_{3000}) + 105D(q,H_{2000},H_{2000},H_{2000}) \\
    &+ 35E(q,q,q,q,H_{3000}) + 105E(q,q,q,H_{2000},H_{2000}) \\
    &+ 21K(q,q,q,q,q,H_{2000}) + L(q,q,q,q,q,q,q)], \\
    H_{6100} &= (5i\omega_0I_n- A)^{-1}[6B(q,H_{5100}) + B(\bar{q},H_{6000}) + 6B(H_{1100},H_{5000}) \\
    &+ 15B(H_{2000},H_{4100}) + 15B(H_{2100},H_{4000}) + 20B(H_{3000},H_{3100}) \\
    &+ 15C(q,q,H_{4100}) + 6C(q,\bar{q},H_{5000}) + 30C(q,H_{1100},H_{4000})\\
    &+ 60C(q,H_{2000},H_{3100}) + 60C(q,H_{2100},H_{3000}) + 15C(\bar{q},H_{2000},H_{4000}) \\
    &+ 10C(\bar{q},H_{3000},H_{3000}) + 60C(H_{1100},H_{2000},H_{3000})\\
    &+ 45C(H_{2000},H_{2000},H_{2100}) + 20D(q,q,q,H_{3100}) + 15D(q,q,\bar{q},H_{4000})\\
    &+ 60D(q,q,H_{1100},H_{3000}) + 90D(q,q,H_{2000},H_{2100}) \\
    &+ 60D(q,\bar{q},H_{2000},H_{3000}) + 90D(q,H_{1100},H_{2000},H_{2000}) \\
    &+ 15D(\bar{q},H_{2000},H_{2000},H_{2000}) + 15E(q,q,q,q,H_{2100}) \\
    &+ 20E(q,q,q,\bar{q},H_{3000}) + 60E(q,q,q,H_{1100},H_{2000}) \\
    &+ 45E(q,q,\bar{q},H_{2000},H_{2000}) + 6K(q,q,q,q,q,H_{1100}) \\
    &+ 15K(q,q,q,q,\bar{q},H_{2000}) + L(q,q,q,q,q,q,\bar{q}) - 30c_1(0)H_{5000}],\\
    H_{5200} &= (3i\omega_0I_n- A)^{-1}[5B(q,H_{4200}) + 2B(\bar{q},H_{5100}) \\
    &+ B(\overline{H}_{2000},H_{5000}) + 10B(H_{1100},H_{4100}) + 5B(\overline{H}_{2100},H_{4000}) \\
    &+ 10B(H_{2000},H_{3200}) + 20B(H_{2100},H_{3100}) + 10B(H_{2200},H_{3000}) \\
    &+ 10C(q,q,H_{3200}) + 10C(q,\bar{q},H_{4100}) + 5C(q, \overline{H}_{2000},H_{4000}) \\
    &+ 40C(q,H_{1100},H_{3100}) + 20C(q,\overline{H}_{2100},H_{3000}) \\
    &+ 30C(q,H_{2000},H_{2200}) + 30C(q,H_{2100},H_{2100}) + C(\bar{q},\bar{q},H_{5000}) \\
    &+ 10C(\bar{q},H_{1100},H_{4000}) + 20C(\bar{q},H_{2000},H_{3100}) + 20C(\bar{q},H_{2100},H_{3000}) \\
    &+ 10C(\overline{H}_{2000},H_{2000},H_{3000}) + 20C(H_{1100},H_{1100},H_{3000}) \\
    &+ 60C(H_{1100},H_{2000},H_{2100}) + 15C(\overline{H}_{2100},H_{2000},H_{2000}) \\
    &+ 10D(q,q,q,H_{2200}) + 20D(q,q,\bar{q},H_{3100}) + 10D(q,q,\overline{H}_{2000},H_{3000})\\
    &+ 60D(q,q,H_{1100},H_{2100}) + 30D(q,q,\overline{H}_{2100},H_{2000}) \\
    &+ 5D(q,\bar{q},\bar{q},H_{4000}) + 40D(q,\bar{q},H_{1100},H_{3000}) + 60D(q,\bar{q},H_{2000},H_{2100}) \\
    &+ 15D(q,\overline{H}_{2000},H_{2000},H_{2000}) + 60D(q,H_{1100},H_{1100},H_{2000}) \\
    &+ 10D(\bar{q},\bar{q},H_{2000},H_{3000}) + 30D(\bar{q},H_{1100},H_{2000},H_{2000}) \\
    &+ 5E(q,q,q,q,\overline{H}_{2100}) + 20E(q,q,q,\bar{q},H_{2100}) + 10E(q,q,q,\overline{H}_{2000},H_{2000})\\
    &+ 20E(q,q,q,H_{1100},H_{1100}) + 10E(q,q,\bar{q},\bar{q},H_{3000})\\
    &+ 60E(q,q,\bar{q},H_{1100},H_{2000}) + 15E(q,\bar{q},\bar{q},H_{2000},H_{2000}) \\
    &+ K(q,q,q,q,q,\overline{H}_{2000}) + 10K(q,q,q,q,\bar{q},H_{1100}) + 10K(q,q,q,\bar{q},\bar{q},H_{2000}) \\
    &+ L(q,q,q,q,q,\bar{q},\bar{q})-(120c_2(0)H_{3000} + (40c_1(0) + 10\overline{c_1}(0))H_{4100})].
\end{align*}

\subsubsection{Parameter-dependent coefficients ODE}

We also need the following coefficients
 \begin{align*}
     H_{3010} &= (3i\omega_0I_n -A)^{-1} [A_1(H_{3000},K_{10}) + 3B(q,H_{2010}) + B(H_{0010},H_{3000}) \nonumber\\
    &+ 3B(H_{1010},H_{2000}) + 3B_1(q,H_{2000},K_{10}) + 3C(q,q,H_{1010})+ 3C(q,H_{0010},H_{2000}) \\
    &+ C_1(q,q,q,K_{10}) + D(q,q,q,H_{0010}) - 3(1 + ib_{1,10})H_{3000}], \\
    H_{4001} &= (4i\omega_0 I_n - A)^{-1}[A_1(H_{4000},K_{01}) + 4B(q,H_{3001}) + B(H_{0001},H_{4000}) \\
    &+ 4B(H_{1001},H_{3000}) + 6B(H_{2000},H_{2001})+ 4B_1(q,H_{3000},K_{01}) \\
    &+ 3B_1(H_{2000},H_{2000},K_{01}) + 6C(q,q,H_{2001}) + 4C(q,H_{0001},H_{3000}) \\
    &+ 12C(q,H_{1001},H_{2000}) + 3C(H_{0001},H_{2000},H_{2000}) \\
    &+ 6C_1(q,q,H_{2000},K_{01}) + 4D(q,q,q,H_{1001}) + 6D(q,q,H_{0001},H_{2000}) \\
    &+ D_1(q,q,q,q,K_{01}) + E(q,q,q,q,H_{0001}) - 4ib_{1,01}H_{4000}], \\
    H_{5001} &= (5i\omega_0I_n - A)^{-1}[A_1(H_{5000},K_{01}) + 5B(q,H_{4001}) + B(H_{0001},H_{5000}) \\
    &+ 5B(H_{1001},H_{4000}) + 10B(H_{2000},H_{3001}) + 10B(H_{2001},H_{3000}) \\
    &+ 5B_1(q,H_{4000},K_{01}) + 10B_1(H_{2000},H_{3000},K_{01}) + 10C(q,q,H_{3001}) \\
    &+ 5C(q,H_{0001},H_{4000}) + 20C(q,H_{1001},H_{3000}) + 30C(q,H_{2000},H_{2001}) \\
    &+ 10C(H_{0001},H_{2000},H_{3000}) + 15C(H_{1001},H_{2000},H_{2000}) \\
    &+ 10C_1(q,q,H_{3000},K_{01}) + 15C_1(q,H_{2000},H_{2000},K_{01}) \\
    &+ 10D(q,q,q,H_{2001}) + 10D(q,q,H_{0001},H_{3000}) + 30D(q,q,H_{1001},H_{2000}) \\
    &+ 15D(q,H_{0001},H_{2000},H_{2000}) + 10D_1(q,q,q,H_{2000},K_{01}) \\
    &+ 5E(q,q,q,q,H_{1001}) + 10E(q,q,q,H_{0001},H_{2000})\\
    &+ E_1(q,q,q,q,q,K_{01}) + K(q,q,q,q,q,H_{0001}) - 5ib_{1,01}H_{5000}],\\
    H_{4101} &= (3i\omega_0I_n - A)^{-1}[A_1(H_{4100},K_{01}) + 4B(q,H_{3101}) + B(\bar{q},H_{4001})\\
    &+ B(H_{0001},H_{4100}) + B(\overline{H}_{1001},H_{4000}) + 4B(H_{1001},H_{3100}) \\
    &+ 4B(H_{1100},H_{3001}) + 4B(H_{1101},H_{3000}) + 6B(H_{2000},H_{2101})\\
    &+ 6B(H_{2001},H_{2100}) + 4B_1(q,H_{3100},K_{01}) + B_1(\bar{q},H_{4000},K_{01})\\
    &+ 4B_1(H_{1100},H_{3000},K_{01}) + 6B_1(H_{2000},H_{2100},K_{01}) + 6C(q,q,H_{2101}) \\
    &+ 4C(q,\bar{q},H_{3001}) + 4C(q,H_{0001},H_{3100}) + 4C(q,\overline{H}_{1001},H_{3000}) \\
    &+ 12C(q,H_{1001},H_{2100}) + 12C(q,H_{1100},H_{2001}) + 12C(q,H_{1101},H_{2000}) \\
    &+ C(\bar{q},H_{0001},H_{4000}) + 4C(\bar{q},H_{1001},H_{3000}) + 6C(\bar{q},H_{2000},H_{2001}) \\
    &+ 4C(H_{0001},H_{1100},H_{3000}) + 6C(H_{0001},H_{2000},H_{2100}) \\
    &+ 3C(\overline{H}_{1001},H_{2000},H_{2000}) + 12C(H_{1001},H_{1100},H_{2000}) \\
    &+ 6C_1(q,q,H_{2100},K_{01}) + 4C_1(q,\bar{q},H_{3000},K_{01}) \\
    &+ 12C_1(q,H_{1100},H_{2000},K_{01}) + 3C_1(\bar{q},H_{2000},H_{2000},K_{01})\\
    &+ 4D(q,q,q,H_{1101}) + 6D(q,q,\bar{q},H_{2001}) + 6D(q,q,H_{0001},H_{2100}) \\
    &+ 6D(q,q,\overline{H}_{1001},H_{2000}) + 12D(q,q,H_{1001},H_{1100}) \\
    &+ 4D(q,\bar{q},H_{0001},H_{3000}) + 12D(q,\bar{q},H_{1001},H_{2000}) \\
    &+ 12D(q,H_{0001},H_{1100},H_{2000}) + 3D(\bar{q},H_{0001},H_{2000},H_{2000}) \\
    &+ 4D_1(q,q,q,H_{1100},K_{01}) + 6D_1(q,q,\bar{q},H_{2000},K_{01})\\
    &+ E(q,q,q,q,\overline{H}_{1001}) + 4E(q,q,q,\bar{q},H_{1001}) + 4E(q,q,q,H_{0001},H_{1100}) \\
    &+ 6E(q,q,\bar{q},H_{0001},H_{2000}) + E_1(q,q,q,q,\bar{q},K_{01}) + K(q,q,q,q,\bar{q},H_{0001}) \\
    &-\left(12(1 + ib_{2,01})H_{3000} + 12c_1(0)H_{3001} + 3ib_{1,01}H_{4100}\right)].
 \end{align*}
Finally, we need the coefficient $H_{3002}$ which can be derived by collecting the $w^3\beta_2^2$ terms. This yields the equation
 \begin{align*}
    H_{3002} &= (3i\omega_0I_n- A)^{-1}[2A_1(H_{3001},K_{01}) + A_1(H_{3000},K_{02}) + 3B(q, H_{2002}) \\
    &+ 2B(H_{0001}, H_{3001}) + B(H_{0002}, H_{3000}) + 6B(H_{1001}, H_{2001}) \\
    &+ 3B(H_{1002}, H_{2000}) + A_2(H_{3000},K_{01}, K_{01}) + 6B_1(q, H_{2001}, K_{01}) \\
    &+ 3B_1(q, H_{2000},K_{02}) + 2B_1(H_{0001}, H_{3000},K_{01}) + 6B_1(H_{1001}, H_{2000},K_{01}) \\
    &+ 3C(q, q, H_{1002}) + 6C(q, H_{0001}, H_{2001}) + 3C(q, H_{0002}, H_{2000}) \\
    &+ 6C(q, H_{1001}, H_{1001}) + C(H_{0001}, H_{0001}, H_{3000}) \\
    &+ 6C(H_{0001}, H_{1001}, H_{2000}) + 3B_2(q, H_{2000}, K_{01}, K_{01}) + C_1(q, q, q, K_{02}) \\
    &+ 6C_1(q, q, H_{1001}, K_{01}) + 6C_1(q, H_{0001}, H_{2000}, K_{01}) + D(q, q, q, H_{0002}) \\
    &+ 6D(q, q, H_{0001}, H_{1001}) + 3D(q, H_{0001}, H_{0001}, H_{2000}) \\
    &+ C_2(q, q, q, K_{01}, K_{01}) + 2D_1(q, q, q, H_{0001}, K_{01}) \\
    &+ E(q, q, q, H_{0001}, H_{0001}) - (3ib_{1,02}H_{3000} + 6ib_{1,01}H_{3001})].
 \end{align*}

\subsubsection{Parameter-independent coefficients DDE}
For our center manifold approximation, we also need the following parameter-independent coefficients 
\begin{align*}
    H_{5000}(\theta) &= e^{5i\omega_0\theta}\Delta^{-1}(5i\omega_0)[5B(\varphi,H_{4000}) + 10B(H_{2000},H_{3000}) + 10C(\varphi,\varphi,H_{3000}) \\
    &+ 15C(\varphi,H_{2000},H_{2000}) + 10D(\varphi,\varphi,\varphi,H_{2000}) + E(\varphi,\varphi,\varphi,\varphi,\varphi)], \\
    H_{6000}(\theta) &= e^{6i\omega_0\theta}\Delta^{-1}(6i\omega_0)[6B(\varphi,H_{5000}) + 15B(H_{2000},H_{4000}) + 10B(H_{3000},H_{3000}) \\
    &+ 15C(\varphi,\varphi,H_{4000}) + 60C(\varphi,H_{2000},H_{3000}) + 15C(H_{2000},H_{2000},H_{2000})\\
    &+ 20D(\varphi,\varphi,\varphi,H_{3000}) + 45D(\varphi,\varphi,H_{2000}.H_{2000}) + 15E(\varphi,\varphi,\varphi,\varphi,H_{2000}) \\
    &+ K(\varphi,\varphi,\varphi,\varphi,\varphi,\varphi)],\\
    H_{5100}(\theta) &= e^{4i\omega_0\theta}\Delta^{-1}(4i\omega_0)[5B(\varphi,H_{4100}) + B(\bar{\varphi},H_{5000}) + 5B(H_{1100},H_{4000}) \\
    &+ 10B(H_{2000},H_{3100}) + 10B(H_{2100},H_{3000}) + 10C(\varphi,\varphi,H_{3100}) + 5C(\varphi,\bar{\varphi},H_{4000}) \\
    &+ 20C(\varphi,H_{1100},H_{3000}) + 30C(\varphi,H_{2000},H_{2100}) + 10C(\bar{\varphi},H_{2000},H_{3000}) \\
    &+ 15C(H_{1100},H_{2000},H_{2000}) + 10D(\varphi,\varphi,\varphi,H_{2100}) + 10D(\varphi,\varphi,\bar{\varphi},H_{3000}) \\
    &+ 30D(\varphi,\varphi,H_{1100},H_{2000})+ 15D(\varphi,\bar{\varphi},H_{2000},H_{2000}) + 5E(\varphi,\varphi,\varphi,\varphi,H_{1100}) \\
    &+ 10E(\varphi,\varphi,\varphi,\bar{\varphi},H_{2000}) + K(\varphi,\varphi,\varphi,\varphi,\varphi,\bar{\varphi})] \\
    &- 20c_1(0)\Delta^{-1}(4i\omega_0)[\Delta'(4i\omega_0) - \theta\Delta(4i\omega_0)]H_{4000}(\theta), \\
    H_{7000}(\theta) &= e^{7i\omega_0\theta}\Delta^{-1}(7i\omega_0)[7B(\varphi,H_{6000}) + 21B(H_{2000},H_{5000}) + 35B(H_{3000},H_{4000}) \\
    &+ 21C(\varphi,\varphi,H_{5000}) + 105C(\varphi,H_{2000},H_{4000}) + 70C(\varphi,H_{3000},H_{3000}) \\
    &+ 105C(H_{2000},H_{2000},H_{3000}) + 35D(\varphi,\varphi,\varphi,H_{4000}) \\
    &+ 210D(\varphi,\varphi,H_{2000},H_{3000}) + 105D(\varphi,H_{2000},H_{2000},H_{2000}) \\
    &+ 35E(\varphi,\varphi,\varphi,\varphi,H_{3000}) + 105E(\varphi,\varphi,\varphi,H_{2000},H_{2000}) \\
    &+ 21K(\varphi,\varphi,\varphi,\varphi,\varphi,H_{2000}) + L(\varphi,\varphi,\varphi,\varphi,\varphi,\varphi,\varphi)], \\
    H_{6100}(\theta) &= e^{5i\omega_0\theta}\Delta^{-1}(5i\omega_0)[6B(\varphi,H_{5100}) + B(\bar{\varphi},H_{6000}) + 6B(H_{1100},H_{5000}) \\
    &+ 15B(H_{2000},H_{4100}) + 15B(H_{2100},H_{4000}) + 20B(H_{3000},H_{3100}) \\
    &+ 15C(\varphi,\varphi,H_{4100}) + 6C(\varphi,\bar{\varphi},H_{5000}) + 30C(\varphi,H_{1100},H_{4000})\\
    &+ 60C(\varphi,H_{2000},H_{3100}) + 60C(\varphi,H_{2100},H_{3000}) + 15C(\bar{\varphi},H_{2000},H_{4000}) \\
    &+ 10C(\bar{\varphi},H_{3000},H_{3000}) + 60C(H_{1100},H_{2000},H_{3000})\\
    &+ 45C(H_{2000},H_{2000},H_{2100}) + 20D(\varphi,\varphi,\varphi,H_{3100}) + 15D(\varphi,\varphi,\bar{\varphi},H_{4000})\\
    &+ 60D(\varphi,\varphi,H_{1100},H_{3000}) + 90D(\varphi,\varphi,H_{2000},H_{2100}) \\
    &+ 60D(\varphi,\bar{\varphi},H_{2000},H_{3000}) + 90D(\varphi,H_{1100},H_{2000},H_{2000}) \\
    &+ 15D(\bar{\varphi},H_{2000},H_{2000},H_{2000}) + 15E(\varphi,\varphi,\varphi,\varphi,H_{2100}) \\
    &+ 20E(\varphi,\varphi,\varphi,\bar{\varphi},H_{3000}) + 60E(\varphi,\varphi,\varphi,H_{1100},H_{2000}) \\
    &+ 45E(\varphi,\varphi,\bar{\varphi},H_{2000},H_{2000}) + 6K(\varphi,\varphi,\varphi,\varphi,\varphi,H_{1100}) \\
    &+ 15K(\varphi,\varphi,\varphi,\varphi,\bar{\varphi},H_{2000}) + L(\varphi,\varphi,\varphi,\varphi,\varphi,\varphi,\bar{\varphi})] \\
    &- 30c_1(0)\Delta^{-1}(5i\omega_0)[\Delta'(5i\omega_0) - \theta\Delta(5i\omega_0)]H_{5000}(\theta), \\
    H_{5200}(\theta) &= e^{3i\omega_0\theta}\Delta^{-1}(3i\omega_0)[5B(\varphi,H_{4200}) + 2B(\bar{\varphi},H_{5100}) \\
    &+ B(\overline{H}_{2000},H_{5000}) + 10B(H_{1100},H_{4100}) + 5B(\overline{H}_{2100},H_{4000}) \\
    &+ 10B(H_{2000},H_{3200}) + 20B(H_{2100},H_{3100}) + 10B(H_{2200},H_{3000}) \\
    &+ 10C(\varphi,\varphi,H_{3200}) + 10C(\varphi,\bar{\varphi},H_{4100}) + 5C(\varphi, \overline{H}_{2000},H_{4000}) \\
    &+ 40C(\varphi,H_{1100},H_{3100}) + 20C(\varphi,\overline{H}_{2100},H_{3000}) \\
    &+ 30C(\varphi,H_{2000},H_{2200}) + 30C(\varphi,H_{2100},H_{2100}) + C(\bar{\varphi},\bar{\varphi},H_{5000}) \\
    &+ 10C(\bar{\varphi},H_{1100},H_{4000}) + 20C(\bar{\varphi},H_{2000},H_{3100}) + 20C(\bar{\varphi},H_{2100},H_{3000}) \\
    &+ 10C(\overline{H}_{2000},H_{2000},H_{3000}) + 20C(H_{1100},H_{1100},H_{3000}) \\
    &+ 60C(H_{1100},H_{2000},H_{2100}) + 15C(\overline{H}_{2100},H_{2000},H_{2000}) \\
    &+ 10D(\varphi,\varphi,\varphi,H_{2200}) + 20D(\varphi,\varphi,\bar{\varphi},H_{3100}) + 10D(\varphi,\varphi,\overline{H}_{2000},H_{3000})\\
    &+ 60D(\varphi,\varphi,H_{1100},H_{2100}) + 30D(\varphi,\varphi,\overline{H}_{2100},H_{2000}) \\
    &+ 5D(\varphi,\bar{\varphi},\bar{\varphi},H_{4000}) + 40D(\varphi,\bar{\varphi},H_{1100},H_{3000}) + 60D(\varphi,\bar{\varphi},H_{2000},H_{2100}) \\
    &+ 15D(\varphi,\overline{H}_{2000},H_{2000},H_{2000}) + 60D(\varphi,H_{1100},H_{1100},H_{2000}) \\
    &+ 10D(\bar{\varphi},\bar{\varphi},H_{2000},H_{3000}) + 30D(\bar{\varphi},H_{1100},H_{2000},H_{2000}) \\
    &+ 5E(\varphi,\varphi,\varphi,\varphi,\overline{H}_{2100}) + 20E(\varphi,\varphi,\varphi,\bar{\varphi},H_{2100}) + 10E(\varphi,\varphi,\varphi,\overline{H}_{2000},H_{2000})\\
    &+ 20E(\varphi,\varphi,\varphi,H_{1100},H_{1100}) + 10E(\varphi,\varphi,\bar{\varphi},\bar{\varphi},H_{3000})\\
    &+ 60E(\varphi,\varphi,\bar{\varphi},H_{1100},H_{2000}) + 15E(\varphi,\bar{\varphi},\bar{\varphi},H_{2000},H_{2000}) \\
    &+ K(\varphi,\varphi,\varphi,\varphi,\varphi,\overline{H}_{2000}) + 10K(\varphi,\varphi,\varphi,\varphi,\bar{\varphi},H_{1100}) + 10K(\varphi,\varphi,\varphi,\bar{\varphi},\bar{\varphi},H_{2000}) \\
    &+ L(\varphi,\varphi,\varphi,\varphi,\varphi,\bar{\varphi},\bar{\varphi})] \\
    &-120c_2(0)\Delta^{-1}(3i\omega_0)[\Delta'(3i\omega_0) - \theta\Delta(3i\omega_0)]H_{3000}(\theta) \\
    &-(40c_1(0) + 10\overline{c_1}(0))e^{3i\omega_0\theta}\Delta^{-1}(3i\omega_0)\Bigl( [\Delta'(3i\omega_0) - \theta\Delta(3i\omega_0)]H_{4100}(0) \nonumber\\
    &+ 6c_1(0)[\Delta''(3i\omega_0) - \theta^2\Delta(3i\omega_0)]H_{3000}(0)\Bigr),\\
    H_{4300}  &= B_{i\omega_0}^{Inv}\Bigl(M_{4300}, -144c_3(0)q\Bigr) -72(2c_2(0) + \overline{c_2}(0))B_{i\omega_0}^{Inv}\Bigl(0,H_{2100}(0) , 2c_1(0)q\Bigr)\\
    &-12i\Im{}\{c_1(0) \}B_{i\omega_0}^{Inv}\Bigl(0,H_{3200}(0), [12c_2(0)q \\
    &+ 6i\Im{}\left\{c_1(0)\right\}H_{2100}(0)],6i\Im{}\left\{c_1(0)\right\}c_1(0)q \Bigr).
\end{align*}

\subsubsection{Parameter-dependent coefficients DDE}
For our center manifold approximation, we also need the following parameter-dependent coefficients 
 \begin{align*}
     H_{3010}(\theta) &= e^{3i\omega_0\theta}\Delta^{-1}(3i\omega_0)[A_1(H_{3000},K_{10}) + 3B(\varphi,H_{2010}) + B(H_{0010},H_{3000}) \nonumber\\
    &+ 3B(H_{1010},H_{2000}) + 3B_1(\varphi,H_{2000},K_{10}) + 3C(\varphi,\varphi,H_{1010})+ 3C(\varphi,H_{0010},H_{2000}) \\
    &+ C_1(\varphi,\varphi,\varphi,K_{10}) + D(\varphi,\varphi,\varphi,H_{0010})] \\
    &- 3(1 + ib_{1,10})\Delta^{-1}(3i\omega_0)[\Delta'(3i\omega_0) - \theta\Delta(3i\omega_0)]H_{3000}(\theta)],\\
    H_{4001}(\theta) &= e^{4i\omega_0\theta}\Delta^{-1}(4i\omega_0)[A_1(H_{4000},K_{01}) + 4B(\varphi,H_{3001}) + B(H_{0001},H_{4000}) \\
    &+ 4B(H_{1001},H_{3000}) + 6B(H_{2000},H_{2001})+ 4B_1(\varphi,H_{3000},K_{01}) \\
    &+ 3B_1(H_{2000},H_{2000},K_{01}) + 6C(\varphi,\varphi,H_{2001}) + 4C(\varphi,H_{0001},H_{3000}) \\
    &+ 12C(\varphi,H_{1001},H_{2000}) + 3C(H_{0001},H_{2000},H_{2000}) \\
    &+ 6C_1(\varphi,\varphi,H_{2000},K_{01}) + 4D(\varphi,\varphi,\varphi,H_{1001}) + 6D(\varphi,\varphi,H_{0001},H_{2000}) \\
    &+ D_1(\varphi,\varphi,\varphi,\varphi,K_{01}) + E(\varphi,\varphi,\varphi,\varphi,H_{0001})]\\
    &- 4ib_{1,01}\Delta^{-1}(4i\omega_0)[\Delta'(4i\omega_0) - \theta\Delta(4i\omega_0)]H_{4000}(\theta), \\
    H_{5001}(\theta) &=  e^{5i\omega_0\theta}\Delta^{-1}(5i\omega_0)[A_1(H_{5000},K_{01}) + 5B(\varphi,H_{4001}) + B(H_{0001},H_{5000}) \\
    &+ 5B(H_{1001},H_{4000}) + 10B(H_{2000},H_{3001}) + 10B(H_{2001},H_{3000}) \\
    &+ 5B_1(\varphi,H_{4000},K_{01}) + 10B_1(H_{2000},H_{3000},K_{01}) + 10C(\varphi,\varphi,H_{3001}) \\
    &+ 5C(\varphi,H_{0001},H_{4000}) + 20C(\varphi,H_{1001},H_{3000}) + 30C(\varphi,H_{2000},H_{2001}) \\
    &+ 10C(H_{0001},H_{2000},H_{3000}) + 15C(H_{1001},H_{2000},H_{2000}) \\
    &+ 10C_1(\varphi,\varphi,H_{3000},K_{01}) + 15C_1(\varphi,H_{2000},H_{2000},K_{01}) \\
    &+ 10D(\varphi,\varphi,\varphi,H_{2001}) + 10D(\varphi,\varphi,H_{0001},H_{3000}) + 30D(\varphi,\varphi,H_{1001},H_{2000}) \\
    &+ 15D(\varphi,H_{0001},H_{2000},H_{2000}) + 10D_1(\varphi,\varphi,\varphi,H_{2000},K_{01}) \\
    &+ 5E(\varphi,\varphi,\varphi,\varphi,H_{1001}) + 10E(\varphi,\varphi,\varphi,H_{0001},H_{2000})\\
    &+ E_1(\varphi,\varphi,\varphi,\varphi,\varphi,K_{01}) + K(\varphi,\varphi,\varphi,\varphi,\varphi,H_{0001})] \\
    &- 5ib_{1,01}\Delta^{-1}(5i\omega_0)[\Delta'(5i\omega_0) - \theta\Delta(5i\omega_0)]H_{5000}(\theta),\\
    H_{4101}(\theta) &= e^{3i\omega_0\theta}\Delta^{-1}(3i\omega_0)[A_1(H_{4100},K_{01}) + 4B(\varphi,H_{3101}) + B(\bar{\varphi},H_{4001})\\
    &+ B(H_{0001},H_{4100}) + B(\overline{H}_{1001},H_{4000}) + 4B(H_{1001},H_{3100}) \\
    &+ 4B(H_{1100},H_{3001}) + 4B(H_{1101},H_{3000}) + 6B(H_{2000},H_{2101})\\
    &+ 6B(H_{2001},H_{2100}) + 4B_1(\varphi,H_{3100},K_{01}) + B_1(\bar{\varphi},H_{4000},K_{01})\\
    &+ 4B_1(H_{1100},H_{3000},K_{01}) + 6B_1(H_{2000},H_{2100},K_{01}) + 6C(\varphi,\varphi,H_{2101}) \\
    &+ 4C(\varphi,\bar{\varphi},H_{3001}) + 4C(\varphi,H_{0001},H_{3100}) + 4C(\varphi,\overline{H}_{1001},H_{3000}) \\
    &+ 12C(\varphi,H_{1001},H_{2100}) + 12C(\varphi,H_{1100},H_{2001}) + 12C(\varphi,H_{1101},H_{2000}) \\
    &+ C(\bar{\varphi},H_{0001},H_{4000}) + 4C(\bar{\varphi},H_{1001},H_{3000}) + 6C(\bar{\varphi},H_{2000},H_{2001}) \\
    &+ 4C(H_{0001},H_{1100},H_{3000}) + 6C(H_{0001},H_{2000},H_{2100}) \\
    &+ 3C(\overline{H}_{1001},H_{2000},H_{2000}) + 12C(H_{1001},H_{1100},H_{2000}) \\
    &+ 6C_1(\varphi,\varphi,H_{2100},K_{01}) + 4C_1(\varphi,\bar{\varphi},H_{3000},K_{01}) \\
    &+ 12C_1(\varphi,H_{1100},H_{2000},K_{01}) + 3C_1(\bar{\varphi},H_{2000},H_{2000},K_{01})\\
    &+ 4D(\varphi,\varphi,\varphi,H_{1101}) + 6D(\varphi,\varphi,\bar{\varphi},H_{2001}) + 6D(\varphi,\varphi,H_{0001},H_{2100}) \\
    &+ 6D(\varphi,\varphi,\overline{H}_{1001},H_{2000}) + 12D(\varphi,\varphi,H_{1001},H_{1100}) \\
    &+ 4D(\varphi,\bar{\varphi},H_{0001},H_{3000}) + 12D(\varphi,\bar{\varphi},H_{1001},H_{2000}) \\
    &+ 12D(\varphi,H_{0001},H_{1100},H_{2000}) + 3D(\bar{\varphi},H_{0001},H_{2000},H_{2000}) \\
    &+ 4D_1(\varphi,\varphi,\varphi,H_{1100},K_{01}) + 6D_1(\varphi,\varphi,\bar{\varphi},H_{2000},K_{01})\\
    &+ E(\varphi,\varphi,\varphi,\varphi,\overline{H}_{1001}) + 4E(\varphi,\varphi,\varphi,\bar{\varphi},H_{1001}) + 4E(\varphi,\varphi,\varphi,H_{0001},H_{1100}) \\
    &+ 6E(\varphi,\varphi,\bar{\varphi},H_{0001},H_{2000}) + E_1(\varphi,\varphi,\varphi,\varphi,\bar{\varphi},K_{01}) + K(\varphi,\varphi,\varphi,\varphi,\bar{\varphi},H_{0001})]\\
    &-12(1 + ib_{2,01})\Delta^{-1}(3i\omega_0)[\Delta'(3i\omega_0) - \theta\Delta(3i\omega_0)]H_{3000}(\theta) \\
    &- 12c_1(0)e^{3i\omega_0\theta}\Delta^{-1}(3i\omega_0)\Bigl( [\Delta'(3i\omega_0) - \theta\Delta(3i\omega_0)]H_{3001}(0) \nonumber\\
    &+ \frac{3}{2}ib_{1,01}[\Delta''(3i\omega_0) - \theta^2\Delta(3i\omega_0)]H_{3000}(0)\Bigr) \\
    &- 3ib_{1,01}e^{3i\omega_0\theta}\Delta^{-1}(3i\omega_0)\Bigl( [\Delta'(3i\omega_0) - \theta\Delta(3i\omega_0)]H_{4100}(0) \nonumber\\
    &+ 6c_1(0)[\Delta''(3i\omega_0) - \theta^2\Delta(3i\omega_0)]H_{3000}(0)\Bigr),\\
    H_{3201}&= B_{i\omega_0}^{Inv}\Bigl(M_{3201}, -12g_{3201}q\Bigr) - 12c_2(0)B_{i\omega_0}^{Inv}\Bigl(0,H_{1001}(0),ib_{1,01}q\Bigr)\\
     &-(18 + 6ib_{2,01})B_{i\omega_0}^{Inv}\Bigl(0,H_{2100}(0), 2c_1(0)q\Bigr) \\
     &-6i\Im{}\{c_1(0)\}B_{i\omega_0}^{Inv}\Bigl(0,H_{2101}(0), [2(1 + ib_{2,01})q + ib_{1,01}H_{2100}(0) \\
     &+ 2c_1(0)H_{1001}(0)], 2ib_{1,01}c_1(0)q\Bigr) \\
     &-ib_{1,01}B_{i\omega_0}^{Inv}\Bigl(0,H_{3200}(0), [12c_2(0)q + 6i\Im{}\left\{c_1(0)\right\}H_{2100}(0)],6i\Im{}\left\{c_1(0)\right\}c_1(0)q \Bigr),\\
    H_{3002}(\theta) &= e^{3i\omega_0\theta}\Delta^{-1}(3i\omega_0)[2A_1(H_{3001},K_{01}) + A_1(H_{3000},K_{02}) + 3B(\varphi, H_{2002}) \\
    &+ 2B(H_{0001}, H_{3001}) + B(H_{0002}, H_{3000}) + 6B(H_{1001}, H_{2001}) \\
    &+ 3B(H_{1002}, H_{2000}) + A_2(H_{3000},K_{01}, K_{01}) + 6B_1(\varphi, H_{2001}, K_{01}) \\
    &+ 3B_1(\varphi, H_{2000},K_{02}) + 2B_1(H_{0001}, H_{3000},K_{01}) + 6B_1(H_{1001}, H_{2000},K_{01}) \\
    &+ 3C(\varphi, \varphi, H_{1002}) + 6C(\varphi, H_{0001}, H_{2001}) + 3C(\varphi, H_{0002}, H_{2000}) \\
    &+ 6C(\varphi, H_{1001}, H_{1001}) + C(H_{0001}, H_{0001}, H_{3000}) \\
    &+ 6C(H_{0001}, H_{1001}, H_{2000}) + 3B_2(\varphi, H_{2000}, K_{01}, K_{01}) + C_1(\varphi, \varphi, \varphi, K_{02}) \\
    &+ 6C_1(\varphi, \varphi, H_{1001}, K_{01}) + 6C_1(\varphi, H_{0001}, H_{2000}, K_{01}) + D(\varphi, \varphi, \varphi, H_{0002}) \\
    &+ 6D(\varphi, \varphi, H_{0001}, H_{1001}) + 3D(\varphi, H_{0001}, H_{0001}, H_{2000}) \\
    &+ C_2(\varphi, \varphi, \varphi, K_{01}, K_{01}) + 2D_1(\varphi, \varphi, \varphi, H_{0001}, K_{01}) \\
    &+ E(\varphi, \varphi, \varphi, H_{0001}, H_{0001}) ]\\
    &- 3ib_{1,02}\Delta^{-1}(3i\omega_0)[\Delta'(3i\omega_0) - \theta\Delta(3i\omega_0)]H_{3000}(\theta)  \\
    &- 6ib_{1,01}e^{3i\omega_0\theta}\Delta^{-1}(3i\omega_0)\Bigl( [\Delta'(3i\omega_0) - \theta\Delta(3i\omega_0)]H_{3001}(0) \nonumber\\
    &+ \frac{3}{2}ib_{1,01}[\Delta''(3i\omega_0) - \theta^2\Delta(3i\omega_0)]H_{3000}(0)\Bigr).
 \end{align*}
\newpage
\section{Terms collected from the homological equation for ODEs}\label{sec:equations_hom}
In this supplement, all the equations collected from the homological equation for the generalized Hopf bifurcation in ODEs \cref{eq:HOM_GH_ODE} are presented. The collected equations are obtained with the help of Mathematica. We use the truncated normal form expressed in terms of $\beta$
\[
\dot{w} =  \lambda(\beta)w + c_1(\beta)w|w|^2 + c_2(\beta)w|w|^4 + c_3(0)w|w|^6, \quad w \in \C,\ \beta \in \R^2,
\]
where we expand
\begin{align*}
    \lambda(\beta) &= i\omega_0 + g_{1010}\beta_1 + g_{1001}\beta_2 + \frac{1}{2}g_{1020}\beta_1^2 + g_{1011}\beta_1\beta_2 +\frac{1}{2}g_{1002}\beta_2^2 + \frac{1}{6}g_{1003}\beta_2^3, \\
    c_1(\beta) &= c_1(0) + g_{2110}\beta_1 + g_{2101}\beta_2 + \frac{1}{2}g_{2120}\beta_1^2 + g_{2111}\beta_1\beta_2 +\frac{1}{2}g_{2102}\beta_2^2  + \frac{1}{6}g_{2103}\beta_2^3, \\
    c_2(\beta) &= c_2(0) + g_{3210}\beta_1 + g_{3201}\beta_2.
\end{align*}
Compared to \cref{eq:wdot} we have that for $\mu = (10),(01)$
\begin{align*}
g_{10\mu} &= \delta_{\mu}^{10} + b_{1,\mu}i,\\
g_{21\mu} &= \delta_{\mu}^{01} + b_{2,\mu}i,
\end{align*}
For $\mu =(20),(02),(11)$ we have $g_{10\mu} = b_{1,\mu}i$ and $g_{21\mu} = b_{2,\mu}i$. The vectorfield $F$ is expanded as \cref{eq:F_Taylor}. Note that $H_{ij\mu} = \overline{H}_{ji\mu}$. 
\subsection{Linear terms}
Collecting linear terms of \cref{eq:HOM_GH_ODE}:
\begin{align*}
    w&: \quad Aq = iw_0q \\
    \bar{w}&: \quad A\bar{q} = -iw_0 \bar{q} \\
    \beta_1&: AH_{0010} = -J_1K_{10} \\
    \beta_2&: AH_{0001} = -J_1K_{01}
\end{align*}

\subsection{Quadratic terms}
Collecting quadratic terms of \cref{eq:HOM_GH_ODE}:
\begin{align*}
    w^2&:  &(2i\omega_0 I_n - A)H_{2000} &= B(q,q) \\
    w\bar{w}&:  &-AH_{1100} &= B(q,\bar{q})\\
    w\beta_1&: &(i\omega_0I_n - A)H_{1010} &=  A_1(q,K_{10}) + B(q,H_{0010}) - g_{1010}q \\
    w\beta_2&: &(i\omega_0I_n - A)H_{1001} &=  A_1(q,K_{01}) + B(q,H_{0001}) - g_{1001}q \\
    \beta_1^2 &: &-AH_{0020} &= J_1K_{20} + 2A_1(H_{0010},K_{10}) + B(H_{0010},H_{0010})+ J_2(K_{10},K_{10})\\
    \beta_2^2 &:  &-AH_{0002} &= J_1K_{02} + 2A_1(H_{0001},K_{01}) + B(H_{0001},H_{0001})+ J_2(K_{01},K_{01}) \\
    \beta_1\beta_2 &: &-AH_{0011} &=  J_1K_{11} + A_1(H_{0010},K_{01}) + A_1(H_{0001},K_{10}) + B(H_{0001},H_{0010}) \\
    &&&+ J_2(K_{01},K_{10})
\end{align*}

\subsection{Qubic terms}
Collecting qubic terms of \cref{eq:HOM_GH_ODE}:
\begin{align*}
    w^3&: &(3i\omega_0I_n - A)H_{3000} &= 3B(q,H_{2000}) + C(q,q,q) \\
    w^2\bar{w} &: &(i\omega_0I_n - A)H_{2100} &=  2B(q,H_{1100}) +B(\bar{q},H_{2000}) + C(q,q,\bar{q}) -2c_1(0)q\\
    \beta_1w^2 &: &(2i\omega_0I_n - A)H_{2010} &=   A_1(H_{2000},K_{10}) +2B(q,H_{1010}) + B(H_{0010},H_{2000}) \\
    &&&+ B_1(q,q,K_{10}) + C(q,q,H_{0010}) - 2g_{1010}H_{2000} \\
    \beta_2w^2 &: &(2i\omega_0I_n - A)H_{2001} &=  A_1(H_{2000},K_{01}) + 2B(q,H_{1001}) + B(H_{0001},H_{2000}) \\
    &&&+ B_1(q,q,K_{01}) + C(q,q,H_{0001}) - 2g_{1001}H_{2000} \\
    w\bar{w}\beta_1 &: &-AH_{1110} &=  A_1(H_{1100},K_{10}) + B(q,H_{0110}) + B(\bar{q},H_{1010}) \\
    &&&+ B(H_{0010},H_{1100}) + B_1(q,\bar{q},K_{10}) + C(q,\bar{q},H_{0010})\\
    &&&-(g_{1010} + \bar{g}_{1010})H_{1100} \\
    w\bar{w}\beta_2 &: &-AH_{1101} &= A_1(H_{1100},K_{01}) + B(q,H_{0101}) + B(\bar{q},H_{1001}) \\
    &&&+ B(H_{0001},H_{1100}) + B_1(q,\bar{q},K_{01}) + C(q,\bar{q},H_{0001}) \\
    &&&-(g_{1001} + \bar{g}_{1001})H_{1100}\\
    w\beta_1^2 &: &(i\omega_0I_n - A)H_{1020} &=  A_1(q,K_{20}) + 2A_1(H_{1010},K_{10}) + B(q,H_{0020}) \\
    &&&+ 2B(H_{0010},H_{1010}) + A_2(q,K_{10},K_{10}) + 2B_1(q,H_{0010},K_{10}) \\
    &&&+ C(q,H_{0010},H_{0010}) - [g_{1020}q + 2g_{1010}H_{1010}]\\
    w\beta_2^2 &: &(i\omega_0I_n - A)H_{1002} &=  A_1(q,K_{02}) + 2A_1(H_{1001},K_{01}) + B(q,H_{0002}) \\
    &&&+ 2B(H_{0001},H_{1001}) + A_2(q,K_{01},K_{01}) + 2B_1(q,H_{0001},K_{01}) \\
    &&&+ C(q,H_{0001},H_{0001}) - [g_{1002}q + 2g_{1001}H_{1001}]\\
    w\beta_1\beta_2 &: &(i\omega_0I_n - A)H_{1011} &=  A_1(q,K_{11}) + A_1(H_{1010},K_{01}) \\
    &&&+ A_1(H_{1001},K_{10}) + B(q,H_{0011}) + B(H_{0001},H_{1010}) \\
    &&&+ B(H_{0010},H_{1001}) + A_2(q,K_{01},K_{10}) + B_1(q,H_{0010},K_{01}) \\
    &&&+ B_1(q,H_{0001},K_{10}) + C(q,H_{0001},H_{0010}) \\
    &&&-[g_{1011}q + g_{1010}H_{1001} + g_{1001}H_{1010}]\\
    \beta_2^3 &: &-AH_{0003}  &= J_1 K_{03}+ 3A_1(H_{0002},K_{01}) + 3A_1(H_{0001},K_{02}) \\
    &&&+ 3B(H_{0001}, H_{0002}) + 3J_2(K_{01}, K_{02}) + 3A_2(H_{0001},K_{01}, K_{01}) \\
    &&&+ 3B_1(H_{0001}, H_{0001},K_{01}) + J_3(K_{01}, K_{01}, K_{01}) \\
    &&&+ C(H_{0001}, H_{0001}, H_{0001})
\end{align*}

\subsection{Quartic terms}
Collecting quartic terms of \cref{eq:HOM_GH_ODE}:
\begin{align*}
    w^4 &: &(4i\omega_0I_n - A)H_{4000} &= 4B(q,H_{3000}) + 3B(H_{2000},H_{2000}) + 6C(q,q,,H_{2000}) \\
    &&&+ D(q,q,q,q)\\
    w^3\bar{w} &: &(2i\omega_0I_n - A)H_{3100} &=  3B(q,H_{2100}) + B(\bar{q},H_{3000}) + 3B(H_{1100},H_{2000}) \\
    &&&+ 3C(q,q,H_{1100}) +  3C(q,\bar{q},H_{2000}) + D(q,q,q,\bar{q}) \\
    &&&-6c_1(0)H_{2000} \\
    w^2\bar{w}^2 &:  &-AH_{2200} &=   2B(q,H_{1200}) + 2B(\bar{q},H_{2100}) + B(H_{0200},H_{2000}) \\
    &&&+ 2B(H_{1100},H_{1100}) + C(q,q,H_{0200}) + 4C(q,\bar{q},H_{1100}) \\
    &&&+ C(\bar{q},\bar{q},H_{2000}) + D(q,q,\bar{q},\bar{q}) -4(c_1(0) +\bar{c}_1(0))H_{1100}\\
    w^3\beta_1 &: &(3i\omega_0I_n - A)H_{3010} &= A_1(H_{3000},K_{10}) + 3B(q,H_{2010}) + B(H_{0010},H_{3000}) \\
    &&&+ 3B(H_{1010},H_{2000}) + 3B_1(q,H_{2000},K_{10}) + 3C(q,q,H_{1010}) \\
    &&&+ 3C(q,H_{0010},H_{2000}) + C_1(q,q,q,K_{10}) + D(q,q,q,H_{0010}) \\
    &&&- 3g_{1010}H_{3000}  \\
    w^3\beta_2 &: &(3i\omega_0I_n - A)H_{3001} &= A_1(H_{3000},K_{01}) + 3B(q,H_{2001}) + B(H_{0001},H_{3000}) \\
    &&&+ 3B(H_{1001},H_{2000}) + 3B_1(q,H_{2000},K_{01}) + 3C(q,q,H_{1001}) \\
    &&&+ 3C(q,H_{0001},H_{2000}) + C_1(q,q,q,K_{01}) + D(q,q,q,H_{0001}) \\
    &&& -3g_{1001}H_{3000} \\
    w^2\bar{w}\beta_1 &: &(i\omega_0I_n - A)H_{2110} &=  A_1(H_{2100},K_{10}) + 2B(q,H_{1110}) + B(\bar{q},H_{2010}) \\
    &&&+ B(H_{0010},H_{2100}) + B(H_{0110},H_{2000}) + 2B(H_{1010},H_{1100}) \\
    &&&+ 2B_1(q,H_{1100},K_{10}) + B_1(\bar{q},H_{2000},K_{10}) + C(q,q,H_{0110}) \\
    &&&+ 2C(q,\bar{q}.H_{1010}) + 2C(q,H_{0010},H_{1100}) + C(\bar{q},H_{0010},H_{2000}) \\
    &&&+ C_1(q,q,\bar{q},K_{10}) + D(q,q,\bar{q},H_{0010})\\
    &&&-[2g_{2110}q + (2g_{1010} + \bar{g}_{1010})H_{2100} + 2c_1(0)H_{1010}]\\
    w^2\bar{w}\beta_2 &: &(i\omega_0I_n - A)H_{2101} &=  A_1(H_{2100},K_{01}) + 2B(q,H_{1101}) + B(\bar{q},H_{2001}) \\
    &&&+ B(H_{0001},H_{2100}) + B(H_{0101},H_{2000}) + 2B(H_{1001},H_{1100}) \\
    &&&+ 2B_1(q,H_{1100},K_{01}) + B_1(\bar{q},H_{2000},K_{01}) + C(q,q,H_{0101}) \\
    &&&+ 2C(q,\bar{q}.H_{1001}) + 2C(q,H_{0001},H_{1100}) + C(\bar{q},H_{0001},H_{2000}) \\
    &&&+ C_1(q,q,\bar{q},K_{01}) + D(q,q,\bar{q},H_{0001}) \\
    &&&-[2g_{2101}q + (2g_{1001} + \bar{g}_{1001})H_{2100} + 2c_1(0)H_{1001}] \\
    w^2\beta_1^2 &: &(2i\omega_0I_n - A)H_{2020} &=  2A_1(H_{2010},K_{10}) + A_1(H_{2000},K_{20}) + 2B(q,H_{1020}) \\
    &&&+ 2B(H_{0010},H_{2010}) + B(H_{0020},H_{2000}) + 2B(H_{1010},H_{1010}) \\
    &&&+ A_2(H_{2000},K_{10},K_{10}) + B_1(q,q,K_{20}) + 4B_1(q,H_{1010},K_{10}) \\
    &&&+ 2B_1(H_{0010},H_{2000},K_{10}) + C(q,q,H_{0020}) + 4C(q,H_{0010},H_{1010}) \\
    &&&+ C(H_{0010},H_{0010},H_{2000}) + B_2(q,q,K_{10},K_{10}) \\
    &&&+ 2C_1(q,q,H_{0010},K_{10}) + D(q,q,H_{0010},H_{0010}) \\
    &&&- [2g_{1020}H_{2000} + 4g_{1010}H_{2010}] \\
    w^2\beta_2^2 &: &(2i\omega_0I_n - A)H_{2002} &=  2A_1(H_{2001},K_{01}) + A_1(H_{2000},K_{02}) + 2B(q,H_{1002}) \\
    &&&+ 2B(H_{0001},H_{2001}) + B(H_{0002},H_{2000}) + 2B(H_{1001},H_{1001}) \\
    &&&+ A_2(H_{2000},K_{01},K_{01}) + B_1(q,q,K_{02}) + 4B_1(q,H_{1001},K_{01}) \\
    &&&+ 2B_1(H_{0001},H_{2000},K_{01}) + C(q,q,H_{0002}) + 4C(q,H_{0001},H_{1001}) \\
    &&&+ C(H_{0001},H_{0001},H_{2000}) + B_2(q,q,K_{01},K_{01}) \\
    &&&+ 2C_1(q,q,H_{0001},K_{01}) + D(q,q,H_{0001},H_{0001}) \\
    &&&- [2g_{1002}H_{2000} + 4g_{1001}H_{2001}]
\end{align*}

\begin{align*}
    w^2\beta_1\beta_2 &: &(2i\omega_0I_n - A)H_{2011} &= A_1(H_{2010},K_{01}) + A_1(H_{2001},K_{10}) + A_1(H_{2000},K_{11}) \\
    &&&+ 2B(q,H_{1011}) + B(H_{0001},H_{2010}) + B(H_{0010},H_{2001}) \\
    &&&+ B(H_{0011},H_{2000}) + 2B(H_{1001},H_{1010}) + A_2(H_{2000},K_{01},K_{10}) \\
    &&&+ B_1(q,q,K_{11}) + 2B_1(q,H_{1010},K_{01}) + 2B_1(q,H_{1001},K_{10}) \\
    &&&+ B_1(H_{0010},H_{2000},K_{01}) + B_1(H_{0001},H_{2000},K_{10}) \\
    &&&+ C(q,q,H_{0011}) + 2C(q,H_{0001},H_{1010}) + 2C(q,H_{0010},H_{1001}) \\
    &&&+ C(H_{0001},H_{0010},H_{2000}) + B_2(q,q,K_{01},K_{10}) \\
    &&&+ C_1(q,q,H_{0010},K_{01}) + C_1(q,q,H_{0001},K_{10}) \\
    &&&+ D(q,q,H_{0001},H_{0010})\\
    &&&-[2g_{1011}H_{2000} + 2g_{1010}H_{2001} + 2g_{1001}H_{2010} ]\\
    w\bar{w}\beta_1^2 &: &-AH_{1120} &=  2A_1(H_{1110},K_{10}) + A_1(H_{1100},K_{20}) + B(q,H_{0120}) \\
    &&&+ B(\bar{q},H_{1020}) + 2B(H_{0010},H_{1110}) + B(H_{0020},H_{1100}) \\
    &&&+ 2B(H_{0110},H_{1010}) + A_2(H_{1100},K_{10},K_{10}) + B_1(q,\bar{q},K_{20}) \\
    &&&+ 2B_1(q,H_{0110},K_{10}) + 2B_1(\bar{q},H_{1010},K_{10}) \\
    &&&+ 2B_1(H_{0010},H_{1100},K_{10})+ C(q,\bar{q},H_{0020}) \\
    &&&+ 2C(q,H_{0010},H_{0110}) + 2C(\bar{q},H_{0010},H_{1010}) \\
    &&&+ C(H_{0010},H_{0010},H_{1100}) + B_2(q,\bar{q},K_{10},K_{10}) \\
    &&&+ 2C_1(q,\bar{q},H_{0010},K_{10}) + D(q,\bar{q},H_{0010},H_{0010}) \\
    &&&-[(g_{1020}+ \bar{g}_{1020})H_{1100} + 2(g_{1010} + \bar{g}_{1010})H_{1110}]\\
    w\bar{w}\beta_2^2 &: &-AH_{1102} &=  2A_1(H_{1101},K_{01}) + A_1(H_{1100},K_{02}) + B(q,H_{0102}) \\
    &&&+ B(\bar{q},H_{1002}) + 2B(H_{0001},H_{1101}) + B(H_{0002},H_{1100}) \\
    &&&+ 2B(H_{0101},H_{1001}) + A_2(H_{1100},K_{01},K_{01}) + B_1(q,\bar{q},K_{02}) \\
    &&&+ 2B_1(q,H_{0101},K_{01}) + 2B_1(\bar{q},H_{1001},K_{01}) \\
    &&&+ 2B_1(H_{0001},H_{1100},K_{01})+ C(q,\bar{q},H_{0002}) \\
    &&&+ 2C(q,H_{0001},H_{0101}) + 2C(\bar{q},H_{0001},H_{1001}) \\
    &&&+ C(H_{0001},H_{0001},H_{1100}) + B_2(q,\bar{q},K_{01},K_{01}) \\
    &&&+ 2C_1(q,\bar{q},H_{0001},K_{01}) + D(q,\bar{q},H_{0001},H_{0001})\\
    &&&-[(g_{1002}+ \bar{g}_{1002})H_{1100} + 2(g_{1001} + \bar{g}_{1001})H_{1101}]\\
    w\bar{w}\beta_1\beta_2 &: &-AH_{1111} &= A_1(H_{1110},K_{01}) + A_1(H_{1101},K_{10}) + A_1(H_{1100},K_{11}) \\
    &&&+ B(q,H_{0111}) + B(\bar{q},H_{1011}) + B(H_{0001},H_{1110}) \\
    &&&+ B(H_{0010},H_{1101}) + B(H_{0011},H_{1100}) + B(H_{0101},H_{1010}) \\
    &&&+ B(H_{0110},H_{1001}) + A_2(H_{1100},K_{01},K_{10}) \\
    &&&+ B_1(q,\bar{q},K_{11}) + B_1(q,H_{0110},K_{01}) \\
    &&&+ B_1(q,H_{0101},K_{10}) + B_1(\bar{q},H_{1010},K_{01}) \\
    &&&+ B_1(\bar{q},H_{1001},K_{10}) + B_1(H_{0010},H_{1100},K_{01}) \\
    &&&+ B_1(H_{0001},H_{1100},K_{10}) + C(q,\bar{q},H_{0011}) \\
    &&&+ C(q,H_{0001},H_{0110}) + C(q,H_{0010},H_{0101}) \\
    &&&+ C(\bar{q},H_{0001},H_{1010}) + C(\bar{q},H_{0010},H_{1001}) \\
    &&&+ C(H_{0001},H_{0010},H_{1100}) + B_2(q,\bar{q},K_{01},K_{10}) \\
    &&&+ C_1(q,\bar{q},H_{0010},K_{01}) + C_1(q,\bar{q},H_{0001},K_{10}) \\
    &&&+ D(q,\bar{q},H_{0001},H_{0010}) -[(g_{1011} + \bar{g}_{1011})H_{1100} \\
    &&&+ (g_{1010} + \bar{g}_{1010})H_{1101} + (g_{1001} + \bar{g}_{1001})H_{1110} ]\\
    w\beta_2^3 &: &(i\omega_0 I_n - A)H_{1003} &= A_1(q,K_{03}) + 3A_1(H_{1002},K_{01}) + 3A_1(H_{1001},K_{02}) \\
    &&&+ B(q, H_{0003}) + 3B(H_{0001}, H_{1002})+ 3B(H_{0002}, H_{1001})\\
    &&&+ 3A_2(q, K_{01}, K_{02})+ 3A_2(H_{1001},K_{01}, K_{01})  \\
    &&&+ 3B_1(q, H_{0002}, K_{01})+ 3B_1(q, H_{0001}, K_{02})\\
    &&&+ 6B_1(H_{0001}, H_{1001},K_{01})+ 3C(q, H_{0001}, H_{0002})\\
    &&&+ 3C(H_{0001}, H_{0001}, H_{1001})+ A_3(q, K_{01}, K_{01}, K_{01}) \\
    &&&+ 3B_2(q, H_{0001}, K_{01}, K_{01}) + 3C_1(q, H_{0001}, H_{0001}, K_{01})   \\
    &&&+ D(q, H_{0001}, H_{0001}, H_{0001}) \\
    &&&- [g_{1003}q + 3g_{1002}H_{1001} + 3g_{1001}H_{1002}]
\end{align*}
\subsection{Quintic terms}
Collecting quintic terms of \cref{eq:HOM_GH_ODE}:
\begin{align*}
    w^5 &: &(5i\omega_0I_n - A)H_{5000} &= 5B(q,H_{4000}) + 10B(H_{2000},H_{3000}) + 10C(q,q,H_{3000}) \\
    &&&+ 15C(q,H_{2000},H_{2000}) + 10D(q,q,q,H_{2000}) \\
    &&&+ E(q,q,q,q,q)\\
    w^3\bar{w}^2 &: &(i\omega_0I_n - A)H_{3200 } &=  3B(q,H_{2200}) + 2B(\bar{q},H_{3100})+ B(H_{0200},H_{3000}) \\
    &&&+ 6B(H_{1100},H_{2100}) + 3B(H_{1200},H_{2000}) + 3C(q,q,H_{1200})\\
    &&&+ 6C(q,\bar{q},H_{2100}) + 3C(q,H_{0200},H_{2000}) + 6C(q,H_{1100},H_{1100}) \\
    &&&+ C(\bar{q},\bar{q},H_{3000}) + 6C(\bar{q},H_{1100},H_{2000}) + D(q,q,q,H_{0200}) \\
    &&&+ 6D(q,q,\bar{q},H_{1100}) + 3D(q,\bar{q},\bar{q},H_{2000}) + E(q,q,q,\bar{q},\bar{q})\\
    &&&-[12c_2(0)q + 6(2c_1(0) + \bar{c_1}(0))H_{2100}]\\
    w^4\bar{w} &: &(3i\omega_0I_n - A)H_{4100} &=  4B(q,H_{3100}) + B(\bar{q},H_{4000}) + 4B(H_{1100},H_{3000}) \\
    &&&+ 6B(H_{2000},H_{2100}) + 6C(q,q,H_{2100}) + 4C(q,\bar{q},H_{3000}) \\
    &&&+ 12C(q,H_{1100},H_{2000}) + 3C(\bar{q},H_{2000},H_{2000}) \\
    &&&+ 4D(q,q,q,H_{1100}) + 6D(q,q,\bar{q},H_{2000}) + E(q,q,q,q,\bar{q}) \\
    &&&- 12c_1(0)H_{3000}  \\
    w^3\bar{w}\beta_2 &: &(2i\omega_0I_n - A)H_{3101} &=  A_1(H_{3100},K_{01}) + 3B(q,H_{2101}) + B(\bar{q},H_{3001}) \\
    &&&+ B(H_{0001},H_{3100}) + B(H_{0101},H_{3000}) + 3B(H_{1001},H_{2100}) \\
    &&&+ 3B(H_{1100},H_{2001}) + 3B(H_{1101},H_{2000}) + 3B_1(q,H_{2100},K_{01}) \\
    &&&+ B_1(\bar{q},H_{3000},K_{01}) + 3B_1(H_{1100},H_{2000},K_{01}) \\
    &&&+ 3C(q,q,H_{1101}) + 3C(q,\bar{q},H_{2001}) \\
    &&&+ 3C(q,H_{0001},H_{2100}) + 3C(q,H_{0101},H_{2000}) \\
    &&&+ 6C(q,H_{1001},H_{1100})+ C(\bar{q},H_{0001},H_{3000}) \\
    &&&+ 3C(\bar{q},H_{1001},H_{2000}) + 3C(H_{0001},H_{1100},H_{2000}) \\
    &&&+ 3C_1(q,q,H_{1100},K_{01}) + 3C_1(q,\bar{q},H_{2000},K_{01}) \\
    &&&+ D(q,q,q,H_{0101})+ 3D(q,q,\bar{q},H_{1001}) \\
    &&&+ 3D(q,q,H_{0001},H_{1100}) + 3D(q,\bar{q},H_{0001},H_{2000})\\
    &&&+ D_1(q,q,q,\bar{q},K_{01}) + E(q,q,q,\bar{q},H_{0001})\\
    &&&-[6g_{2101}H_{2000} + 6c_1(0)H_{2001} + (3g_{1001} + \bar{g}_{1001})H_{3100}]\\
    w^2\bar{w}\beta_1^2 &: &(i\omega_0I_n - A)H_{2120} &= 2A_1(H_{2110},K_{10}) + A_1(H_{2100},K_{20}) \\
   &&&+ 2B(q,H_{1120}) + B(\bar{q},H_{2020}) + 2B(H_{0010},H_{2110}) \\
   &&&+ B(H_{0020},H_{2100}) + 2B(H_{0110},H_{2010}) + B(H_{0120},H_{2000}) \\
   &&&+ 4B(H_{1010},H_{1110})+ 2B(H_{1020},H_{1100}) \\
   &&&+ A_2(H_{2100},K_{10},K_{10}) + 4B_1(q,H_{1110},K_{10}) \\
   &&&+ 2B_1(q,H_{1100},K_{20}) + 2B_1(\bar{q},H_{2010},K_{10}) \\
   &&&+ B_1(\bar{q},H_{2000},K_{20})  + 2B_1(H_{0010},H_{2100},K_{10}) \\
   &&&+ 2B_1(H_{0110},H_{2000},K_{10}) + 4B_1(H_{1010},H_{1100},K_{10}) \\
   &&&+ C(q,q,H_{0120}) + 2C(q,\bar{q},H_{1020}) \\
   &&&+ 4C(q,H_{0010},H_{1110}) + 2C(q,H_{0020},H_{1100}) \\
   &&&+ 4C(q,H_{0110},H_{1010}) + 2C(\bar{q},H_{0010},H_{2010}) \\
   &&&+ C(\bar{q},H_{0020},H_{2000}) + 2C(\bar{q},H_{1010},H_{1010}) \\
   &&&+ C(H_{0010},H_{0010},H_{2100}) + 2C(H_{0010},H_{0110},H_{2000}) \\
   &&&+ 4C(H_{0010},H_{1010},H_{1100}) + 2B_2(q,H_{1100},K_{10},K_{10}) \\
   &&&+ B_2(\bar{q},H_{2000},K_{10},K_{10}) + C_1(q,q,\bar{q},K_{20}) \\
   &&&+ 2C_1(q,q,H_{0110},K_{10}) + 4C_1(q,\bar{q},H_{1010},K_{10}) \\
   &&&+ 4C_1(q,K_{10},H_{0010},H_{1100}) + 2C_1(\bar{q},K_{10},H_{0010},H_{2000}) \\
   &&&+ D(q,q,\bar{q},H_{0020}) +2D(q,q,H_{0010},H_{0110}) \\
   &&&+ 4D(q,\bar{q},H_{0010},H_{1010}) + 2D(q,H_{0010},H_{0010},H_{1100}) \\
   &&&+ D(\bar{q},H_{2000},H_{0010},H_{0010}) + C_2(q,q,\bar{q},K_{10},K_{10}) \\
   &&&+ 2D_1(q,q,\bar{q},K_{10},H_{0010}) + E(q,q,\bar{q},H_{0010},H_{0010})\\
   &&& -[2g_{2120}q + 4g_{2110}H_{1010} + 2c_1(0)H_{1020} \\
   &&&+ (2g_{1020} + \bar{g}_{1020})H_{2100}   + (4g_{1010} + 2\bar{g}_{1010})H_{2110}]\\
   w^2 \bar{w}^2 \beta_2 &: &-AH_{2201} &= A_1(H_{2200},K_{01}) + 2B(q, H_{1201}) + 2B(\bar{q},H_{2101}) \\
    &&&+ B(H_{0001},H_{2200})+ 2B(H_{0101},H_{2100}) +B(H_{0200},H_{2001}) \\
    &&&+ B(H_{0201},H_{2000})+ 2B(H_{1001},H_{1200}) + 4B(H_{1100},H_{1101}) \\
    &&&+ 2B_1(q, H_{1200}, K_{01}) + B_1(H_{0200},H_{2000},K_{01}) \\
    &&&+ 2B_1(H_{1100},H_{1100},K_{01})+ 2B_1(\bar{q}, H_{2100},K_{01}) \\
    &&&+ C(q,q, H_{0201}) + 4C(q,\bar{q},H_{1101}) + 2C(q, H_{0001}, H_{1200}) \\
    &&&+ 4C(q, H_{0101}, H_{1100} ) +2C(q,H_{0200},H_{1001}) + C(\bar{q},\bar{q},H_{2001}) \\
    &&&+ 2C(\bar{q},H_{0001},H_{2100}) + 2C(\bar{q},H_{0101},H_{2000}) \\
    &&&+ 4C(\bar{q},H_{1001},H_{1100})+C( H_{0001},H_{0200},H_{2000}) \\
    &&&+ 2C(H_{0001},H_{1100},H_{1100} ) + C_1(q,q,H_{0200},K_{01}) \\
    &&&+ 4C_1(q,\bar{q},H_{1100},K_{01}) + C_1(\bar{q},\bar{q},H_{2000},K_{01})
    \\
    &&&+ 2D(q,q,\bar{q},H_{0101}) + D(q,q,H_{0001},H_{0200}) \\
    &&&+ 2D(q,\bar{q},\bar{q},H_{1001}) + 4D(q,\bar{q},H_{0001},H_{1100})\\
    &&&+ D(\bar{q},\bar{q},H_{0001},H_{2000})+ D_1(q,q,\bar{q},\bar{q},K_{01})\\
    &&&+ E(q,q,\bar{q},\bar{q},H_{0001})  - [4 (\overline{c_1}(0) + c_1(0) )H_{1101}\\
    &&&+4( \bar{g}_{2101} +  g_{2101}) H_{1100}+2 (\bar{g}_{1001} + g_{1001}) H_{2200}]\\
    w^4\beta_2 &: &(4i\omega_0 I_n - A)H_{4001} &= A_1(H_{4000},K_{01}) + 4B(q,H_{3001}) + B(H_{0001},H_{4000}) \\
    &&&+ 4B(H_{1001},H_{3000}) + 6B(H_{2000},H_{2001})+ 4B_1(q,H_{3000},K_{01}) \\
    &&&+ 3B_1(H_{2000},H_{2000},K_{01}) + 6C(q,q,H_{2001}) \\
    &&&+ 4C(q,H_{0001},H_{3000}) + 12C(q,H_{1001},H_{2000}) \\
    &&&+ 3C(H_{0001},H_{2000},H_{2000}) + 6C_1(q,q,H_{2000},K_{01}) \\
    &&&+ 4D(q,q,q,H_{1001}) + 6D(q,q,H_{0001},H_{2000}) \\
    &&&+ D_1(q,q,q,q,K_{01}) + E(q,q,q,q,H_{0001}) - 4g_{1001}H_{4000}\\
    w^2\bar{w}\beta_2^2 &: &(i\omega_0I_n - A)H_{2102} &= 2A_1(H_{2101},K_{01}) + A_1(H_{2100},K_{02}) + 2B(q,H_{1102}) \\
   &&&+ B(\bar{q},H_{2002}) + 2B(H_{0001},H_{2101}) + B(H_{0002},H_{2100}) \\
   &&&+ 2B(H_{0101},H_{2001}) + B(H_{0102},H_{2000}) + 4B(H_{1001},H_{1101})\\
   &&&+ 2B(H_{1002},H_{1100}) + A_2(K_{01},K_{01},H_{2100}) \\
   &&&+ 4B_1(q,H_{1101},K_{01}) + 2B_1(q,H_{1100},K_{02}) \\
   &&&+ 2B_1(\bar{q},H_{2001},K_{01}) + B_1(\bar{q},H_{2000},K_{02})  \\
   &&&+ 2B_1(H_{0001},H_{2100},K_{01}) + 2B_1(H_{0101},H_{2000},K_{01}) \\
   &&&+ 4B_1(H_{1001},H_{1100},K_{01}) + C(q,q,H_{0102}) \\
   &&&+ 2C(q,\bar{q},H_{1002}) + 4C(q,H_{0001},H_{1101}) + 2C(q,H_{0002},H_{1100}) \\
   &&&+ 4C(q,H_{0101},H_{1001}) + 2C(\bar{q},H_{0001},H_{2001}) \\
   &&&+ C(\bar{q},H_{0002},H_{2000}) + 2C(\bar{q},H_{1001},H_{1001}) \\
   &&&+ C(H_{0001},H_{0001},H_{2100}) + 2C(H_{0001},H_{0101},H_{2000}) \\
   &&&+ 4C(H_{0001},H_{1001},H_{1100}) + 2B_2(q,H_{1100},K_{01},H_{01}) \\
   &&&+ B_2(\bar{q},H_{2000},K_{01},K_{01}) + C_1(q,q,\bar{q},K_{02}) \\
   &&&+ 2C_1(q,q,H_{0101},K_{01}) + 4C_1(q,\bar{q},H_{1001},K_{01}) \\
   &&&+ 4C_1(q,H_{0001},H_{1100},K_{01}) + 2C_1(\bar{q},H_{0001},H_{2000},K_{01}) \\
   &&&+ D(q,q,\bar{q},H_{0002}) +2D(q,q,H_{0001},H_{0101}) \\
   &&&+ 4D(q,\bar{q},H_{0001},H_{1001}) + 2D(q,H_{0001},H_{0001},H_{1100}) \\
   &&&+ D(\bar{q},H_{2000},H_{0001},H_{0001}) + C_2(q,q,\bar{q},K_{01},K_{01}) \\
   &&&+ 2D_1(q,q,\bar{q},H_{0001},K_{01}) + E(q,q,\bar{q},H_{0001},H_{0001})\\
   &&&-[2g_{2102}q + 4g_{2101}H_{1001} + 2c_1(0)H_{1002}  \\
   &&&+ (2g_{1002} + \bar{g}_{1002})H_{2100}  + (4g_{1001} + 2\bar{g}_{1001})H_{2101}] \\
   w^3\beta_2^2 &: &(3i\omega_0I_n - A)H_{3002} &= 2A_1(H_{3001},K_{01}) + A_1(H_{3000},K_{02}) + 3B(q, H_{2002}) \\
    &&&+ 2B(H_{0001}, H_{3001}) + B(H_{0002}, H_{3000}) + 6B(H_{1001}, H_{2001}) \\
    &&&+ 3B(H_{1002}, H_{2000}) + A_2(H_{3000},K_{01}, K_{01}) \\
    &&&+ 6B_1(q, H_{2001}, K_{01}) + 3B_1(q, H_{2000},K_{02}) \\
    &&&+ 2B_1(H_{0001}, H_{3000},K_{01}) + 6B_1(H_{1001}, H_{2000},K_{01}) \\
    &&&+ 3C(q, q, H_{1002}) + 6C(q, H_{0001}, H_{2001}) \\
    &&&+ 3C(q, H_{0002}, H_{2000}) + 6C(q, H_{1001}, H_{1001}) \\
    &&&+ C(H_{0001}, H_{0001}, H_{3000}) + 6C(H_{0001}, H_{1001}, H_{2000}) \\
    &&&+ 3B_2(q, H_{2000}, K_{01}, K_{01}) + C_1(q, q, q, K_{02})\\
    &&&+ 6C_1(q, q, H_{1001}, K_{01}) + 6C_1(q, H_{0001}, H_{2000}, K_{01}) \\
    &&&+ D(q, q, q, H_{0002}) + 6D(q, q, H_{0001}, H_{1001}) \\
    &&&+ 3D(q, H_{0001}, H_{0001}, H_{2000}) + C_2(q, q, q, K_{01}, K_{01}) \\
    &&&+ 2D_1(q, q, q, H_{0001}, K_{01}) + E(q, q, q, H_{0001}, H_{0001}) \\
    &&&- [3g_{1002}H_{3000} + 6g_{1001}H_{3001}]\\
    w^2\bar{w}\beta_1\beta_2 &: &(i\omega_0I_n - A)H_{2111} &= 
   A_1(H_{2110},K_{01})+A_1(H_{2101},K_{10})+A_1(H_{2100},K_{11})\\
    &&&+2 B(q,H_{1111}) +B(\bar{q},H_{2011})+B(H_{0001},H_{2110})\\
    &&&+B(H_{0010},H_{2101})+B(H_{0011},H_{2100})+B(H_{0101},H_{2010})\\
    &&&+B(H_{0110},H_{2001})+B(H_{0111},H_{2000})+2 B(H_{1001},H_{1110})\\
    &&&+2 B(H_{1010},H_{1101})+2 B(H_{1011},H_{1100})\\
    &&&+ A_2(H_{2100},K_{01},K_{10})+2 B_1(q,H_{1110},K_{01})\\
    &&&+2 B_1(q,H_{1101},K_{10})+2B_1(q,H_{1100},K_{11})\\&&&+B_1(\bar{q},H_{2010},K_{01})+B_1(\bar{q},H_{2001},K_{10})\\
    &&&+B_1(\bar{q},H_{2000},K_{11})+B_1(H_{0010},H_{2100},K_{01})\\
    &&&+B_1(H_{0110},H_{2000},K_{01})+2B_1(H_{1010},H_{1100},K_{01})\\
    &&&+B_1(H_{0001},H_{2100},K_{10})+B_1(H_{0101},H_{2000},K_{10})\\
    &&&+2B_1(H_{1001},H_{1100},K_{10})+C(q,q,H_{0111})\\
    &&&+2C(q,\bar{q},H_{1011})+2 C(q,H_{0001},H_{1110})\\
    &&&+2 C(q,H_{0010},H_{1101})+2C(q,H_{0011},H_{1100})\\
    &&&+2 C(q,H_{0101},H_{1010})+2 C(q,H_{0110},H_{1001})\\
    &&&+C(\bar{q},H_{0001},H_{2010})+C(\bar{q},H_{0010},H_{2001})\\
    &&&+C(\bar{q},H_{0011},H_{2000})+2C(\bar{q},H_{1001},H_{1010})\\
    &&&+C(H_{0001},H^{0010},H_{2100})+C(H_{0001},H_{0110},H_{2000})\\
    &&&+2 C(H_{0001},H_{1010},H_{1100})+C(H_{0010},H_{0101},H_{2000})\\
    &&&+2C(H_{0010},H_{1001},H_{1100})+2B_2(q,H_{1100},K_{01},K_{10})\\
    &&&+B_2(\bar{q},H_{2000},K_{01},K_{10})+C_1(q,q,\bar{q},K_{11})\\
    &&&+C_1(q,q,H_{0110},K_{01})+C_1(q,q,H_{0101},K_{10})\\
    &&&+2C_1(q,\bar{q},H_{1010},K_{01})+2C_1(q,\bar{q},H_{1001},K_{10})\\
    &&&+2C_1(q,H_{0010},H_{1100},K_{01})+2C_1(q,H_{0001},H_{1100},K_{10})\\
    &&&+C_1(\bar{q},H_{0010},H_{2000},K_{01})+C_1(\bar{q},H_{0001},H_{2000},K_{10})\\
    &&&+D(q,q,\bar{q},H_{0011})+D(q,q,H_{0001},H_{0110})\\
    &&&+D(q,q,H_{0010},H_{0101})+2 D(q,\bar{q},H_{0001},H_{1010})\\
    &&&+2 D(q,\bar{q},H_{0010},H_{1001})+2 D(q,H_{0001},H_{0010},H_{1100})\\
    &&&+D(\bar{q},H_{0001},H_{0010},H_{2000})+C_2(q,q,\bar{q},K_{01},K_{10})\\
    &&&+D_1(q,q,\bar{q},H_{0010},K_{01})+D_1(q,q,\bar{q},H_{0001},K_{10})\\
    &&&+E(q,q,\bar{q},H_{0001},H_{0010})-[2g_{2111}q+2g_{2110}H_{1001}\\
    &&&+2 g_{2101} H_{1010}+2c_1(0)H_{1011} +\left(2 g_{1011} + \bar{g}_{1011}\right) H_{2100}\\
    &&&+ \left(2 g_{1010} + \bar{g}_{1010} \right) H_{2101}+\left(2 g_{1001}+\bar{g}_{1001} \right) H_{2110}]\\
    w^2\beta_2^3 &: &(2i\omega_0I_n - A)H_{2003} &= 3A_1(H_{2002},K_{01}) + 3A_1( H_{2001},K_{02}) \\
    &&&+ A_1(H_{2000},K_{03}) + 2B(q, H_{1003}) + 3B(H_{0001}, H_{2002}) \\
    &&&+ 3B(H_{0002}, H_{2001}) + B(H_{0003}, H_{2000}) + 6B(H_{1001}, H_{1002}) \\
    &&&+3A_2(H_{2001},K_{01}, K_{01}) + 3A_2(H_{2000},K_{01}, K_{02}) \\
    &&&+ B_1(q, q, K_{03}) + 6B_1(q, H_{1002}, K_{01}) \\
    &&&+ 6B_1(q, H_{1001}, K_{02}) + 6B_1(H_{0001}, H_{2001},K_{01}) \\
    &&&+ 3B_1(H_{0002}, H_{2000},K_{01}) + 6B_1(H_{1001}, H_{1001},K_{01}) \\
    &&&+3B_1(H_{0001}, H_{2000},K_{02}) + C(q, q, H_{0003}) \\
    &&&+ 6C(q, H_{0001}, H_{1002}) + 6C(q, H_{0002}, H_{1001}) \\
    &&&+ 3C(H_{0001}, H_{0001}, H_{2001}) + 3C(H_{0001}, H_{0002}, H_{2000})\\
    &&&+ 6C(H_{0001}, H_{1001}, H_{1001}) + A_3(H_{2000},K_{01}, K_{01}, K_{01}) \\
    &&&+ 3B_2(q, q, K_{01}, K_{02}) + 6B_2(q, H_{1001}, K_{01}, K_{01}) \\
    &&&+ 3B_2( H_{0001}, H_{2000},K_{01}, K_{01})+ 3C_1(q, q, H_{0002}, K_{01}) \\
    &&&+ 3C_1(q, q, H_{0001}, K_{02}) + 12C_1(q, H_{0001}, H_{1001}, K_{01}) \\
    &&&+ 3C_1(H_{0001}, H_{0001}, H_{2000},K_{01}) + 3D(q, q, H_{0001}, H_{0002})\\
    &&&+ 6D(q, H_{0001}, H_{0001}, H_{1001}) + D(H_{0001}, H_{0001}, H_{0001}, H_{2000}) \\
    &&&+ B_3(q, q, K_{01}, K_{01}, K_{01}) + 3C_2(q, q, H_{0001}, K_{01}, K_{01})\\
    &&&+3D_1(q, q, K_{01}, H_{0001}, H_{0001}) + E(q, q, H_{0001}, H_{0001}, H_{0001})\\
    &&&-[2g_{1003}H_{2000} + 6g_{1002}H_{2001} + 6g_{1001}H_{2002}] \\
    w\bar{w}\beta_2^3 &: & - AH_{1103} &= 3A_1(H_{1102},K_{01}) + 3A_1(H_{1101},K_{02}) +A_1(H_{1100},K_{03}) \\
    &&&+ B(q, H_{0103}) + B(\bar{q}, H_{1003}) + 3B(H_{0001}, H_{1102}) \\
    &&&+ 3B(H_{0002}, H_{1101})+ B(H_{0003}, H_{1100}) + 3B(H_{0101}, H_{1002}) \\
    &&&+ 3B(H_{0102}, H_{1001}) +3A_2(H_{1101},K_{01}, K_{01}) \\
    &&&+ 3A_2(H_{1100},K_{01}, K_{02}) + B_1(q, \bar{q}, K_{03}) \\
    &&&+ 3B_1(q, H_{0102}, K_{01}) + 3B_1(q, H_{0101}, K_{02}) \\
    &&&+ 3B_1(\bar{q}, H_{1002}, K_{01})+ 3B_1(\bar{q}, H_{1001}, K_{02}) \\
    &&&+ 6B_1(H_{0001}, H_{1101},K_{01})+ 3B_1( H_{0002}, H_{1100},K_{01}) \\
    &&&+ 6B_1(H_{0101}, H_{1001},K_{01}) + 3B_1(H_{0001},H_{1100},K_{02})\\
    &&&+ C(q, \bar{q}, H_{0003})+ 3C(q, H_{0001}, H_{0102}) \\
    &&&+ 3C(q, H_{0002}, H_{0101}) + 3C(\bar{q}, H_{0001}, H_{1002}) \\
    &&&+ 3C(\bar{q}, H_{0002}, H_{1001}) + 3C(H_{0001}, H_{0001}, H_{1101})\\
    &&&+ 3C(H_{0001}, H_{0002}, H_{1100}) + 6C(H_{0001}, H_{0101}, H_{1001}) \\
    &&&+ A_3(H_{1100},K_{01}, K_{01}, K_{01}) + 3B_2(q, \bar{q}, K_{01}, K_{02}) \\
    &&&+ 3B_2(q, H_{0101}, K_{01}, K_{01}) + 3B_2(\bar{q}, H_{1001}, K_{01}, K_{01}) \\
    &&&+ 3B_2(H_{0001}, H_{1100},K_{01}, K_{01}) + 3C_1(q, \bar{q}, H_{0002}, K_{01}) \\
    &&&+ 3C_1(q, \bar{q}, H_{0001}, K_{02})+ 6C_1(q, H_{0001}, H_{0101}, K_{01}) \\
    &&&+ 6C_1(\bar{q}, H_{0001}, H_{1001}, K_{01}) + 3C_1(H_{0001}, H_{0001}, H_{1100},K_{01}) \\
    &&&+ 3D(q, \bar{q}, H_{0001}, H_{0002}) + 3D(q, H_{0001}, H_{0001}, H_{0101})\\
    &&&+ 3D(\bar{q}, H_{0001}, H_{0001}, H_{1001}) + D(H_{0001}, H_{0001}, H_{0001}, H_{1100}) \\
    &&&+ B_3(q, \bar{q}, K_{01}, K_{01}, K_{01}) + 3C_2(q, \bar{q}, H_{0001}, K_{01}, K_{01}) \\
    &&&+ 3D_1(q, \bar{q}, H_{0001}, H_{0001}, K_{01}) +E(q, \bar{q}, H_{0001}, H_{0001}, H_{0001}) \\
    &&&-[(g_{1003} + \overline{g}_{1003})H_{1100} + 3(g_{1002} + \overline{g}_{1002})H_{1101} \\
    &&&+ 3(g_{1001} + \overline{g}_{1001})H_{1102}]
\end{align*}

\subsection{Sixth order terms}
Collecting sixth order terms of \cref{eq:HOM_GH_ODE}:
\begin{align*}
    w^6 &: &(6i\omega_0I_n - A)H_{6000} &= 6B(q,H_{5000}) + 15B(H_{2000},H_{4000}) \\
    &&&+ 10B(H_{3000},H_{3000}) + 15C(q,q,H_{4000}) \\
    &&&+ 60C(q,H_{2000},H_{3000}) + 15C(H_{2000},H_{2000},H_{2000})\\
    &&&+ 20D(q,q,q,H_{3000}) + 45D(q,q,H_{2000}.H_{2000}) \\
    &&&+ 15E(q,q,q,q,H_{2000}) + K(q,q,q,q,q,q) \\
    w^5\bar{w} &: &(4i\omega_0I_n - A)H_{5100} &= 5B(q,H_{4100}) + B(\bar{q},H_{5000}) + 5B(H_{1100},H_{4000}) \\
    &&&+ 10B(H_{2000},H_{3100}) + 10B(H_{2100},H_{3000}) \\
    &&&+ 10C(q,q,H_{3100}) + 5C(q,\bar{q},H_{4000}) \\
    &&&+ 20C(q,H_{1100},H_{3000}) + 30C(q,H_{2000},H_{2100}) \\
    &&&+ 10C(\bar{q},H_{2000},H_{3000}) + 15C(H_{1100},H_{2000},H_{2000}) \\
    &&&+ 10D(q,q,q,H_{2100}) + 10D(q,q,\bar{q},H_{3000}) \\
    &&&+ 30D(q,q,H_{1100},H_{2000})+ 15D(q,\bar{q},H_{2000},H_{2000}) \\
    &&&+ 5E(q,q,q,q,H_{1100}) + 10E(q,q,q,\bar{q},H_{2000}) \\
    &&&+ K(q,q,q,q,q,\bar{q}) - 20c_1(0)H_{4000} \\
    w^4\bar{w}^2 &: &(2i\omega_0I_n - A)H_{4200} &=  4B(q,H_{3200}) + 2B(\bar{q},H_{4100}) + B(H_{0200},H_{4000}) \\
    &&&+ 8B(H_{1100},H_{3100}) + 4B(H_{1200},H_{3000}) \\
    &&&+ 6B(H_{2000},H_{2200}) + 6B(H_{2100},H_{2100}) \\
    &&&+ 6C(q,q,H_{2200}) + 8C(q,\bar{q},H_{3100}) \\
    &&&+ 4C(q,H_{0200},H_{3000}) + 24C(q,H_{1100},H_{2100}) \\
    &&&+ 12C(q,H_{1200},H_{2000}) + C(\bar{q},\bar{q},H_{4000}) \\
    &&&+ 8C(\bar{q},H_{1100},H_{3000})+ 12C(\bar{q},H_{2000},H_{2100}) \\
    &&&+ 3C(H_{0200},H_{2000},H_{2000}) + 12C(H_{1100},H_{1100},H_{2000}) \\
    &&&+ 4D(q,q,q,H_{1200}) + 12D(q,q,\bar{q},H_{2100}) \\
    &&&+ 6D(q,q,H_{0200},H_{2000}) + 12D(q,q,H_{1100},H_{1100}) \\
    &&&+ 4D(q,\bar{q},\bar{q},H_{3000}) + 24D(q,\bar{q},H_{1100},H_{2000}) \\
    &&&+ 3D(\bar{q},\bar{q},H_{2000},H_{2000}) + E(q,q,q,q,H_{0200}) \\
    &&&+ 8E(q,q,q,\bar{q},H_{1100})+ 6E(q,q,\bar{q},\bar{q},H_{2000}) \\
    &&&+ K(q,q,q,q,\bar{q},\bar{q})\\
    &&&-[48c_2(0)H_{2000} + 8(3c_1(0) + \overline{c_1}(0))H_{3100}]\\
    w^3\bar{w}^3 &: &-AH_{3300} &= 3B(q,H_{2300}) + 3B(\bar{q},H_{3200})+ 3B(H_{0200},H_{3100}) \\
    &&&+ B(H_{0300},H_{3000}) + 9B(H_{1100},H_{2200}) \\
    &&&+ 9B(H_{1200},H_{2100}) + 3B(H_{1300},H_{2000}) \\
    &&&+ 3C(q,q,H_{1300}) + 9C(q,\bar{q},H_{2200}) \\
    &&&+ 9C(q,H_{0200},H_{2100}) + 3C(q,H_{0300},H_{2000}) \\
    &&&+ 18C(q,H_{1100},H_{1200}) + 3C(\bar{q},\bar{q},H_{3100}) \\
    &&&+ 3C(\bar{q},H_{0200},H_{3000}) + 18C(\bar{q},H_{1100},H_{2100}) \\
    &&&+ 9C(\bar{q},H_{1200},H_{2000}) + 9C(H_{0200},H_{1100},H_{2000}) \\
    &&&+ 6C(H_{1100},H_{1100},H_{1100}) + D(q,q,q,H_{0300}) \\
    &&&+ 9D(q,q,\bar{q},H_{1200}) + 9D(q,q,H_{0200},H_{1100}) \\
    &&&+ 9D(q,\bar{q},\bar{q},H_{2100}) + 9D(q,\bar{q},H_{0200},H_{2000}) \\
    &&&+ 18D(q,\bar{q},H_{1100},H_{1100}) + D(\bar{q},\bar{q},\bar{q},H_{3000}) \\
    &&&+ 9D(\bar{q},\bar{q},H_{1100},H_{2000}) + 3E(q,q,q,\bar{q},H_{0200}) \\
    &&&+ 9E(q,q,\bar{q},\bar{q},H_{1100}) + 3E(q,\bar{q},\bar{q},\bar{q},H_{2000}) \\
    &&&+ K(q,q,q,\bar{q},\bar{q},\bar{q}) \\
    &&&-[36(c_2(0) + \overline{c_2}(0))H_{1100} + 18(c_1(0) + \overline{c_1}(0))H_{2200} ] \\
    w^3\bar{w}^2\beta_2 &: &(i\omega_0I_n - A)H_{3201} &=  A_1(H_{3200},K_{01}) + 3B(q,H_{2201}) + 2B(\bar{q},H_{3101}) \\
    &&&+ B(H_{0001},H_{3200})+ 2B(H_{0101},H_{3100}) + B(H_{0200},H_{3001}) \\
    &&&+ B(H_{0201},H_{3000}) + 3B(H_{1001},H_{2200}) + 6B(H_{1100},H_{2101}) \\
    &&&+ 6B(H_{1101},H_{2100})+ 3B(H_{1200},H_{2001}) + 3B(H_{1201},H_{2000}) \\
    &&&+ 3B_1(q,H_{2200},K_{01}) + 2B_1(\bar{q},H_{3100},K_{01}) \\
    &&&+ B_1(H_{0200},H_{3000},K_{01}) + 6B_1(H_{1100},H_{2100},K_{01})\\
    &&&+ 3B_1(H_{1200},H_{2000},K_{01}) + 3C(q,q,H_{1201}) \\
    &&&+ 6C(q,\bar{q},H_{2101}) + 3C(q,H_{0001},H_{2200}) \\
    &&&+ 6C(q,H_{0101},H_{2100}) + 3C(q,H_{0200},H_{2001})\\
    &&&+ 3C(q,H_{0201},H_{2000}) + 6C(q,H_{1001},H_{1200}) \\
    &&&+ 12C(q,H_{1100},H_{1101}) + C(\bar{q},\bar{q},H_{3001}) \\
    &&&+ 2C(\bar{q},H_{0001},H_{3100}) + 2C(\bar{q},H_{0101},H_{3000}) \\
    &&&+ 6C(\bar{q},H_{1001},H_{2100}) + 6C(\bar{q},H_{1100},H_{2001}) \\
    &&&+ 6C(\bar{q},H_{1101},H_{2000})+ C(H_{0001},H_{0200},H_{3000}) \\
    &&&+ 6C(H_{0001},H_{1100},H_{2100}) + 3C(H_{0001},H_{1200},H_{2000})\\
    &&&+ 6C(H_{0101}, H_{1100},H_{2000}) + 3C(H_{0200},H_{1001},H_{2000}) \\
    &&&+ 6C(H_{1001},H_{1100},H_{1100}) + 3C_1(q,q,H_{1200},K_{01}) \\
    &&&+ 6C_1(q,\bar{q},H_{2100},K_{01})+ 3C_1(q,H_{0200},H_{2000},K_{01}) \\
    &&&+ 6C_1(q,H_{1100},H_{1100},K_{01}) + C_1(\bar{q},\bar{q},H_{3000},K_{01}) \\
    &&&+ 6C_1(\bar{q},H_{1100},H_{2000},K_{01}) + D(q,q,q,H_{0201}) \\
    &&&+ 6D(q,q,\bar{q},H_{1101}) + 3D(q,q,H_{0001},H_{1200}) \\
    &&&+ 6D(q,q,H_{0101},H_{1100}) + 3D(q,q,H_{0200},H_{1001})\\
    &&&+ 3D(q,\bar{q},\bar{q},H_{2001}) + 6D(q,\bar{q},H_{0001},H_{2100}) \\
    &&&+ 6D(q,\bar{q},H_{0101},H_{2000})+ 12D(q,\bar{q},H_{1001},H_{1100}) \\
    &&&+ 3D(q,H_{0001},H_{0200},H_{2000})+ 6D(q,H_{0001},H_{1100},H_{1100}) \\
    &&&+ D(\bar{q},\bar{q},H_{0001},H_{3000}) + 3D(\bar{q},\bar{q},H_{1001},H_{2000})\\
    &&&+ 6D(\bar{q},H_{0001},H_{1100},H_{2000}) + D_1(q,q,q,H_{0200},K_{01}) \\
    &&&+ 6D_1(q,q,\bar{q},H_{1100},K_{01}) + 3D_1(q,\bar{q},\bar{q},H_{2000},K_{01}) \\
    &&&+ 2E(q,q,q,\bar{q},H_{0101}) + E(q,q,q,H_{0001},H_{0200}) \\
    &&&+ 3E(q,q,\bar{q},\bar{q},H_{1001}) + 6E(q,q,\bar{q},H_{0001},H_{1100}) \\
    &&&+ 3E(q,\bar{q},\bar{q},H_{0001},H_{2000}) + E_1(q,q,q,\bar{q},\bar{q},K_{01}) \\
    &&&+ K(q,q,q,\bar{q},\bar{q},H_{0001})-[12g_{3201}q + 12c_2(0)H_{1001}  \\
    &&&+(12g_{2101} + 6\bar{g}_{2101})H_{2100} + (12c_1(0) + 6\overline{c_1}(0))H_{2101} \\
    &&&+ (3g_{1001} + 2\bar{g}_{1001})H_{3200}]\\
    w^5\beta_2 &: &(5i\omega_0 I_n - A)H_{5001} &= A_1(H_{5000},K_{01}) + 5B(q,H_{4001}) + B(H_{0001},H_{5000}) \\
    &&&+ 5B(H_{1001},H_{4000}) + 10B(H_{2000},H_{3001}) \\
    &&&+ 10B(H_{2001},H_{3000}) + 5B_1(q,H_{4000},K_{01}) \\
    &&&+ 10B_1(H_{2000},H_{3000},K_{01}) + 10C(q,q,H_{3001}) \\
    &&&+ 5C(q,H_{0001},H_{4000}) + 20C(q,H_{1001},H_{3000}) \\
    &&&+ 30C(q,H_{2000},H_{2001}) + 10C(H_{0001},H_{2000},H_{3000}) \\
    &&&+ 15C(H_{1001},H_{2000},H_{2000}) + 10C_1(q,q,H_{3000},K_{01}) \\
    &&&+ 15C_1(q,H_{2000},H_{2000},K_{01}) + 10D(q,q,q,H_{2001}) \\
    &&&+ 10D(q,q,H_{0001},H_{3000}) + 30D(q,q,H_{1001},H_{2000}) \\
    &&&+ 15D(q,H_{0001},H_{2000},H_{2000})+ 10D_1(q,q,q,H_{2000},K_{01}) \\
    &&&+ 5E(q,q,q,q,H_{1001}) + 10E(q,q,q,H_{0001},H_{2000})\\
    &&&+ E_1(q,q,q,q,q,K_{01}) + K(q,q,q,q,q,H_{0001}) \\
    &&&- 5g_{1001}H_{5000} \\
    w^4\bar{w}\beta_2 &: &(3i\omega_0I_n - A)H_{4101} &= A_1(H_{4100},K_{01}) + 4B(q,H_{3101}) + B(\bar{q},H_{4001})\\
    &&&+ B(H_{0001},H_{4100}) + B(H_{0101},H_{4000}) + 4B(H_{1001},H_{3100}) \\
    &&&+ 4B(H_{1100},H_{3001}) + 4B(H_{1101},H_{3000}) + 6B(H_{2000},H_{2101})\\
    &&&+ 6B(H_{2001},H_{2100}) + 4B_1(q,H_{3100},K_{01}) + B_1(\bar{q},H_{4000},K_{01})\\
    &&&+ 4B_1(H_{1100},H_{3000},K_{01}) + 6B_1(H_{2000},H_{2100},K_{01}) \\
    &&&+ 6C(q,q,H_{2101}) + 4C(q,\bar{q},H_{3001}) + 4C(q,H_{0001},H_{3100}) \\
    &&&+ 4C(q,H_{0101},H_{3000}) + 12C(q,H_{1001},H_{2100}) \\
    &&&+ 12C(q,H_{1100},H_{2001}) + 12C(q,H_{1101},H_{2000}) \\
    &&&+ C(\bar{q},H_{0001},H_{4000}) + 4C(\bar{q},H_{1001},H_{3000}) \\
    &&&+ 6C(\bar{q},H_{2000},H_{2001}) + 4C(H_{0001},H_{1100},H_{3000}) \\
    &&&+ 6C(H_{0001},H_{2000},H_{2100}) + 3C(H_{0101},H_{2000},H_{2000}) \\
    &&&+ 12C(H_{1001},H_{1100},H_{2000}) + 6C_1(q,q,H_{2100},K_{01}) \\
    &&&+ 4C_1(q,\bar{q},H_{3000},K_{01}) + 12C_1(q,H_{1100},H_{2000},K_{01}) \\
    &&&+ 3C_1(\bar{q},H_{2000},H_{2000},K_{01})+ 4D(q,q,q,H_{1101}) \\
    &&&+ 6D(q,q,\bar{q},H_{2001}) + 6D(q,q,H_{0001},H_{2100}) \\
    &&&+ 6D(q,q,H_{0101},H_{2000}) + 12D(q,q,H_{1001},H_{1100}) \\
    &&&+ 4D(q,\bar{q},H_{0001},H_{3000}) + 12D(q,\bar{q},H_{1001},H_{2000}) \\
    &&&+ 12D(q,H_{0001},H_{1100},H_{2000}) + 3D(\bar{q},H_{0001},H_{2000},H_{2000}) \\
    &&&+ 4D_1(q,q,q,H_{1100},K_{01}) + 6D_1(q,q,\bar{q},H_{2000},K_{01})\\
    &&&+ E(q,q,q,q,H_{0101}) + 4E(q,q,q,\bar{q},H_{1001}) \\
    &&&+ 4E(q,q,q,H_{0001},H_{1100}) + 6E(q,q,\bar{q},H_{0001},H_{2000}) \\
    &&&+ E_1(q,q,q,q,\bar{q},K_{01}) + K(q,q,q,q,\bar{q},H_{0001}) \\
    &&&-[12g_{2101}H_{3000} + 12c_1(0)H_{3001} + (4g_{1001} + \bar{g}_{1001})H_{4100}]\\
    w^2\bar{w}\beta_2^3 &: &(i\omega_0I_n - A)H_{2103} &= 3A_1(H_{2102},K_{01}) + 3A_1(H_{2101},K_{02}) \\
    &&&+A_1(H_{2100},K_{03}) + 2B(q, H_{1103}) + B(\bar{q}, H_{2003}) \\
    &&&+ 3B(H_{0001}, H_{2102}) + 3B(H_{0002}, H_{2101}) + B(H_{0003}, H_{2100}) \\
    &&&+ 3B(H_{0101}, H_{2002})+ 3B(H_{0102}, H_{2001}) + B(H_{0103}, H_{2000}) \\
    &&&+ 6B(H_{1001}, H_{1102}) + 6B(H_{1002}, H_{1101}) + 2B(H_{1003}, H_{1100}) \\
    &&&+ 3A_2(H_{2101},K_{01}, K_{01}) + 3A_2(H_{2100},K_{01}, K_{02}) \\
    &&&+ 6B_1(q, H_{1102}, K_{01}) + 6B_1(q, H_{1101}, K_{02}) \\
    &&&+ 2B_1(q, H_{1100}, K_{03}) + 3B_1(\bar{q}, H_{2002}, K_{01}) \\
    &&&+ 3B_1(\bar{q}, H_{2001}, K_{02}) + B_1(\bar{q}, H_{2000}, K_{03}) \\
    &&&+ 6B_1(H_{0001}, H_{2101},K_{01}) + 3B_1( H_{0002}, H_{2100},K_{01}) \\
    &&&+ 6B_1(H_{0101}, H_{2001},K_{01}) + 3B_1(H_{0102}, H_{2000},K_{01}) \\
    &&&+ 12B_1(H_{1001}, H_{1101},K_{01}) + 6B_1(H_{1002}, H_{1100},K_{01})\\
    &&&+ 3B_1(H_{0001}, H_{2100},K_{02}) + 3B_1(H_{0101}, H_{2000},K_{02}) \\
    &&&+ 6B_1(H_{1001}, H_{1100},K_{02}) + C(q,q,H_{0103}) + 2C(q,\bar{q},H_{1003}) \\
    &&& + 6C(q,H_{0001},H_{1102}) + 6C(q,H_{0002},H_{1101}) \\
    &&&+ 2C(q,H_{0003},H_{1100})+ 6C(q,H_{0101},H_{1002}) \\
    &&&+ 6C(q,H_{0102},H_{1001}) + 3C(\bar{q},H_{0001},H_{2002}) \\
    &&&+ 3C(\bar{q},H_{0002},H_{2001}) + C(\bar{q},H_{0003},H_{2000}) \\
    &&&+ 6C(\bar{q},H_{1001},H_{1002})+3C(H_{0001},H_{0001},H_{2101}) \\
    &&&+ 3C(H_{0001},H_{0002},H_{2100})+6C(H_{0001},H_{0101},H_{2001}) \\
    &&&+ 3C(H_{0001},H_{0102},H_{2000}) +12 C(H_{0001},H_{1001},H_{1101}) \\
    &&&+ 6 C(H_{0001},H_{1002},H_{1100}) +3C(H_{0002},H_{0101},H_{2000}) \\
    &&&+ 6C(H_{0002},H_{1001},H_{1100}) + 6C(H_{0101},H_{1001},H_{1001})\\
    &&&+ A_3(H_{2100},K_{01},K_{01},K_{01}) + 6B_2(q,H_{1101},K_{01},K_{01}) \\
    &&&+ 6B_2(q,H_{1100},K_{01},K_{02}) + 3B_2(\bar{q},H_{2001},K_{01},K_{01}) \\
    &&&+ 3B_2(\bar{q},H_{2000},K_{01},K_{02}) +3B_2(H_{0001},H_{2100},K_{01},K_{01}) \\
    &&&+ 3B_2(H_{0101},H_{2000},K_{01},K_{01})+ 6B_2(H_{1001},H_{1100},K_{01},K_{01}) \\
    &&&+ C_1(q,q,\bar{q},K_{03}) +3C_1(q,q,H_{0102},K_{01}) \\
    &&&+ 3C_1(q,q,H_{0101},K_{02}) + 6C_1(q,\bar{q},H_{1002},K_{01})\\
    &&&+ 6C_1(q,\bar{q},H_{1001},K_{02}) + 12C_1(q,H_{0001},H_{1101},K_{01}) \\
    &&&+ 6C_1(q,H_{0002},H_{1100},K_{01}) + 12C_1(q,H_{0101},H_{1001},,K_{01}) \\
    &&&+ 6C_1(q,H_{0001},H_{1100},K_{02}) + 6C_1(\bar{q},H_{0001},H_{2001},K_{01}) \\
    &&&+ 3C_1(\bar{q},H_{0002},H_{2000},K_{01})+ 6 C_1(\bar{q},H_{1001},H_{1001},K_{01}) \\
    &&&+ 3C_1(\bar{q},H_{0001},H_{2000},K_{02}) + 3C_1(H_{0001},H_{0001},H_{2100},K_{01}) \\
    &&&+ 6C_1(H_{0001},H_{0101},H_{2000},K_{01}) \\
    &&&+ 12 C_1(H_{0001},H_{1001},H_{1100},K_{01})\\
    &&&+ D(q,q,\bar{q},H_{0003}) + 3 D(q,q,H_{0001},H_{0102})\\
    &&&+ 3D(q,q,H_{0002},H_{0101}) + 6 D(q,\bar{q},H_{0001},H_{1002})\\
    &&&+ 6 D(q,\bar{q},H_{0002},H_{1001}) + 6 D(q,H_{0001},H_{0001},H_{1101}) \\
    &&&+ 6 D(q,H_{0001},H_{0002},H_{1100}) + 12 D(q,H_{0001},H_{0101},H_{1001}) \\
    &&&+ 3 D(\bar{q},H_{0001},H_{0001},H_{2001}) + 3 D(\bar{q},H_{0001},H_{0002},H_{2000}) \\
    &&&+ 6 D(\bar{q},H_{0001},H_{1001},H_{1001}) + D(H_{0001},H_{0001},H_{0001},H_{2100}) \\
    &&&+ 3D(H_{0001},H_{0001},H_{0101},H_{2000}) \\
    &&&+ 6D(H_{0001},H_{0001},H_{1001},H_{1100})\\
    &&&+2B_3(q,H_{1100},K_{01},K_{01},K_{01}) + B_3(\bar{q},H_{2000},K_{01},K_{01},K_{01}) \\
    &&&+3C_2(q,q,\bar{q},K_{01},K_{02}) + 3C_2(q,q,H_{0101},K_{01},K_{01}) \\
    &&&+6C_2(q,\bar{q},H_{1001},K_{01},K_{01}) + 6C_2(q,H_{0001},H_{1100},K_{01},K_{01}) \\
    &&&+3C_2(\bar{q},H_{0001},H_{2000},K_{01},K_{01}) +3D_1(q,q,\bar{q},H_{0002},K_{01}) \\
    &&&+ 3D_1(q,q,\bar{q},H_{0001},K_{02}) +6D_1(q,q,H_{0001},H_{0101},K_{01}) \\
    &&&+12D_1(q,\bar{q},H_{0001},H_{1001},K_{01}) \\
    &&&+6D_1(q,H_{0001},H_{0001},H_{1100},K_{01}) \\
    &&&+3D_1(\bar{q},H_{0001},H_{0001},H_{2000},K_{01}) +  3E(q,q,\bar{q},H_{0001},H_{0002}) \\
    &&&+ 3E(q,q,H_{0001},H_{0001},H_{0101}) + 6E(q,\bar{q},H_{0001},H_{0001},H_{1001}) \\
    &&&+ 2E(q,H_{0001},H_{0001},H_{0001},H_{1100}) \\
    &&&+ E(\bar{q},H_{0001},H_{0001},H_{0001},H_{2000}) \\
    &&&+C_3(q,q,\bar{q},K_{01},K_{01},K_{01}) + 3E_1(q,q,\bar{q},H_{0001},H_{0001},K_{01}) \\
    &&&+ K(q,q,\bar{q},H_{0001},H_{0001},H_{0001}) -[2g_{2103}q + 6g_{2102}H_{1001} \\
    &&&+ 6g_{2101}H_{1002} + 2c_1(0)H_{1003} + (2g_{1003} + \overline{g}_{1003})H_{2100} \\
    &&&+ (6g_{1002} + 3\overline{g}_{1002})H_{2101} + (6g_{1001} + 3\overline{g}_{1001})H_{2102}]
\end{align*}

\subsection{Seventh order terms}
Collecting seventh order terms of \cref{eq:HOM_GH_ODE}:
\begin{align*}
    w^7 &: &(7i\omega_0I_n - A)H_{7000} &= 7B(q,H_{6000}) + 21B(H_{2000},H_{5000}) + 35B(H_{3000},H_{4000}) \\
    &&&+ 21C(q,q,H_{5000}) + 105C(q,H_{2000},H_{4000}) \\
    &&&+ 70C(q,H_{3000},H_{3000}) + 105C(H_{2000},H_{2000},H_{3000}) \\
    &&&+ 35D(q,q,q,H_{4000}) + 210D(q,q,H_{2000},H_{3000}) \\
    &&&+ 105D(q,H_{2000},H_{2000},H_{2000}) + 35E(q,q,q,q,H_{3000}) \\
    &&&+ 105E(q,q,q,H_{2000},H_{2000}) + 21K(q,q,q,q,q,H_{2000}) \\
    &&&+ L(q,q,q,q,q,q,q)\\
    w^6\bar{w} &: &(5i\omega_0I_n - A)H_{6100} &= 6B(q,H_{5100}) + B(\bar{q},H_{6000}) + 6B(H_{1100},H_{5000}) \\
    &&&+ 15B(H_{2000},H_{4100}) + 15B(H_{2100},H_{4000}) + 20B(H_{3000},H_{3100}) \\
    &&&+ 15C(q,q,H_{4100}) + 6C(q,\bar{q},H_{5000}) + 30C(q,H_{1100},H_{4000})\\
    &&&+ 60C(q,H_{2000},H_{3100}) + 60C(q,H_{2100},H_{3000}) \\
    &&&+ 15C(\bar{q},H_{2000},H_{4000}) + 10C(\bar{q},H_{3000},H_{3000}) \\
    &&&+ 60C(H_{1100},H_{2000},H_{3000})+ 45C(H_{2000},H_{2000},H_{2100}) \\
    &&&+ 20D(q,q,q,H_{3100}) + 15D(q,q,\bar{q},H_{4000})\\
    &&&+ 60D(q,q,H_{1100},H_{3000}) + 90D(q,q,H_{2000},H_{2100}) \\
    &&&+ 60D(q,\bar{q},H_{2000},H_{3000}) + 90D(q,H_{1100},H_{2000},H_{2000}) \\
    &&&+ 15D(\bar{q},H_{2000},H_{2000},H_{2000}) + 15E(q,q,q,q,H_{2100}) \\
    &&&+ 20E(q,q,q,\bar{q},H_{3000}) + 60E(q,q,q,H_{1100},H_{2000}) \\
    &&&+ 45E(q,q,\bar{q},H_{2000},H_{2000}) + 6K(q,q,q,q,q,H_{1100}) \\
    &&&+ 15K(q,q,q,q,\bar{q},H_{2000}) + L(q,q,q,q,q,q,\bar{q}) - 30c_1(0)H_{5000}\\
    w^5\bar{w}^2 &: &(3i\omega_0I_n - A)H_{5200} &= 5B(q,H_{4200}) + 2B(\bar{q},H_{5100}) + B(H_{0200},H_{5000}) \\
    &&&+ 10B(H_{1100},H_{4100}) + 5B(H_{1200},H_{4000}) + 10B(H_{2000},H_{3200}) \\
    &&&+ 20B(H_{2100},H_{3100}) + 10B(H_{2200},H_{3000}) \\
    &&&+ 10C(q,q,H_{3200}) + 10C(q,\bar{q},H_{4100}) + 5C(q, H_{0200},H_{4000}) \\
    &&&+ 40C(q,H_{1100},H_{3100}) + 20C(q,H_{1200},H_{3000}) \\
    &&&+ 30C(q,H_{2000},H_{2200}) + 30C(q,H_{2100},H_{2100}) + C(\bar{q},\bar{q},H_{5000}) \\
    &&&+ 10C(\bar{q},H_{1100},H_{4000}) + 20C(\bar{q},H_{2000},H_{3100}) \\
    &&&+ 20C(\bar{q},H_{2100},H_{3000}) + 10C(H_{0200},H_{2000},H_{3000}) \\
    &&&+ 20C(H_{1100},H_{1100},H_{3000}) + 60C(H_{1100},H_{2000},H_{2100}) \\
    &&&+ 15C(H_{1200},H_{2000},H_{2000}) + 10D(q,q,q,H_{2200}) \\
    &&&+ 20D(q,q,\bar{q},H_{3100}) + 10D(q,q,H_{0200},H_{3000})\\
    &&&+ 60D(q,q,H_{1100},H_{2100}) + 30D(q,q,H_{1200},H_{2000}) \\
    &&&+ 5D(q,\bar{q},\bar{q},H_{4000}) + 40D(q,\bar{q},H_{1100},H_{3000}) \\
    &&&+ 60D(q,\bar{q},H_{2000},H_{2100}) + 15D(q,H_{0200},H_{2000},H_{2000}) \\
    &&&+ 60D(q,H_{1100},H_{1100},H_{2000}) + 10D(\bar{q},\bar{q},H_{2000},H_{3000}) \\
    &&&+ 30D(\bar{q},H_{1100},H_{2000},H_{2000}) + 5E(q,q,q,q,H_{1200}) \\
    &&&+ 20E(q,q,q,\bar{q},H_{2100}) + 10E(q,q,q,H_{0200},H_{2000})\\
    &&&+ 20E(q,q,q,H_{1100},H_{1100}) + 10E(q,q,\bar{q},\bar{q},H_{3000})\\
    &&&+ 60E(q,q,\bar{q},H_{1100},H_{2000}) + 15E(q,\bar{q},\bar{q},H_{2000},H_{2000}) \\
    &&&+ K(q,q,q,q,q,H_{0200}) + 10K(q,q,q,q,\bar{q},H_{1100}) \\
    &&&+ 10K(q,q,q,\bar{q},\bar{q},H_{2000}) + L(q,q,q,q,q,\bar{q},\bar{q})\\
    &&&-[120c_2(0)H_{3000} + (40c_1(0) + 10\overline{c_1}(0))H_{4100}]\\
    w^4\bar{w}^3 &: &(i\omega_0I_n - A)H_{4300} &= 4B(q,H_{3300}) + 3B(\bar{q},H_{4200}) + 3B(H_{0200},H_{4100}) \\
    &&&+ B(H_{0300},H_{4000})+ 12B(H_{1100},H_{3200}) + 12B(H_{1200},H_{3100}) \\
    &&&+ 4B(H_{1300},H_{3000}) + 6B(H_{2000},H_{2300}) + 18B(H_{2100},H_{2200}) \\
    &&&+ 6C(q,q,H_{2300}) + 12C(q,\bar{q},H_{3200}) + 12C(q,H_{0200},H_{3100}) \\
    &&&+ 4C(q,H_{0300},H_{3000}) + 36C(q,H_{1100},H_{2200}) \\
    &&&+ 36C(q,H_{1200},H_{2100}) + 12C(q,H_{1300},H_{2000})\\
    &&&+ 3C(\bar{q},\bar{q},H_{4100}) + 3C(\bar{q},H_{0200},H_{4000}) + 24C(\bar{q},H_{1100},H_{3100}) \\
    &&&+ 12C(\bar{q},H_{1200},H_{3000}) + 18C(\bar{q},H_{2000},H_{2200}) \\
    &&&+ 18C(\bar{q},H_{2100},H_{2100}) + 12C(H_{0200},H_{1100},H_{3000}) \\
    &&&+ 18C(H_{0200},H_{2000},H_{2100}) + 3C(H_{0300},H_{2000},H_{2000}) \\
    &&&+ 36C(H_{1100},H_{1100},H_{2100}) + 36C(H_{1100},H_{1200},H_{2000}) \\
    &&&+ 4D(q,q,q,H_{1300}) + 18D(q,q,\bar{q},H_{2200}) \\
    &&&+ 18D(q,q,H_{0200},H_{2100}) + 6D(q,q,H_{0300},H_{2000}) \\
    &&&+ 36D(q,q,H_{1100},H_{1200}) + 12D(q,\bar{q},\bar{q},H_{3100}) \\
    &&&+ 12D(q,\bar{q},H_{0200},H_{3000}) + 72D(q,\bar{q},H_{1100},H_{2100}) \\
    &&&+ 36D(q,\bar{q},H_{1200},H_{2000}) + 36D(q,H_{0200},H_{1100},H_{2000}) \\
    &&&+ 24D(q,H_{1100},H_{1100},H_{1100}) + D(\bar{q},\bar{q},\bar{q},H_{4000}) \\
    &&&+ 12D(\bar{q},\bar{q},H_{1100},H_{3000}) + 18D(\bar{q},\bar{q},H_{2000},H_{2100}) \\
    &&&+ 9D(\bar{q},H_{0200},H_{2000},H_{2000}) + 36D(\bar{q},H_{1100},H_{1100},H_{2000}) \\
    &&&+ E(q,q,q,q,H_{0300}) + 12E(q,q,q,\bar{q},H_{1200})\\
    &&&+12E(q,q,q,H_{0200},H_{1100}) + 18E(q,q,\bar{q},\bar{q},H_{2100}) +\\
    &&& 18E(q,q,\bar{q},H_{0200},H_{2000}) + 36E(q,q,\bar{q},H_{1100},H_{1100}) \\
    &&&+ 4E(q,\bar{q},\bar{q},\bar{q},H_{3000}) + 36E(q,\bar{q},\bar{q},H_{1100},H_{2000}) \\
    &&&+ 3E(\bar{q},\bar{q},\bar{q},H_{2000},H_{2000}) + 3K(q,q,q,q,\bar{q},H_{0200}) \\
    &&&+ 12K(q,q,q,\bar{q},\bar{q},H_{1100}) + 6K(q,q,\bar{q},\bar{q},\bar{q},H_{2000}) \\
    &&&+ L(q,q,q,q,\bar{q},\bar{q},\bar{q})-[144c_3(0)q + 72(2c_2(0) \\
    &&&+ \overline{c_2}(0))H_{2100} + 12(3c_1(0) + 2\overline{c_1}(0))H_{3200}]
\end{align*}

\bibliographystyle{siamplain}
\bibliography{references}

\end{document}

\fi

\ifarxiv
\clearpage
\setcounter{page}{1}
\renewcommand{\thepage}{\roman{page}}
\pdfbookmark[0]{References}{references}
\renewcommand\refname{References for main text and supplement}
\fi

\bibliographystyle{siamplain}
\bibliography{references}

\end{document}